\documentclass[12pt,english,twoside]{article}
\usepackage[T1]{fontenc}
\usepackage[latin9]{inputenc}
\usepackage[a4paper]{geometry}
\usepackage{fancyhdr}
\usepackage{verbatim}
\usepackage{amsbsy}
\usepackage{amstext}
\usepackage{setspace}
\usepackage{amssymb}
\usepackage{esint}

\makeatletter
\usepackage{graphicx}

\newtheorem{mytheo}{Theorem}
\newtheorem{mylem}{Lemma}
\newtheorem{myremark}{Remark}

\newtheorem{mydef}{Definition}

\newcommand{\myshorttitle}{L\'evy-driven synchronization models}

\renewcommand{\Re}{\,\mathrm{Re}\,}

\usepackage{cite}

\makeatother

\usepackage{babel}

\begin{document}

\newcommand{\rdate}{\mbox{\normalsize  September 9, 2014}}

\title{Intrinsic scales for high-dimensional L\' evy-driven models with
non-Markovian synchronizing updates}

\author{Anatoly Manita\thanks{{\bf Address:} Department of Probability, Faculty of Mathematics and Mechanics, Lomonosov Moscow State University, 119991, Moscow, Russia. \qquad {\bf E-mail:}  {\tt manita@mech.math.msu.su} \protect\\ \hspace*{1\parindent} This work is supported by Russian Foundation for Basic Research (grant 12-01-00897).}}

\date{\rdate }

\maketitle
\newcommand{\keyword}[1]{\relax} 

%
\begin{comment}
abccomment %
\end{comment}
{}

 \newcommand{\ds}{\displaystyle } 

\global\long\def\sP{\mathsf{P}}
 \global\long\def\sE{\mathsf{E}\,}
 \global\long\def\Lapl#1#2{\mathcal{L}(#1;#2)}
 \global\long\def\lap#1{\widetilde{#1}}
 \global\long\def\Law#1{\mathcal{P}_{#1}}

\global\long\def\RR{\mathbb{R}}
 \global\long\def\dd{1}
 \global\long\def\MADD{\boldsymbol{D}}
 \global\long\def\ESvd{2S_{v}^{2}}
 \global\long\def\tn{t_{N}}
 \global\long\def\tN{\tn}

\global\long\def\myf{\varphi}
 \global\long\def\myps{\psi}
 \global\long\def\gp{\gamma}
 \global\long\def\ap{a}
 \global\long\def\uun#1{U_{n,#1}}
 \global\long\def\mst{\star}

\global\long\def\mykon{C}
 \global\long\def\de{\delta}
 \global\long\def\kapo{\varkappa_{1}}
 \global\long\def\kapd{\varkappa_{2}}
 \global\long\def\kap{\varkappa}
 \global\long\def\GG{G}
 \global\long\def\cc{c}
 \global\long\def\ccn{\cc_{N}}
 \global\long\def\fgn{g_{N}}
 \global\long\def\emf#1{f^{(#1)}}
 \global\long\def\emV#1{V^{(#1)}}
 \global\long\def\mkn{k_{N}}
 \global\long\def\hn{\mkn}
 \global\long\def\prir{\Delta}
 \global\long\def\toh{{\textstyle \frac{1}{2}}}
 \global\long\def\ln{l_{N}}
\global\long\def\aan{a_{N}}
 \global\long\def\wcn{w_{N}(c)}

\global\long\def\dM{M}
 \global\long\def\Kn{K_{N}}
 \global\long\def\Ln{L_{N}}
 \global\long\def\elb{\ell}
 \global\long\def\mD{D_{N}}
 \global\long\def\vln{\overrightarrow{l_{N}}}
 \global\long\def\mA{A_{N}}
 \global\long\def\Id{\mbox{Id}}
 \global\long\def\vM{U}
 \global\long\def\spM{\mathcal{B}}
 \global\long\def\mC{C}
 \global\long\def\mB{B}
 \global\long\def\delem{l}
 \global\long\def\da{r}

\global\long\def\diel{b^{(N)}}
 \global\long\def\diam#1{\,\mbox{diag}\left(\left\{  #1\right\}  \right)}
 \global\long\def\betn{\beta^{(N)}}

\global\long\def\zc{z^{\circ}}
 \global\long\def\kQo{Q_{0,N}}
 \global\long\def\kQe{Q_{1,N}}
 \global\long\def\un{u_{N}}
 \global\long\def\men{\eta_{N}}
 \global\long\def\myalf{\alpha}

%%%  Above  synchro-dom-semimark-text-sept.lyx macros  

 \hyphenation{syn-chro-ni-za-tion syn-chro-ni-za-tions syn-chro-ni-zing cha-racte-ristic com-po-nents com-po-nent coor-di-na-tes con-centra-ted intro-du-ced con-si-de-ra-tion con-figura-tion pro-cesses super-posi-tion per-tur-ba-tions per-tur-ba-tion}

%%%  Local macros  

\global\long\def\prf{\varphi}
 \global\long\def\prfs{\prf_{1,s}}
 \global\long\def\prfi{\prf_{2}}
\global\long\def\hvF{\overline{F}}
 \global\long\def\mst{\ast}
 \global\long\def\poles#1{\mathcal{P}\left(#1\right)}
 \global\long\def\roots#1{\mathcal{R}\left(#1\right)}
 \global\long\def\kapn{\varkappa_{N}}
 \global\long\def\res{\mathop{\mbox{res}}}
 \global\long\def\mtet{\vartheta}
 \global\long\def\degp#1{\mbox{mult}{}_{#1}}
 \global\long\def\gamn{\gamma_{N}}
 \global\long\def\Cl#1{\mathcal{K}_{#1}}
 \global\long\def\zCl#1{\mathcal{K}_{#1}^{*}}
 \global\long\def\tCl#1#2{\mathcal{K}_{#1,#2}^{2}}
 \global\long\def\CC{\mathbb{C}}
 \global\long\def\prop{K}
 \global\long\def\kg{\gamma}
 \global\long\def\aln{\alpha_{N}}
 \global\long\def\ben{\beta_{N}}

\global\long\def\xc{x^{\circ}}
 \global\long\def\bb{\boldsymbol{a}}
 \global\long\def\chf#1#2#3{\phi^{#1}(#2;#3)}
 \global\long\def\skpr#1#2{\left\langle #1,#2\right\rangle }
 \global\long\def\etho{\boldsymbol{\eta}}
 \global\long\def\eLH{\etho^{\circ}}
 \global\long\def\zto{\boldsymbol{\zeta}}
 \global\long\def\zLH{\zto^{\circ}}
 \global\long\def\mPi#1#2{\Pi_{#1}^{(#2)}}
 \global\long\def\mPiT#1{\Pi_{#1}^{S}}
 \global\long\def\mPiTI#1{\Pi_{#1}}
 \global\long\def\tT{T}
 \global\long\def\ptT{\underline{\tT}}
 \global\long\def\olN{\overline{1,N}}
 \global\long\def\jolN{j=\olN}
 \global\long\def\nN{\mathcal{N}_{N}}
 \global\long\def\mS{S}
 \global\long\def\ggn{\theta_{1,N}}
 \global\long\def\ggo{\theta_{N}}
 \global\long\def\teton{\theta_{3,N}}
 \global\long\def\tetdn{\theta_{2,N}}
 \global\long\def\ahn{h}
 \global\long\def\myUps{\Upsilon}
 \global\long\def\myps{\boldsymbol{\psi}}
 \global\long\def\myrho{\boldsymbol{\rho}}
 \global\long\def\yy{\boldsymbol{y}}
 \global\long\def\myrst{\myrho_{\mbox{{\scriptsize st}}}}
 \global\long\def\mVar{\boldsymbol{D}}
 \global\long\def\fg{g}
 \global\long\def\sDlt{\underline{\Delta}}
 \global\long\def\nEp#1{\#^{#1}}
 \global\long\def\pTau{\underline{\tau}}
 \global\long\def\fa{a}
 \global\long\def\Var{\mbox{{Var}}}
\global\long\def\Levy{\mbox{{Lévy}}}
 \global\long\def\qn{q_{N}}
 \global\long\def\wn{w_{N}}
 \global\long\def\skb{\beta_{N}}
 \global\long\def\eqd{\stackrel{d}{=}}
 \global\long\def\ID{\mathcal{ID}}
 \global\long\def\LP{\mathcal{LP}}
\global\long\def\St#1#2{\mathcal{S}\left(#1,#2\right)}
 \global\long\def\DOA{\mbox{{\bf DOA}}}
 \global\long\def\DONA{\mbox{{\bf DONA}}}
 \global\long\def\GDOA{\mbox{{\bf GDOA}}}
 \global\long\def\mxi{\xi}
 \global\long\def\mySect{\S}
 \global\long\def\Mk{K}
 \global\long\def\mk{\kappa}
 \global\long\def\plmax#1{\left(#1\right)_{+}}

~\vspace*{-6ex} 
\begin{abstract}
We propose stochastic $N$-component synchronization models $\left(x_{1}(t),\ldots,x_{N}(t)\right)$,
$x_{j}\in\mathbb{R}^{d}$, $t\in\mathbb{R}_{+}$, whose dynamics is
described by L\' evy processes and synchronizing jumps. We prove
that symmetric models reach synchronization in a stochastic sense:
differences between components $d_{kj}^{(N)}(t)=x_{k}(t)-x_{j}(t)$
have limits in distribution as $t\rightarrow\infty$. We give conditions
of existence of natural (intrinsic) space scales for large synchronized
systems, i.e., we are looking for such sequences $\left\{ b_{N}\right\} $
that distribution of $d_{kj}^{(N)}(\infty)/b_{N}$ converges to some
limit as~$N\rightarrow\infty$. It appears that such sequence exists
if the L\' evy process enters a domain of attraction of some stable
law. For Markovian synchronization models based on $\alpha$-stable
L\' evy processes this results holds for any finite~$N$ in the
precise form with $b_{N}=(N-1)^{1/\alpha}$. For non-Markovian models
similar results hold only in the asymptotic sense. The class of limiting
laws includes the Linnik distributions. We also discuss generalizations
of these theorems to the case of non-uniform matrix-based intrinsic
scales. The central point of our proofs is a representation of characteristic
functions of $d_{kj}^{(N)}(t)$ via probability distribution of a
superposition of $N$ independent renewal processes.\medskip{}

\end{abstract}

\tableofcontents{}

\section{Introduction }

Time evolution of a multicomponent system with synchronization $x(t)=(x_{1}(t),\ldots,x_{N}(t))$,
$t\in\mathbb{R}_{+},$ consists of two parts: a free dynamics and
a spontaneous synchronizing interaction between components. $x_{j}(t)\in\mathbb{R}^{d}$
denotes the state of the component~$j$ at time~$t$. The \emph{synchronizing
interaction} is possible only at some random epochs $0<T_{1}<T_{2}<\cdots\,\,$
and has the form of instantaneous jumps $(x_{1},\ldots,x_{N})\rightarrow(x'_{1},\ldots,x'_{N})$
where the new configuration $(x'_{1},\ldots,x'_{N})$ is such that
$\left\{ x'_{1},\ldots,x'_{N}\right\} \varsubsetneq\left\{ x_{1},\ldots,x_{N}\right\} $.
The most important example is a pairwise synchronizing interaction
when for a randomly chosen pair $(j_{1},j_{2})$ the component~$j_{2}$
changes its state to the value~$x_{j_{1}}$: \begin{equation}
x_{j_{2}}(T_{n}+0)=x_{j_{1}}(T_{n}),\qquad x_{j}(T_{n}+0)=x_{j}(T_{n}),\quad j\not=j_{2}.\label{eq:pairw-int-j1-j2}\end{equation}
 The \emph{free dynamics} means that all components evolve independently
between successive epochs of interaction. 

The pairwise synchronizing interaction~(\ref{eq:pairw-int-j1-j2})
can be interpreted as follows: the com\-po\-nent~$j_{1}$ generates
a message containing information about its current state~$x_{j_{1}}$
and sends it to the component~$j_{2}$; the message reaches the destination
instantly; after receiving the message the component~$j_{2}$ reads
it and adjusts its state $x_{j_{2}}$ to the value $x_{j_{1}}$ recorded
in the message.

In this paper we consider stochastic synchronization systems which
are essentially more general than many previously studied mathematical
models~\cite{manita-TVP-large-ident-eng,Man-F-T,Man-Rome,asmda-2013-man}.
For instance, the paper~\cite{Man-F-T} studies a symmetric system
of $N$ identical Brownian particles with pairwise synch\-roni\-za\-tion.
More precisely, in~\cite{Man-F-T} the free dynamics of a single
component is the usual Wiener process with diffusion coefficient $\sigma>0$
and the sequence $\left\{ T_{n}\right\} $ is a Poisson flow of intensity~$\delta>0$.
For breavity we will refer to this system as {}``$\mathcal{B}\mathcal{M}_{N}(\sigma,\delta)$-model''.
The Markovian synchronization model of~\cite{Man-F-T} is very interesting
because many important questions relevant to its long-time behavior
can be answered in an explicit form~\cite{asmda-2013-man,Man-F-T}.
The {}``$\mathcal{B}\mathcal{M}_{N}$-model'' appears to be also
useful for constructing more sophisticated systems, for example, models
of clock synchronization in wireless sensor networks~\cite{Man-Msc-MM-2012}.
Nevertheless, the Markovian assumption is not realistic for many modern
applications. In the present paper we propose a large class non-Markovian
synchronization models. The free dynamics of components will be driven
by multi-dimen\-sional $\Levy$ processes. In particular, this assumption
permit to consider heavy tail cases. In the current paper the random
sequence $\left\{ T_{n}\right\} $ is such that, in general, the inter-event
intervals $\left\{ T_{n+1}-T_{n}\right\} _{n=1}^{\infty}$ are not
independent. Hence the sequence $\left\{ T_{n}\right\} $ is not even
a renewal process. Obviously, in this situation we cannot have any
profit from the Markov processes theory. We need to develop new specific
methods. Before discussing these methods and describing our main results
we would like to say a few words about applications that motivate
introducing the synchronizing interaction between components.

Synchronization models have their origins in computer science~\cite{Bert-Tsits}.
The key idea of asyn\-chro\-nous parallel and distributed algorithms
is to use many computing units (processors etc.) to do some common
job. Most of the time the computing units work independently but sometimes
they need to share information. The exchange of information is realized
by means of a so called message-passing mechanism~\cite{JefWitk,Bert-Tsits}.
During its work, a computing unit sends \emph{timestamped} messages
to other units. After receiving a message the computing unit analyzes
the received data and sometimes adjusts its current state to be in
agreement with other processors. Such adjustments can be interpreted
as syn\-chro\-ni\-zing jumps. Usually in these models the variable~$x_{j}$
denotes a \emph{local time} of the processor~$j$.

Similar problems arise for wireless sensor networks (WSNs)~\cite{SundBuyKshem-review,SimSpagBarStrog}.
In such networks the nodes (sensors) are almost autonomous. Each sensor
is equipped with a non-perfect noisy clock. To work with data collected
by different nodes the network needs a common notion of time. There
exist many \emph{clock synchronization protocols~}\cite{SundBuyKshem-review}
designed for wireless sensor networks. Most of them are based on the
message-passing mechanism.

The first mathematical paper on stochastic synchronization models
was~\cite{MitraMitrani}. Mitra and Mitrani studied a two-dimensional
system which corresponds to parameters $N=2$, $d=1$, $T_{k}=k$
in terms of the above general description. Multi-dimensional models
of distributed computations were proposed by many authors. Unlike~\cite{MitraMitrani}
some of their papers~\cite{AkChFuSer,GuAkFu,MadWalMes,T1,T3,nicol1991,MadWalMes1999}
were focused on numerical simulations and had only auxiliary mathematical
sections. Another papers~\cite{KuSh,ShKuRe,GrMalP,GrShSt} were devoted
to very specific parallel algorithms. The first rigorous treatment
of a multi-dimensional mathematical model with time stamp synchronization
was done in~\cite{man-Shch}. In~\cite{man-Shch,manita-TVP-large-ident-eng,Man-F-T,Man-Rome,asmda-2013-man,malysh-Manita-v}
different $N$-component synchronization systems were considered as
stochastic particle systems with special interaction. Such interpretation
is useful for invoking physical intuition. It should be noted however
that the synchronizing interaction was never studied before in the
framework of traditional interacting particle systems~\cite{Liggett-book}.
Of course, we may also describe the place of the stochastic synchronization
models in purely probabilistic terms as special perturbations of multi-dimensional
random walks.

Stochastic synchronization models with \emph{large} number of components
are of special interests. The goal is to analyze their behavior as
both the number of components $N$ and the time $t$ go to infinity.
Before formulating this general problem in precise terms it is necessary
to understand what kind of a long time behavior we can expect from
a stochastic synchronization system. The word {}``synchronization''
can be used in two senses. In a local sense we speak about \emph{synchronization}
(or equating) of some components as the results of a single synchronization
jump. In a global sense we may ask whenever the total $N$-component
system will \emph{synchronize} as $t\rightarrow\infty$ and what is
the meaning of this synchronization. Of course, this question should
be considered only for {}``irreducible'' multi-component systems
that cannot be divided into two noninteacting subsystems. It is clear
that due to the random nature of dynamics the so called \emph{perfect
synchronization} $(x_{1}=\cdots=x_{N})$ is not possible. Moreover,
as it was explained in~\cite{Man-F-T} for the {}``$\mathcal{B}\mathcal{M}_{N}$-model'',
the stochastic process $x(t)$ does not have even a \emph{limit in
law} as $t\rightarrow\infty$. Nevertheless, according to~\cite{Man-F-T}
the long time stabilization in law is expected for $x(t)$ considered
in a moving coordinate system related, for example, to a tagged particle
or to the center of mass. Note that differences $d_{ij}(x):=x_{i}-x_{j}$
are the same in both the absolute and the moving coordinate systems.
Hence all $x_{i}(t)-x_{j}(t)$ are expected to have limits in law
as $t\rightarrow\infty$. In~\cite{asmda-2013-man} this statement
was proved for the symmetric {}``$\mathcal{B}\mathcal{M}_{N}$-model''
in dimension $d=1$. Moreover, it was also proved that $\left(x_{i}(\infty)-x_{j}(\infty)\right)/\sqrt{(N-1)N}$
has a symmetric Laplace distribution which parameter does not depend
on~$N$. This means that if $t$ is large then components of $x(t)$
form a {}``collective'' which typical space size is of order $N$.
In this sense, one says that $N$ is the typical\emph{ space scale}
for the synchronized system. Note that coordinates of the center of
mass are not stochastically bounded as $t\rightarrow\infty$. It is
worth pointing out that a joint distribution of $\left(x_{i}(\infty)-x_{j}(\infty),\,1\leq i<j\leq N\right)$
cannot be found explicitly and the study of its properties for large
$N$ is a challenging problem. Another interesting problem concerning
synchronization models with \emph{large} number of components is related
to a {}``prestationary'' evolution of $x(t)$. The problem is to
find different \emph{time scales} ($t=t_{N}\rightarrow\infty$ as
$N\rightarrow\infty$) on which the synchronization system $x(t_{N})$
demonstrates completely different qualitative behaviour. The complete
description of times scales was obtained for several models~\cite{mal-man-TVP,manita-TVP-large-ident-eng,Man-Cheb190-e,Man-F-T,manita-questa}.
For example, in~\cite{Man-F-T} it was shown that the {}``$\mathcal{B}\mathcal{M}_{N}$-model''
passes three different phases before it reaches the final synchronization.
The model of clock synchronization in~WSNs (see~\cite{manita-questa})
has 5 different consecutive phases of qualitative behaviour. As it
was explained in~\cite{Man-F-T} and \cite[Sect.~5]{manita-questa},
each phase in evolution of a stochastic synchronization system is
a cumulative result of competition between two opposite tendencies:
with the course of time the free dynamics increases the {}``desynch\-roni\-za\-tion''
in the system while the interaction tries to decrease it.

In the present paper we study multi-component models $x(t)=(x_{1}(t),\ldots,x_{N}(t))$
with pairwise synchronizing interaction. These models generalize the
{}``$\mathcal{B}\mathcal{M}_{N}$-model'' of \cite{Man-F-T,asmda-2013-man}
in several directions. It is assumed that the free dynamics of components
are general $\Levy$ processes with values in $\mathbb{R}^{d}$. This
assumption makes our models very flexible. $\Levy$ processes have
independent and stationary increments. Probability distributions of
these increments may have heavy tails. Note that many modern stochastic
models in finance~\cite{Rach-Ht-finance,Shiryaev-Fin_mat}, insurance~\cite{Resnick-H-tails-ph},
data networks~~\cite{Ht-WL-syst,Khokh-Pagano-traffic,MRRS-Ntwk-traffic},
physics~\cite{Sh-Levy-Flights-ph} etc.\  use heavy-tailed $\Levy$
processes. The theory of such processes is well developped and we
will take advantage of it. We will also see that the stable $\Levy$
processes and domains of attractions of stable laws play an important
role in asymptotic analysis of synchronized system with large number
of components~$N$. Assumptions about the sequence $\boldsymbol{T}=\left\{ T_{n}\right\} $
of synchronization epochs are very natural in the context of multi-component
systems. It is assumed that each component $j$ generates messages
at epochs of some renewal process $\boldsymbol{\tau}^{(j)}=\left\{ \tau_{m}^{(j)}\right\} $
independently of other components. Hence the point process $\left\{ T_{n}\right\} $
is the superposition of $N$ renewal processes: $\boldsymbol{T}=\cup_{j}\boldsymbol{\tau}^{(j)}$.
In general, the superposition of renewal processes no longer forms
a renewal process and therefore an analysis of $\boldsymbol{T}=\left\{ T_{n}\right\} $
is a difficult task. There exists a huge number of studies in this
field~\cite{Grig-Lit,Gned-Kovalenko,Cinlar-Agnew,Cinlar-Sup,Zakusilo-Melesh,Lam-Lehoczky},
most of them are devoted to limit theorems. Unfortunately, none of
them is applicable to our situation. 

The paper is organized as follows. In $\mySect$~\ref{sub:syncr-ep-T-M}--\ref{sub:gen-synchro-model}
we introduce a general synchronization model. The precise definition
of the pairwise interaction is given in $\mySect$~\ref{sub:syncr-ep-T-M}
in terms of parameters $F_{k}(s)$, $k=\olN$, and $R=\left(r_{kj}\right)_{k,j=1}^{N}$
where $F_{k}(s):=\sP\left(\tau_{q+1}^{(k)}-\tau_{q}^{(k)}\leq s\right)$,
$s\in\RR$, is the c.d.f.\  of inter-event intervals in the flow
$\boldsymbol{\tau}^{(k)}$ and $R$ is a routing matrix used for choosing
message destinations. To introduce free dynamics we recall some classic
results from the $\Levy$ processes theory. The free dynamics is determined
by a set of $\Levy$ exponents $\left\{ \eLH_{j}(\lambda),\,\, j=1,\ldots,N\right\} $,
$\lambda\in\RR^{d}$, (see~$\mySect$~\ref{sub:free-evol-nonsym}
). We show that such approach includes, as examples of free dynamics,
Brownian motions, random walks in~$\RR^{d}$ and, in particular,
random walks with heavy-tailed jumps. The general $N$-component synchronization
model with the above parameters will be denoted by $\mathcal{GG}_{N}\left(\left\{ \eLH_{j}\right\} _{j=1}^{N};\left\{ F_{j}\right\} _{j=1}^{N},R\right)$.
While the free dynamics, the flows $\boldsymbol{\tau}^{(k)}$, the
random routing and the initial configuration~$x(0)$ are assumed
to be independent the stochastic process $\left(x(t),\, t\in\RR_{+}\right)$
is very complicated and, in general, \emph{non-Markovian}. The only
exception is the situation when all c.d.f.\  $F_{j}(s)$ correspond
to exponential distributions: $F_{j}(s)=\plmax{1-\exp\left(-s/m_{j}\right)}$,
$m_{k}>0$. Under such assumption the point process $\boldsymbol{T}=\left\{ T_{n}\right\} $
is a Poisson flow and the process $\left(x(t),\, t\in\RR_{+}\right)$
is Markovian. In this case we will use notation $\mathcal{GM}_{N}\left(\left\{ \eLH_{j}\right\} _{j=1}^{N};\left\{ F_{j}\right\} _{j=1}^{N},R\right)$.

The paper is focused on symmetric synchronization models whose definition
is given in~$\mySect$~\ref{sec:sym-synch-mod-N}. The symmetry
assumption means that evolutions of all components follows the same
probabilistic rules with the same parameters ($\eLH_{j}=\eLH$, $F_{j}=F$
$\forall j$) and the routing~$R$ is uniform. We will use short
notation $\mathcal{GGS}_{N}(\eLH;F)$ for the general symmetric model
and $\mathcal{GMS}_{N}(\eLH;m)$ for the Markovian symmetric model
with c.d.f.\  $F(s)=\plmax{1-\exp\left(-s/m\right)}$. For the general
model $\mathcal{GGS}_{N}(\eLH;F)$ we assume that distribution of
the inter-arrival intervals in the flows $\boldsymbol{\tau}^{(j)}$
has a \emph{rational Laplace transform} (see RPFN class in~$\mySect$~\ref{sub:inter-event-symmetr}).
This class of distribution was discussed in 1955 by Cox~\cite{Cox-1955-Cmplx-pr}.
It is large enough to cover a variety of applications in queueing
theory. These distributions are very convenient for analytical treatment
and numerical simulations. Moreover, any probability distribution
on $\RR_{+}$ can be approximated arbitrarily close (in terms of weak
convergence) by distributions with rational Laplace transforms. As
it was shown in~\cite{Asmus-Bladt} this class of probability laws
coincides with the \emph{ME} (matrix-exponential) \emph{distributions}.
It contains as proper subsets the phase-type distributions~\cite{Neuts-book,O-Cinneide-1999-Ph-open},
the Coxian distributions~\cite{KPW-Loss-models-futher}, the general
Erlangian distributions~\cite{Cox,Lipsky} etc. We believe that most
of results of our paper remain true for more general class of distributions
but such generalization would make some of our proofs much longer. 

In $\mySect$~\ref{sub:Lim-distrib} we assume that $N$ is fixed
and $t\rightarrow+\infty$. Under general assumptions on the free
dynamics of the symmetric model $\mathcal{GGS}_{N}(\eLH;F)$ in Theorem~\ref{t-N-fixed-limit-in-t}
we prove existence of the limit in law for the differences $x_{k}(t)-x_{j}(t)$,
\[
d_{kj}^{(N)}(t):=x_{k}(t)-x_{j}(t)\,\stackrel{d}{\longrightarrow}\, d_{kj}^{(N)}(\infty).\]
Next step is to study the distribution of $d_{kj}^{(N)}(\infty)$
for large values of $N$. This problem has different answers for Markovian
and non-Markovian cases. Theorem~\ref{t:chi-N-inf-lambda} is devoted
to characteristic function $\chi_{N}(\infty;\lambda)$ of the limiting
law. Under general assumptions on $\left\{ \tT_{n}\right\} $ it gives
the following asymptotic representation of $\chi_{N}(\infty;\lambda)$
for large~$N$: \[
\chi_{N}(\infty;\lambda):=\sE\exp\, i\skpr{\lambda}{d_{kj}^{(N)}(\infty)}=\frac{1}{1+\ggn\etho(\lambda)}+\tetdn(\lambda).\]
 Here $\etho(\lambda)=-2\Re\eLH(\lambda)$, the real sequence $\left\{ \ggn\right\} $
is such that $\ggn\sim\frac{1}{2}mN$ as $N\rightarrow\infty$, \[
m:=\sE\left(\tau_{q+1}^{(j)}-\tau_{q}^{(j)}\right)=\int_{0}^{\infty}y\, dF_{j}(y),\]
and the sequence of functions $\left\{ \tetdn(\lambda)\right\} $
vanishes uniformly in $\lambda\in\RR^{d}$. This representation is
of great importance for subsequent sections. 

It appears (Theorem~\ref{t-Markov-chi-expl} in $\mySect$~\ref{sub:Markov-case})
that for the Markovian symmetrical model $\mathcal{GMS}_{N}(\eLH;m)$
we have $\ggn=\frac{1}{2}(N-1)m$ and $\tetdn(\lambda)\equiv0$. This
implies (Theorem~\ref{t-Markov-stable} in $\mySect$~\ref{sub:Markov-case})
that if the free dynamics of the model is driven by an $\alpha$-stable
$\Levy$ process then the probability distribution of $d_{jk}^{(N)}(\infty)/\left(N-1\right)^{1/\alpha}$
has the characteristic function \[
\frac{1}{1+\frac{1}{2}m\etho(\lambda)}\]
and hence it does not depend on~$N$. In this case we may say that
the sychronized system possesses an intrinsic space scale $\left(N-1\right)^{1/\alpha}\sim N^{1/\alpha}$.
Indeed, since typical distances between components of the synchronized
system are of order~$N^{1/\alpha}$ it is natural to consider this
system on a new space scale with a new unit which is equal to $N^{1/\alpha}$
old units.

For non-Markovian models~$\mathcal{GGS}_{N}(\eLH;F)$ the function
$\tetdn(\lambda)$ is necessarily nonzero (see~$\mySect$~\ref{sub:chi-f-decomp}).
Therefore we cannot expect such nice result on the existence of the
intrinsic scale for any fixed~$N$ as in the Markovian case. Nevertheless
similar results hold in the asymptotic sense (when $N\rightarrow\infty$)
if we make additional assumptions about the free dynamics. For asymptotic
results it is not strictly necessary to assume that the free dynamics
is a stable $\Levy$ process. It is sufficiently to take the free
dynamics from the domain of attraction of some stable law~in $\RR^{d}$
($\mySect$~\ref{sub:Intr-scales}). The theory of attraction to
stable laws is classical and well developped~\cite{Gned-Kolm-1968-eng,Rvacheva,Meer-Sikor}.
Theorem~\ref{t:asympt-char-f} states that if the free dynamics belongs
to the domain of attraction of some stable law and $\left\{ b_{n}\right\} $
is a corresponding normalizing sequence then the distribution of $d_{jk}^{(N)}(\infty)/b_{N}$
weakly converges as $N\rightarrow\infty$ to some distribution~$Q_{\infty,\infty}(dx)$
on~$\RR^{d}$. A situation when the attracting stable law has the
index of stability~$\alpha$ and the normalizing sequence is $b_{n}=n^{1/\alpha}$
is known as the normal attraction~\cite{Gned-Kolm-1968-eng}. Hence
the distribution of $d_{jk}^{(N)}(\infty)/b_{N}$ is asymptotically
not depending on~$N$. So we may say that $b_{N}$ is the \emph{instrinsic
space scale} of a \emph{large} synchronized $N$-component system.
It is interesting to note that the limit distribution $Q_{\infty,\infty}(dx)$
belongs to the class of symmetric geometric stable distributions~\cite{Mittnik-Rach}
(see Remark~\ref{rem-Geom-stable-d} in~$\mySect$~\ref{sub:Linnik-distrib}).
In particular, this class contains the Laplace distribution and the
famous Linnik distribution~\cite{Linnik53}.

In~$\mySect$~\ref{sub:matrix-Jurek} we generalize Theorem~\ref{t:asympt-char-f}
to the case of matrix-based scales when intrinsic space transformations
have the form of linear operators $d_{jk}^{(N)}(\infty)\mapsto B_{N}d_{jk}^{(N)}(\infty)$
for some special $d\times d$ matrices~$B_{N}$. We show that existence
of such intrinsic matrix scales is related to the problem of attraction
to \emph{operator stable} laws in~$\RR^{d}$ \cite{Meer-Sikor,KMPS-Op-Geo-stable-laws}.
In the case $B_{N}=N^{-B}$ these \emph{non-uniform} scales can be
described in terms of Jurek coordinates~\cite{Jurek-Mason,Meer-Sikor}.
Hence in dimensions~$d>1$ the class of $N$-component synchronization
systems discussed in~$\mySect$~\ref{sub:matrix-Jurek} is much
wider than the class of models of~$\mySect$~\ref{sub:Intr-scales}
with {}``scalar'' intrinsic scales. 

Section~\ref{sec:Proofs-step-by-step} contains proofs of all theorems.
These proof use the representation of the characteristic function
$\chi_{N}(t;\lambda):=\sE\exp\, i\skpr{\lambda}{d_{kj}^{(N)}(t)}$
in terms of the $\Levy$ exponent $\etho(\lambda)$ and generation
functions related to the superposition~$\boldsymbol{T}=\cup_{j}\boldsymbol{\tau}^{(j)}$
of the renewal processes~$\boldsymbol{\tau}^{(j)}$ (Lemmas~\ref{l-repr-chi}
and~\ref{l-IN-t-lambda} in~$\mySect$~\ref{sub:repr-ch-f}). To
get this representation we need a chain of auxiliary results on the
free dynamics and the interaction (Lemmas~\ref{lem:V0-Sx-nm}--\ref{l:rec-eq-Vn-k}
in~$\mySect$~\ref{sub:Lemmas-dynam}--\ref{sub:rec-equats}). These
lemmas are similar to their analogues proved for Markovian models~in~\cite{Man-F-T,manita-TVP-large-ident-eng}.
Nevertheless, the proof of Lemma~\ref{l:rec-eq-Vn-k} meets additional
difficulties related to the involved nature of the sequence of synchronization
epochs~$\boldsymbol{T}$. The symmetry assumption is very essential
for the proof of Lemma~\ref{l:rec-eq-Vn-k}. Note that Lemma~\ref{lem:V0-Sx-nm}
can be generalized for symmetric synchronizing {}``multi-particle''
interactions (see~\cite{manita-TVP-large-ident-eng} and $\mySect$~\ref{sub:perturb-indep})
which are more general than the pairwise interactions. This possibility
opens the way to an obvious generalization of the present paper.

The representation for the characteristic function $\chi_{N}(t;\lambda)$
provided by Lemmas~\ref{l-repr-chi} and~\ref{l-IN-t-lambda} gives
an explicit formula for Markovian models $\mathcal{GMS}_{N}(\eLH;m)$
($\mySect$~\ref{sub:the-Mark-case}). Therefore Theorems~\ref{t-Markov-chi-expl}
and~\ref{t-Markov-stable} (including the convergence in~Theorem~\ref{t-N-fixed-limit-in-t})
easily follow from that explicit formula. 

The non-Markovian case $\mathcal{GGS}_{N}(\eLH;F)$ is more complicated.
Even the existence of the limit $\ds\lim_{t\rightarrow\infty}\chi_{N}(t;\lambda)$
in~Theorem~\ref{t-N-fixed-limit-in-t} is not evident. At first
look a special adaptation of the classic Key Renewal Theorem (KRT)
might be helpful for calculating such limits. But as it is explained
in~$\mySect$~\ref{sub:around-KRT} it is very unlikely that the
classical sufficient conditions for the KRT could be effectively checked
in our concrete problem. So we restrict ourself to renewal processes~$\boldsymbol{\tau}^{(j)}$
with the ME distribution of inter-event intervals. Keeping in mind
this assumption, in $\mySect$~\ref{sub:class-Kt} we develop some
simple rules for manipulating expressions arising in Lemmas~\ref{l-repr-chi}
and~\ref{l-IN-t-lambda}. These rules permit us to get a short proof
of~Theorem~\ref{t-N-fixed-limit-in-t} in an {}``algebraic manner''.
Theorem~\ref{t-N-fixed-limit-in-t} follows from Lemmas~\ref{l-IN-infty-1}
and~\ref{l-IN-infty-2} which proofs are given in~$\mySect$~\ref{sub:class-Kt}.
Lemma~\ref{l-IN-infty-2} also provides an integral representation
for the limiting characteristic function~$\chi_{N}(\infty;\lambda)$.
This representation will be useful for proving Theorems~\ref{t:chi-N-inf-lambda}
and~\ref{t:asympt-char-f} in~\ref{sub:chi-f-decomp}. The method
of these proofs is based on using the Laplace transform for generating
functions. It reduces to an analysis of singularities of rational
complex functions. Such approach is standard in the context of the
classical renewal theory~\cite{Cox,Gned-Kovalenko}. But it is necessary
to pay attention to coefficients in decomposions (Lemma~\ref{l-c0-d0-bounds})
because they depend on~$N$. The problem is to find singularities
giving the principal asymptotics~(Lemma~\ref{l-fi-fi2-princ-exp})
and to obtain precise bounds for the coefficients. $\mySect$~\ref{sub:chi-f-large-N}
completes proofs of Theorems~\ref{t:chi-N-inf-lambda} and~\ref{t:asympt-char-f}.

\section{Model. Definitions. Assumptions. Notation}

In $\mySect$~\ref{sub:perturb-indep} for explanatory purposes only
we describe a general approach to constructing a large class of stochastic
synchronization models. We try to show that different existing synchronization
models may be considered within the unified framework of special perturbations
of simple stochastic evolutions. A definition of our model and precise
assumptions are given in~$\mySect$~\ref{sub:syncr-ep-T-M}--\ref{sub:gen-synchro-model}.

\subsection{Perturbation of independent dynamics by synchronization}

\label{sub:perturb-indep}

Imagine there is some system consisting of $N$ components which are
labeled by the set $\nN=\left\{ 1,\ldots,N\right\} $. First we introduce
independent dynamics of the components.

Let $\left(\xc_{1}(t),\, t\in\RR_{+}\right),\,\ldots,\left(\xc_{N}(t),\, t\in\RR_{+}\right)$
be independent stochastic processes taking their values in $\RR^{d}$.
Assume that each process $\left(\xc_{j}(t),\, t\in\RR_{+}\right)$
has independent increments. We interpret the variable $\xc_{j}\in\RR^{d}$
as a state of the component~$j$ and the set of processes $\xc(t)=\left(\xc_{1}(t),\ldots,\xc_{N}(t)\right)$
as a free dynamics of the system. 

Next we add a perturbation to the system. We modify the evolution
$\xc$ by introducing a special interaction between components. This
interaction happens at random times and consists in a partial synchro\-niza\-tion
of component states. 

For any map $M:\,\nN\rightarrow\nN$ define $\nu_{M}=\mbox{card}\, M\nN$
which is the number of different elements in the image $M\nN=\left\{ M(j):\, j\in\nN\right\} .$
Consider also a set of fixed points $U_{M}=\left\{ j:\, M(j)=j\right\} $.
The map $M$ is called a \emph{synchronization map} if $\nu_{M}<N$
and $\mbox{card}\, U_{M}=\nu_{M}$. Denote by $\mathcal{M}_{N}$ a
set of all synchronization maps of the set $\nN$.

Let $\left\{ \tT_{n},\, n\in\mathbb{Z}_{+}\right\} $ be a random
sequence\[
0\equiv\tT_{0}<\tT_{1}<\cdots<\tT_{n}<\cdots\]
 and $\left\{ M_{n},\, n\in\mathbb{N}\right\} $ be a sequence of
$\mathcal{M}_{N}$-valued random variables. We do not assume that
$\left\{ \tT_{n}\right\} $ and $\left\{ M_{n}\right\} $ are independent.
Consider a new stochastic process $x(t)=\left(x_{1}(t),\ldots,x_{N}(t)\right)\in\left(\RR^{d}\right)^{N}$
which paths are determined by the following relations \begin{equation}
x(t)-x(0)=\xc(t)-\xc(0),\quad t\in[0,\tT_{1}],\qquad x(\tT_{n}+0)=\left(x\circ M_{n}\right)(\tT_{n})\,,\quad n\geq1,\label{eq:pert-1}\end{equation}
\begin{equation}
x(t)=x(\tT_{n}+0)+\left(\xc(t)-\xc(\tT_{n}+0)\right)\,\qquad t\in(\tT_{n},\tT_{n+1}],\label{eq:pert-2}\end{equation}
where $y=\left(x\circ M\right)$ is the vector $y=(y_{1},\ldots,y_{N})\in\left(\RR^{d}\right)^{N}$
with coordinates $y_{j}=x_{M(j)}$, $j\in\nN$. The correspondence
$x\mapsto y=\left(x\circ M\right)$ between points of the configuration
space~$\left(\RR^{d}\right)^{N}$ will be called a \emph{synchronization
jump}. In some sense the process $x(t)$ is the special perturbation
of the free dynamics $\xc(t)$. We will call the process $x(t)$ a
\emph{stochastic synchronization system}. 

We always assume that \emph{initial configuration} $x(0)$ is \emph{independent}
of $\xc(\cdot)$, $\left\{ \tT_{n}\right\} $ and $\left\{ M_{n}\right\} $.

Sometimes another terminology is useful. We can speak about interacting
particle systems (instead of multi-component systems) and consider
$x_{j}(t)$ as a coordinate of $j$-th particle. In \cite{manita-TVP-large-ident-eng}
we studied a system of $N$ identical particles moving as independent
random walks (free dynamics $\xc(t)$) and interacting by means of
special $m$-particle synchronizations happened at epochs $\left\{ \tT_{n}\right\} $
of some Poisson flow. In that case all synchronizing maps $M_{n}$
satisfy the condition $\nu_{M}=N-m+l$ for some $l\leq m/2$. Multiparticle
synchronizations ($m>2$) will not be considered further in this paper.
Starting from $\mySect$~\ref{sub:syncr-ep-T-M} we consider only
pairwise interactions.

Hence the $N$-component stochastic synchronization system $x(t)=\left(x_{1}(t),\ldots,x_{N}(t)\right)$
is determined by specifying the following ingredients:
\begin{description}
\item [{(F)}] the free dynamics $\xc(t)=\left(\xc_{1}(t),\ldots,\xc_{N}(t)\right)$ 
\item [{(T)$+$(M)}] the random flow of synchronization epochs $\left\{ \tT_{n}\right\} $
and the sequence of synchronization maps $\left\{ M_{n}\right\} $
\item [{(I)}] the initial distribution of $x(0)=\left(x_{1}(0),\ldots,x_{N}(0)\right)$ 
\end{description}
The above assumptions on (F), (T)$+$(M) and (I) need to be precised
when defining a concrete model. In some models it is convenient to
consider a marked point process~\cite{bremaud,Daley-Vere-Jones}
\[
\left(\tT_{1},\mk_{1}\right),\ldots,\left(\tT_{n},\mk_{n}\right),\ldots\]
with a finite set of marks $\Mk$ and a marked sequence of synchronization
maps $\left\{ M_{n}^{(\mk_{n})}\right\} $. The interaction (T)$+$(M)
is build by a two-stage construction: first, the generation of the
sequence $\left\{ \left(\tT_{n},\mk_{n}\right)\right\} _{n=1}^{\infty}$,
and then the generation of conditionally independent maps $\left\{ M_{n}^{(\mk_{n})}\right\} $.
For nonsymmetric models probability distributions~of \textbf{$M_{n}^{(\mk)}$
}may be different for different marks~$\mk$. Such situation will
be considered in the current paper, see $\mySect$~\ref{sub:syncr-ep-T-M}
for details.

The above ingredients (F) and (T)$+$(M) may be correlated. For example,
papers~\cite{man-Shch} and~\cite{Malyshkin} were devoted to particular
models in which the probability distribution of $M_{n}$ depends on
$x(\tT_{n})$. 

In models studied in the present paper the free dynamics $\xc(t)$
and the couple ($\left\{ \tT_{n}\right\} $,$\left\{ M_{n}\right\} $)
are independent.

\subsection{Assumptions on synchronization epochs and synchronization maps}

\label{sub:syncr-ep-T-M}

In this paper we consider a pairwise synchronization which is based
on the well known message-passing mechanism \cite{JefWitk,Bert-Tsits}.
This means that components of the system can share the data with other
components by sending and receiving messages containing information
about a current state of the sender. Below we will use terminology
of particle systems and speak about particles instead of components.

Each particle~$k$ has its own sequence of times \[
0<\tau_{1}^{(k)}<\tau_{2}^{(k)}<\cdots\,\]
when it sends messages to other particles. For convenience we put
$\tau_{0}^{(k)}\equiv0$. The choice of recipients will be discussed
below. Denote $\Delta_{n}^{(k)}=\tau_{n}^{(k)}-\tau_{n-1}^{(k)}\,$. 

Let the random variables $\left(\Delta_{n}^{(k)},\, n\in\mathbb{N}\right)$
be independent and identically distributed. This means that $\Pi_{t}^{(k)}=\max\left\{ n:\,\tau_{n}^{(k)}\leq t\right\} $,
$t\geq0$, is a simple renewal process. Assume that for any $k$ a
c.d.f.\  $F_{k}(s)=\sP\left\{ \Delta_{n}^{(k)}\leq s\right\} $ is
continuous. We assume also that the renewal processes $\left(\Pi_{t}^{(k)},\, t\geq0\right)$,
$k=1,\ldots,N$, are independent. Consider events $C_{k_{1},k_{2}}=\left\{ \exists n,m:\,\tau_{n}^{(k_{1})}=\tau_{m}^{(k_{2})}\right\} $.
It follows that \[
\sP\left(\bigcup_{1\leq k_{1}<k_{2}\leq N}C_{k_{1},k_{2}}\right)=0.\]
Consider a point process \[
0=\tT_{0}<\tT_{1}<\tT_{2}<\cdots\,\]
generated by the superposition of the renewal processes $\Pi_{t}^{(k)},$
$k=1,\ldots,N$. In general, inter-arrival times $\tT_{q}-\tT_{q-1}$
are not independent. Denote $\mPiT t=\sum_{j=1}^{N}\mPi tj$. In other
words, $\mPiT t=\max\left\{ m:\,\,\tT_{m}\leq t\,\right\} $.

Fix some $N\times N$ matrix $R=\left(r_{ij}\right)_{i,j=1}^{N}$,
$r_{ii}=0,$ $r_{ij}\geq0$, $\sum_{j=1}^{N}r_{ij}=1$. We define
the interaction between particles of $x(t)$ by means of synchronization
jumps which occur at times of the point process $\left\{ \tT_{q}\right\} $.
Namely, for any point $\tT_{q}$ there exists a unique (random) pair
$(j_{1},n)$, $j_{1}\in\left\{ 1,\ldots,N\right\} ,$ $n\in\mathbb{N}$,
such that $\tT_{q}=\tau_{n}^{(j_{1})}$. It means that at time $\tT_{q}$
the particle~$j_{1}$ sends a message to some another particle $j_{2}$
which is chosen independently with probability $r_{j_{1}j_{2}}$.
The message contains information on the current value of~$x_{j_{1}}$.
Messages reach their destinations instantly. After receiving the message
from~$j_{1}$ the particle~$j_{2}$ ajusts its coordinate to the
value~$x_{j_{1}}$: $x_{j_{2}}(\tT_{q}+0)=x_{j_{1}}(\tT_{q})$. This
is the only jump in the system at the time~$\tT_{q}$: $x_{j}(\tT_{q}+0)=x_{j}(\tT_{q})$
for all $j\not=j_{2}.$ Define a map $\mS_{j_{1}j_{2}}:\,\nN\rightarrow\nN$
as follows \[
\mS_{j_{1}j_{2}}(j)=\left\{ \begin{array}{cc}
j, & j\not=j_{2},\\
j_{1}, & j=j_{2}\,.\end{array}\right.\]
We see that if $\tT_{q}=\tau_{n}^{(j_{1})}$ then the random synchronization
map $M_{q}$ is such that \[
\sP\left\{ M_{q}=\mS_{j_{1}j_{2}}\right\} =r_{j_{1}j_{2}}\]
 for all $j_{2}\not=j_{1}$. In particular, $\nu_{M_{q}}=N-1$.

Hence the synchronization is determined by the following parameters:
$F_{k}(s)$, $k=1,\ldots,N$, and the matrix $R=\left(r_{ij}\right)_{i,j=1}^{N}$.

As it was mentioned in Subsection~\ref{sub:perturb-indep} between
receiving of subsequent messages the particles evolve according to
the free dynamics.

Note that the above defined random sequences $\left\{ \tT_{q}\right\} $
and $\left\{ M_{q}\right\} $ correspond to the formal scheme of $\mySect$~\ref{sub:perturb-indep}.
Namely, $\left\{ \tT_{q}\right\} $ can be obtained from a marked
point process $\left\{ \left(\tT_{q},\mk_{q}\right)\right\} $ where
the set of marks $\Mk$ is $\left\{ 1,\ldots,N\right\} $ and $\mk_{q}$
is such that $\tT_{q}=\tau_{n}^{(\mk_{q})}$ for some~$n$. The probability
distribution of $M_{q}^{(\mk_{q})}$ depends on the mark~$\mk_{q}$
because the values $\mS_{\mk_{q}j}$, $\jolN$, are taken with probabilities~$r_{\mk_{q}j}$.

\subsection{Free evolution }

\label{sub:free-evol-nonsym}

Assume that $\xc_{1}(t),\ldots,\xc_{N}(t)$ are independent $\Levy$
processes. This means that here we make an assumption stronger than
the independence of increments condition (see Subsection~\ref{sub:perturb-indep}).
The $\Levy$ processes theory is well developped (see, for example,
\cite{Sato-Levy,Applebaum-Levy-book}) and we want to make use of
it. We recall basic definitions and introduce some notation.

\begin{mydef}

A stochastic process $\left(\xc_{j}(t),\, t\in\RR_{+}\right)$ is
called a $\Levy$ process if 
\begin{itemize}
\item it starts from the origin: $\xc_{j}(0)=0\in\RR^{d}$ 
\item it has independent and stationary increments
\item it is stochastically continuous.
\end{itemize}
\end{mydef}

Let $y_{1}$ and $y_{2}$ be two vectors in~$\RR^{d}$, $y_{m}=(y_{m}^{1},\ldots,y_{m}^{d})$,
$m=1,2$. Denote by $\skpr{y_{1}}{y_{2}}$ their scalar product, i.e.,
$\skpr{y_{1}}{y_{2}}=\sum_{l=1}^{d}y_{1}^{l}y_{2}^{l}$. If $Y$ is
a random vector in $\RR^{d}$ then $\psi_{Y}(\lambda)$ denotes its
characteristic function: \[
\psi_{Y}(\lambda)=\exp\left(i\skpr{\lambda}Y\right),\quad\lambda\in\RR^{d}.\]
The random vector $Y$ is said to be \emph{infinitely divisible} if
for all $n\in\mathbb{N}$ there exist i.i.d. random vectors $Z_{1}^{(1)},\ldots,Z_{n}^{(n)}$
such that\[
Y\eqd Z_{1}^{(1)}+\cdots+Z_{n}^{(n)}.\]
As usual the notation $V_{1}\eqd V_{2}$ means that random vectors
$V_{1}$ and $V_{2}$ have the same distribution. The fundamental
result established by $\Levy$ and Khinchine states that $\psi_{Y}(\lambda)=\exp\rho_{Y}(\lambda)$
where the function $\rho_{Y}:\,\RR^{d}\rightarrow\CC$ can be represented
in a special form known as the $\Levy$-Khinchine formula~\cite{Applebaum-Levy-book,Gned-Kolm-1968-eng,Sato-Levy}.
We will not use here this formula explicitely. When we need to say
that $Y$ has an infinitely divisible distribution with the $\Levy$
exponent $\rho_{Y}(\lambda)$ we will simply write $Y\sim\ID\left(\rho_{Y}(\lambda)\right)$.

It it clear that increments of a $\Levy$ process are infinitely divisible.
In the sequel we will use the following classical result~\cite{Applebaum-Levy-book}.
Let $\chf j{t-s}{\lambda}$ the characteristic function of the increment
$\xc_{j}(t)-\xc_{j}(s)$~: \[
\chf j{t-s}{\lambda}=\sE\exp\left(i\skpr{\lambda}{\xc_{j}(t)-\xc_{j}(s)}\right),\qquad0\leq s\leq t.\]
Then $\chf jt{\lambda}=e^{t\eLH_{j}(\lambda)}\,,$ $t\geq0$, with
some function $\eLH_{j}:\,\RR^{d}\rightarrow\CC$ having the $\Levy$-Khinchine
form. For such $\Levy$ process $\left(\xc_{j}(t),\, t\geq0\right)$
we will use  a short notation $\xc_{j}\sim\LP\left(\eLH_{j}\right)$. 

We see that the set of $\Levy$ exponents $\left\{ \eLH_{j}(\lambda),\,\, j=1,\ldots,N\right\} $
completely determines free dynamics of our model.

\paragraph{Examples of the free dynamics driven by $\Levy$ processes.}
\begin{itemize}
\item Each component $\xc_{j}(t)$ is a $d$-dimensional \emph{Brownian
motion} with a constant drift: \[
d\,\xc_{j}(t)=\sigma_{j}dB_{j}(t)+b_{j}dt,\]
where $\sigma_{j}$ is a real $d\times d$ matrix, $b_{j}\in\RR^{d}$,
$B_{j}(t)=(B_{j}^{1}(t),\ldots,B_{j}^{d}(t))$ and $B_{j}^{1},\ldots,B_{j}^{d}$
are independent standard Wiener processes with values in $\RR^{1}$.
This case corresponds to the function \begin{equation}
\eLH_{j}(\lambda)=i\skpr{b_{j}}{\lambda}-\frac{1}{2}\skpr{\sigma_{j}\sigma_{j}^{T}\lambda}{\lambda}\,.\label{eq:Win-xo}\end{equation}

\item \emph{Random walks in} $\RR^{d}$. The component $\xc_{j}(t)$ is
a continuous time jump Markov process with generator \begin{equation}
\left(L_{j}f\right)(y)=\beta_{j}\int_{\RR^{d}}\left(f(y+q)-f(y)\right)\mu_{j}(dq),\qquad f\in C_{b}(\RR^{d},\RR),\label{eq:Lj-f-y}\end{equation}
where $C_{b}(\RR^{d},\RR)$ is the Banach space of bounded continuous
functions $f:\,\RR^{d}\rightarrow\RR$, $\beta_{j}>0$ is the intensity
of jumps and a probability measure $\mu_{j}$ is the distribution
of jumps. It is easy to see that in this case \[
\eLH_{j}(\lambda)=\beta_{j}\int_{\RR^{d}}\left(e^{i\skpr{\lambda}q}-1\right)\mu_{j}(dq)\,.\]

\item \emph{Random walks in} $\mathbb{Z}^{d}$. This is a subcase of~(\ref{eq:Lj-f-y})
with the measures $\mu_{j}(dq)$ supported in $\mathbb{Z}^{d}$: \[
S_{j}=\mbox{supp}\,\mu_{j}\subset\mathbb{Z}^{d}.\]
Then \[
\left(L_{j}f\right)(y)=\beta_{j}\sum_{q\in S_{j}}\left(f(y+q)-f(y)\right)\mu_{j}(\left\{ q\right\} ),\quad\quad f\in C_{b}(\mathbb{Z}^{d},\RR),\]
and \[
\eLH_{j}(\lambda)=\beta_{j}\sum_{q\in S_{j}}\left(e^{i\skpr{\lambda}q}-1\right)\mu_{j}(\left\{ q\right\} )\,.\]

\item Consider the following particular subcase of (\ref{eq:Lj-f-y}) \begin{equation}
\mu_{j}(dq)=\frac{C_{\bb}\mathbf{1}_{\left\{ \left|q\right|\geq1\right\} }\, dq}{\left|q\right|^{d+\bb}},\qquad q\in\RR^{d},\label{eq:mu-j-q-d-al}\end{equation}
where the papameter $\bb$ is positive and $C_{\bb}$ is a normalizing
factor. Evidently, the distribution~(\ref{eq:mu-j-q-d-al}) has a
finite expectation iff $\bb>1$. Moreover, it has a finite variance
iff $\bb>2$. 
\end{itemize}

\subsection{General synchronization models}

\label{sub:gen-synchro-model}

Assume that $\xc_{1}(t),\ldots,\xc_{N}(t)$ satisfy to assumptions
of $\mySect$~\ref{sub:free-evol-nonsym} and ($\left\{ \tT_{n}\right\} $,$\left\{ M_{n}\right\} $)
satisfy to assumptions of $\mySect$~\ref{sub:syncr-ep-T-M}. Assume
also that the free dynamics $x^{\circ}$, the pair ($\left\{ \tT_{n}\right\} $,$\left\{ M_{n}\right\} $)
and an initial configuration $x(0)$ are \emph{independent}. A stochastic
process $x(t)=(x_{1}(t),\ldots,x_{N}(t))\in\left(\RR^{d}\right)^{N}$
defined by~(\ref{eq:pert-1})--(\ref{eq:pert-2}) will be called
an $N$-\emph{component synchronization system}. To specify parameters
of the model we will use notation $\mathcal{GG}_{N}\left(\left\{ \eLH_{j}\right\} _{j=1}^{N};\left\{ F_{j}\right\} _{j=1}^{N},R\right)$.

We list some simple properties of the\emph{ general} model $\mathcal{GG}_{N}\left(\left\{ \eLH_{j}\right\} _{j=1}^{N};\left\{ F_{j}\right\} _{j=1}^{N},R\right)$:
\begin{itemize}
\item Under assumptions of $\mySect$~\ref{sub:syncr-ep-T-M}--\ref{sub:free-evol-nonsym}
the process $x(t)$ is stochastically continuous\-.
\item $x(t)$ is not a process with independent increments. 
\item The process $x(t)$ is neither Markovian nor semi-Markovian.
\end{itemize}
The lack of Markovian property is explained by the complicated structure
of the sequence~$\left\{ \tT_{n}\right\} $. However there is an
important exclusion.

\begin{myremark}  \label{rem:Markov-semi-M}

If all $F_{k}(s)$ correspond to exponential distributions, \[
F_{k}(s)=\plmax{1-\exp\left(-s/m_{k}\right)},\quad s\in\RR,\quad m_{k}>0,\]
then $x(t)$ is a Markov process. Indeed, in this case the point process
$\left\{ \tT_{n}\right\} $ is a Poissonian flow as the superposition
of independent Poissonian flows $\left\{ \tau_{n}^{(j)}\right\} $,
$\jolN$.

\end{myremark} 

Sometimes we will denote the Markovian model by $\mathcal{GM}_{N}\left(\left\{ \eLH_{j}\right\} _{j=1}^{N};\left\{ F_{j}\right\} _{j=1}^{N},R\right)$.

\section{Symmetric models: main results}

In this paper we mainly study a \emph{symmetric} synchronization model
which will be introduced in Subsection~\ref{sec:sym-synch-mod-N}

\subsection{Symmetry assumptions}

\label{sec:sym-synch-mod-N}

The general synchronization system was introduced in Subsections~\ref{sub:syncr-ep-T-M}--\ref{sub:gen-synchro-model}.
Here we add more assumptions to define symmetric model.

\smallskip 

\textbf{Free dynamics. }We assume that all functions $\eLH_{j}$,
$\jolN$, defining the independent $\Levy$ processes $\xc_{j}(t)$
are equal: $\eLH_{j}(\lambda)\equiv\eLH(\lambda)$.

\smallskip 

\textbf{Synchronization epochs.} $F_{j}(y)=F(y)$ for all $\jolN$.

\smallskip 

\textbf{Routing matrix. }Senders choose destinations for their messages
uniformly: $r_{jk}=1/(N-1)$ for all $k\not=j$, $r_{jj}=0$.

In other words, the symmetric model means that all components are
identical. Their evolutions follow the same probabilistic rules with
the same parameters.

\smallskip 

For any random vector $z=(z_{1},\ldots,z_{N})\in\left(\RR^{d}\right)^{N}$
with components $z_{j}\in\RR^{d}$ we denote by $\Law z$ the distribution
law of $z$. Hence $\Law z$ is some probability measure on $\left(\RR^{d}\right)^{N}$.
Let $\pi:\,(1,\ldots,N)\rightarrow(i_{1},\ldots,i_{N})$ be an arbitrary
permutation. The permutation~$\pi$ generates a map on $\left(\RR^{d}\right)^{N}$:
\[
\pi\star(z_{1},\ldots,z_{N})=(z_{i_{1}},\ldots,z_{i_{N}})\,.\]

\textbf{Initial distribution. }Assume that the initial distribution
$\Law{x(0)}$ is invariant with respect to permutations of indices,
i.e., \begin{equation}
\Law{\pi\star x(0)}=\Law{x(0)}\label{eq:P-pi-x0}\end{equation}
for all~$\pi$. Note that the denenerated case when all components
start from the origin, i.e., \[
x_{i}(0)=0\quad\mbox{for all}\quad i=1,\ldots,N,\]
is a particular example of the assumption~(\ref{eq:P-pi-x0}).

As it was already mentioned in Subsection~\ref{sub:perturb-indep}
we always assume that\emph{ }initial configuration $x(0)$ is \emph{independent}
of the free dynamics and synchronizations.

If all above assumptions hold then for all $t>0$ the distribution
of the $N$-component system $x(t)$ remains invariant with respect
to permutations of indices. In such case we will simply call the process
$x(t)$ a symmetric synchronization model. 

Note that, in general, the symmetric model $x(t)$ is not Markovian
nor semi-Markovian stochastic process. The only exception is the situation
discussed in Remark~\ref{rem:Markov-semi-M}. 

The \emph{general} (non-Markovian) \emph{symmetric} model will be
denoted by $\mathcal{GGS}_{N}(\eLH;F)$. For \emph{Markovian} \emph{symmetric}
model we use notation $\mathcal{GMS}_{N}(\eLH;m)$ where $m>0$ is
the mean  of the exponential distribution with c.d.f.\  $F(s)=\plmax{1-\exp\left(-s/m\right)}$.

\subsection{Desynchronization between components}

Since the free dynamics of different components are independent our
stochastic system will never reach the \emph{perfect synchronization
regime} when states of all components $x_{1}(t)$, $\ldots$, $x_{N}(t)$
become equal after some (possibly random) time $t_{0}$. Such phenomenon
is impossible due to the stochastic nature of the dynamics. What we
can expect\emph{ }is a long time stabilization of synchronization
errors in the distributional sense. To get some control over magnitudes
of the synchronization errors we will consider differences $d_{jk}^{(N)}(t)=x_{j}(t)-x_{k}(t)$
between states of any pair ($j$,$k$) at time $t$.

Let a probability measure $P_{N,t}(dx)$ on $\left(\RR^{d},\mathcal{B}(\RR^{d}))\right)$
be the distribution of $d_{jk}^{(N)}(t)=x_{j}(t)-x_{k}(t)$. By the
symmetry assumptions it is the same for all $j\not=k$. Similarly,
the characteristic function of $d_{jk}^{(N)}(t)$, \begin{equation}
\chi_{N,j,k}(t;\lambda)=\sE\exp\left(i\skpr{\lambda}{x_{j}(t)-x_{k}(t)}\right),\quad\quad\lambda\in\RR^{d},\label{eq:chi-N-jk}\end{equation}
does not depend on $j$ and $k$ for the symmetric model. So we will
omit indices $j,k$ and use notation $\chi_{N}(t;\lambda)$, \[
\chi_{N}(t;\lambda)=\int_{\RR^{d}}e^{i\skpr{\lambda}x}P_{N,t}(dx)\,.\]
Our aim is to study the characteristic function $\chi_{N}(t;\lambda)$
for large $t$ and $N$. Main results will be presented in Subsections~\ref{sub:Lim-distrib}
and~\ref{sub:Intr-scales}.

\subsection{Assumptions on inter-event interval distribution}

\label{sub:inter-event-symmetr}

The independent renewal processes $\left(\Pi_{t}^{(k)},\, t\geq0\right)$,
$k=1,\ldots,N$, defined in Subsection~\ref{sub:syncr-ep-T-M} are
identically distributed in the symmetric model. Up to the end of this
paper we will assume that the below conditions holds.

\smallskip \par\noindent  \textbf{Assumption P1.} The probability
distribution function $F$ is absolutely continuous: \begin{eqnarray*}
F(s) & = & \sP\left\{ \Delta_{n}^{(k)}\leq s\right\} =\int_{0}^{s}p(s')\, ds',\quad\quad s\geq0,\\
F(s) & = & 0,\quad\quad p(s)\,=\,0,\quad\quad s<0.\end{eqnarray*}

\smallskip 

Note that this assumption concerns only inter-event intervals in each
$\Pi_{t}^{(k)}$. The point process $\left\{ \tT_{n}\right\} $ which
is the superposition of $\Pi_{t}^{(k)}$, $k=1,\ldots,N$, is very
complicated.

Given a function $q=q(s)$ such that $q(s)=0$ for $s<0$, we denote
by $q^{*}(z)$ its Laplace tranform~\cite{Cox},\[
q^{*}(z)=\int_{0}^{+\infty}e^{-zs}q(s)\, ds,\qquad z\in\CC.\]
If $q(s)$ is a probability density function then $q^{*}(z)$ is well
defined at least in the complex half-plane $\left\{ z:\,\Re z\geq0\right\} $.

Before introducing the next assumption we discuss a special class
of complex functions $f=f(z)$, $f:\,\mathbb{C}\rightarrow\mathbb{C}$.
We say that $f(z)$ is a \emph{RPF-function} if it can be represented
as a proper fraction $f(z)=\frac{P(z)}{Q(z)}$ where $P(z)$ and $Q(z)$
are some polynomials such that $\mbox{deg}\, P<\mbox{deg}\, Q$. Note
that summation and mutiplication of RPF-functions again give a RPF-function.
Evidently, any RPF-function has finite number of poles and is vanishing
as $z\rightarrow\infty$. Such functions can be written as \begin{equation}
f(z)=\sum_{j=1}^{v}\sum_{k=1}^{n_{j}}\left(z-z_{j}\right)^{-k}\, c_{j,k}\,,\label{eq:Laurent-RPF}\end{equation}
where $n_{j}\geq1$ are natural numbers, $z_{j}\in\mathbb{C}$ are
poles of $f$ and $c_{j,k}\in\mathbb{C}$. The representation~(\ref{eq:Laurent-RPF})
is just the sum of principal parts of Laurent expansions about poles,
the number $n_{j}$ is the order of the pole $z_{j}$. 

If all poles $z_{i}$ have \emph{strictly negative} real parts ($\Re z_{i}<0$)
we say that the function~$f$ belongs to \emph{the class RPFN}. 

\smallskip \par\noindent  \textbf{Assumption P2.} The probability
density function $p(s)$ is such that its Laplace transform $p^{*}(z)$
is a RPFN-function.

\smallskip  

As it was already mentioned in Introduction the probability distributions
satisfying to the Assumption~P2 are exactly the ME distributions~\cite{Asmus-Bladt}.
An important role of distributions with rational Laplace tranform
for the queueing theory was discovered by Cox in~\cite{Cox-1955-Cmplx-pr}.

In particular, Assumption P2 implies the existence of an exponential
moment

\[
\sE\exp\left(\delta\Delta_{n}^{(k)}\right)=\int_{0}^{\infty}\exp(\delta u)\, p(u)\, du<\infty\]
for some $\delta>0$ and hence the existence of all moments \begin{equation}
m_{r}=\sE\left(\Delta_{n}^{(k)}\right)^{r}=\int_{0}^{+\infty}s^{r}p(s)\, ds\,,\quad\quad r\in\mathbb{N}.\label{eq:moments-interv}\end{equation}
For shortness we will use also notation $m$ for the mean: $m=m_{1}=\ds\int s\, p(s)\, ds$.

The function $p(s)$ is a probability density hence $p^{*}(0)=1$.
If $p(s)$ satisfies Assumption~P2 then the equation \begin{equation}
1-p^{*}(z)=0\label{eq:1-pz-z-0}\end{equation}
 has a finite number of roots. Let $\left\{ r_{0},\, r_{1},\ldots,r_{q}\right\} $,
$r_{0}=0$, be the set of different roots of the equation~(\ref{eq:1-pz-z-0}).

\begin{mylem} \label{l--roots-1-p-0} All numbers $r_{1},\ldots,r_{q}$
belong to the subplane $\Re z<0$.

\end{mylem}

\begin{flushleft}
\emph{Proof}. Since $\left(p^{*}\right)'(0)=-m<0$ the root $r_{0}=0$
is simple. Note that $\left|p^{*}(v)\right|<1$ for any $v\in\left\{ z\in\CC:\,\Re z\geq0\right\} \backslash\left\{ 0\right\} $.
Indeed, \[
\left|p^{*}(a+ib)\right|\leq p^{*}(a)\qquad a,b\in\RR.\]
Evidently, $p^{*}(a)<1$ for $a>0$ so $\left|p^{*}(z)\right|<1$
if $\Re z>0$. Moreover, $\left|p^{*}(ib)\right|<1$, $b\not=0$,
as a characteristic function of a non-lattice distribution~\cite{SHIRYAEV-Probab-eng}.
\hfill  $\square$
\par\end{flushleft}

\smallskip \par\noindent  \textbf{Assumption P3.} The roots $r_{1},\ldots,r_{q}$
are \textbf{simple} that is $\left(p^{*}\right)'(r_{j})\not=0$. 

Assumption P3 is not necessary for the main results but it makes some
proofs shorter. Obviously, Assumption P3 corresponds to the general
case situation.

\subsection{Limiting distributions}

\label{sub:Lim-distrib}

We consider the symmetric synchronization model of Subsection~\ref{sec:sym-synch-mod-N}
under Assumptions~P1 and P2. Recall that $P_{N,t}(dx)$ and $\chi_{N}(t;\lambda)$
denote the distribution law and the characteristic function of $d_{jk}^{(N)}(t)=x_{j}(t)-x_{k}(t)$.

\begin{mytheo}  \label{t-N-fixed-limit-in-t}

For any fixed $N$ the distribution of $d_{jk}^{(N)}(t)=x_{j}(t)-x_{k}(t)$
has a (weak) limit as $t\rightarrow\infty$: \[
P_{N,t}\stackrel{w}{\rightarrow}P_{N,\infty}\,.\]

\end{mytheo} 

This theorem follows from the $\Levy$ continuity theorem (Theorem~3.6.2
in~\cite{Lukacs}) and the next lemma. 

\begin{mylem}  \label{l-ch-f-limit}

For any fixed $N$ the family of characteristic functions $\left\{ \chi_{N}(t;\lambda),\, t\geq0\right\} $
con\-ver\-gen\-ces to some function $\chi_{N}(+\infty;\lambda)$
as $t\rightarrow+\infty$ and, moreover, this convergence is uniform
in~$\lambda\in\RR^{d}$.

\end{mylem}  It is well known \cite[Th. 3.6.2]{Lukacs} that the
function $\chi_{N}(+\infty;\lambda)$ is the characteristic function
of the limiting distribution $P_{N,\infty}$.

In the next theorem we need additional Assumptions~P3.

\begin{mytheo}  \label{t:chi-N-inf-lambda}Let Assumptions~P1-P3
hold. The characteristic function $\chi_{N}(+\infty;\lambda)$ admits
the following representation\[
\chi_{N}(+\infty;\lambda)=\frac{1}{1+\ggn\etho(\lambda)}+\tetdn(\lambda).\]
 Here $\etho(\lambda)=-2\Re\eLH(\lambda)$, the real sequence $\left\{ \ggn\right\} $
is such that $\ggn\sim\frac{1}{2}mN$ as $N\rightarrow\infty$ and
the sequence of functions $\left\{ \tetdn(\lambda)\right\} $ vanishes
uniformly in $\lambda$: \begin{equation}
\sup_{\lambda\in\RR^{d}}\left|\tetdn(\lambda)\right|\rightarrow0\quad\quad(N\rightarrow\infty)\,.\label{eq:theta-2-N}\end{equation}

\end{mytheo} 

Theorem~\ref{t:chi-N-inf-lambda} is proved in Subsection~\ref{sub:chi-f-decomp}.
We will see from Subsection~\ref{sub:chi-f-large-N} that for the
sequence of functions $\left\{ \tetdn(\lambda)\right\} $ a result
stronger than~(\ref{eq:theta-2-N}) holds. Namely, there exists a
real sequence $\left\{ \ggo\right\} $ such that $\ggo\rightarrow0$
as $N\rightarrow\infty$, \begin{equation}
\sup_{\lambda\in\RR^{d}}\left|\tetdn(\lambda)\right|\leq\ggo,\label{eq:theta-uniforn-in-eta}\end{equation}
 and $\left\{ \ggo\right\} $ is the same for any function $\etho=\etho(\lambda)\geq0$.

\subsection{Intrinsic scales for synchronized $N$-component systems}

\label{sub:Intr-scales}

Distributions of the differences $d_{jk}^{(N)}=x_{j}-x_{k}$ are important
from practical and theoretical viewpoints because many reasonable
synchronization error estimates are functions of $d_{jk}^{(N)}$.
When we consider the symmetric $N$-component system for large $N$
we may ask about a proper space scale which depends on $N$ and corresponds
to typical values of the synchronization errors. It appears that probabilistic
properties of the free dynamics have an important impact on the typical
scale of the synchronized system.

\subsubsection*{Stable random vectors}

We need to remind some classical facts about stable distributions~\cite{Sam-Taqqu,Uchaj-Zolot}. 

\begin{mydef} \label{def-UUnUD}

A random vector $U\in\RR^{d}$ has a stable distribution if there
exist an $\alpha\in(0,2]$ and a sequence $\left\{ D_{n}\right\} $
of nonrandom vectors in $\RR^{d}$ such that for any $n\in\mathbb{N}$
\begin{equation}
U_{1}+\cdots+U_{n}\eqd n^{1/\alpha}U+D_{n}\label{eq:UUnUD}\end{equation}
where $U_{1},\ldots,U_{n}$ are independent copies of $U$.

\end{mydef}

\begin{mydef}

The vector $U\in\RR^{d}$ is called strictly stable if (\ref{eq:UUnUD})
holds with $D_{n}=0$. 

\end{mydef}

Recall that the probability distribution of a random vector $V$ is
called symmetric if $V\eqd-V$. A symmetric stable vector is strictly
stable.

The stable laws are infinitely divisible~\cite{Sato-Levy,Applebaum-Levy-book}.
Hence the characteristic function $\myps_{U}(\lambda)=\sE\exp\left(i\skpr{\lambda}U\right)$
of a stable vector~$U$ has the form $\myps_{U}(\lambda)=\exp\zto_{U}(\lambda)$.
Therefore the distribution of the stable vector $U$ is completely
determined by the function $\zto_{U}(\lambda)$. We will denote the
stable distribution defined in~(\ref{eq:UUnUD})  by $\St{\alpha}{\zto_{U}(\lambda)}$
and write $U\sim\St{\alpha}{\zto_{U}(\lambda)}$. Note that the parameter~$\alpha$
is also determined by $\zto_{U}(\lambda)$. The presence of~$\alpha$
in $\St{\alpha}{\zto_{U}(\lambda)}$ is not necessary but it makes
the notation more informative. The number $\alpha$ is called the
index of stability. It is evident that $(-U)\sim\St{\alpha}{\overline{\zto_{U}(\lambda)}}$
where $\overline{z}$ denotes the complex conjugation of~$z$. 

The general form of the function $\zto_{U}(\lambda)$ is known~\cite{Sam-Taqqu,Uchaj-Zolot}
but we will not use it. We simply note that~(\ref{eq:UUnUD}) can
be rewritten as \begin{equation}
\exp\left(n\zto_{U}(\lambda)\right)=\exp\left(i(\lambda,D_{n})+\zto_{U}(n^{1/\alpha}\lambda)\right),\qquad\lambda\in\RR^{d}.\label{eq:psi-stable}\end{equation}
In the case $\alpha=2$ the stable laws are exaclty the $d$-dimensional
Gaussian distributions.

Stable laws are the only possible limiting distributions of scalar-normalized
sums of i.i.d. random vectors. The following definition is equivalent
to Definition~\ref{def-UUnUD} (see~\cite{Sam-Taqqu}).

\begin{mydef} \label{def-DOA-stable} 

A random vector $U\in\RR^{d}$ is stable if it has a domain of attraction,
i.e., if there is a random vector $V$ and sequences of positive numbers
$\left\{ b_{n}\right\} $ and nonrandom vectors $\left\{ C_{n}\right\} $,
$C_{n}\in\RR^{d}$, such that \begin{equation}
\frac{V_{1}+\cdots+V_{n}}{b_{n}}+C_{n}\stackrel{d}{\longrightarrow}U\label{eq:VVb-CU}\end{equation}
where $V_{1},\ldots V_{n},\ldots$ are independent copies of~$V$
and the notation $\stackrel{d}{\longrightarrow}$ denotes convergence
in distribution. 

\end{mydef}

In the situation of Definition~\ref{def-DOA-stable} the random vector~$V$
is said to be in the \emph{domain of attraction} of the stable vector~$U$.
Following the book~\cite{Meer-Sikor} we will write $V\in\DOA(U)$.
In the case when the normalizing sequence $\left\{ b_{n}\right\} $
has the form $b_{n}=n^{1/\alpha}$ we say that $V$ belongs to the
\emph{domain of normal attraction} of~$U$ and write $V\in\DONA(U)$.
Sometimes we will put in these notation distributions instead of random
vectors. Evidently, $\DOA(U)\supset\DONA(U)\ni U$.

The exhaustive study of domains of attraction for one-dimensional
stable laws were presented in~\cite{Gned-Kolm-1968-eng}.  In dimensions
$d\geq2$ the first results about domains of attraction belong to
Rvacheva~\cite{Rvacheva}, the disciple of B.V.~Gnedenko.

We will need the next simple facts following directly from (\ref{eq:VVb-CU})
and Definition~\ref{def-DOA-stable}.

\begin{mylem} \label{l-V-V-i-iii}

Let $V'$ be an independent copy of some random vector $V$. Let a
random vector $U$ be stable with the index~$\alpha$.
\begin{description}
\item [{\phantom{ii}(i)}] If $V\in\DOA(U)$ then $V-V'\in\DOA(\underline{U})$
where $\underline{U}\sim\St{\alpha}{2\Re\zto_{U}(\lambda)}$. Moreover,
the normalyzing sequence $\left\{ b_{n}\right\} $ in~(\ref{eq:VVb-CU})
is the same for $V$ and $V-V'$.
\item [{\phantom{i}(ii)}] The statement (i) remains true if we replace
$\DOA$ by $\DONA$. 
\item [{(iii)}] Assume aditionally that $V$ is infinitely divisible: $V\in\ID(\rho(\lambda))$.
Then $V-V'$ is infinitely divisible too: $V-V'\in\ID(2\Re\rho(\lambda))$
\end{description}
\end{mylem}

\subsubsection*{Infinite divisible laws in the domains of attraction}

Let $\yy(t)\in\RR^{d}$, $t\geq0$, be a $\Levy$ process with the
characteristic function of increments $\chf{\yy}t{\lambda}=e^{t\myrho(\lambda)}$,
$\lambda\in\RR^{d}$, i.e., $\yy\sim\LP(\myrho(\lambda))$ in notation
of Subsection~\ref{sub:free-evol-nonsym}. Let $\St{\alpha}{\myrst(\lambda)}$
be some stable distribution with the index of stability $\alpha$,
$0<\alpha\leq2$. 

\begin{mydef}

We say that the $\Levy$ process $\yy=(\yy(t),\, t\geq0)$\emph{ }\textbf{belongs
to the domain of attraction} of the stable law $\St{\alpha}{\myrst(\lambda)}$
if \[
\yy(1)\in\DOA(\St{\alpha}{\myrst(\lambda)}).\]
 We say that $\yy=(\yy(t),\, t\geq0)$\emph{ }belongs to the domain
of\emph{ }\textbf{normal attraction} of the stable law $\St{\alpha}{\myrst(\lambda)}$
if $\yy(1)\in\DONA(\St{\alpha}{\myrst(\lambda)})$.

\end{mydef}

\begin{myremark}

Recall that a $\Levy$ process $\yy=(\yy(t),\, t\geq0)$ is called
\emph{stable} if each $\yy(t)$ is stable. In this case, evidently,
the process $\yy$ belongs to the domain of \emph{normal} attraction
of $\yy(1)$.

\end{myremark}

According to assumptions of Subsections~\ref{sub:free-evol-nonsym}
and \ref{sec:sym-synch-mod-N} $\xc_{j}\sim\LP(\eLH(\lambda))$, i.e.,
the free dynamics of any component of $x(t)$ is the $\Levy$ process
with the common $\Levy$ exponent $\eLH:\,\RR^{d}\rightarrow\CC$.

\paragraph{Assumption~D.}

There exist a stable law $\St{\alpha}{\zLH(\lambda)}$ in~$\RR^{d}$
such that any component $\xc_{j}(t)$ of the free dynamics $\xc(t)=(\xc_{1}(t),\ldots,\xc_{N}(t))$
belongs to the \emph{domain of attraction} of $\St{\alpha}{\zLH(\lambda)}$. 

\smallskip 

\noindent According to Definition~\ref{def-DOA-stable} under Assumption~D
there exist sequences $\left\{ b_{n}\right\} $ and $\left\{ C_{n}\right\} $
such that for all $\lambda\in\RR^{d}$ \begin{equation}
\exp\left(n\eLH(\lambda/b_{n})+i\skpr{C_{n}}{\lambda}\right)\rightarrow\exp\zLH(\lambda)\quad\mbox{as}\quad n\rightarrow\infty\,.\label{eq:dom-n-attr-zeta-2}\end{equation}

\paragraph{Assumption~DN.}

There exist a stable law $\St{\alpha}{\zLH(\lambda)}$ in~$\RR^{d}$
such that any component $\xc_{j}(t)$ of the free dynamics $\xc(t)=(\xc_{1}(t),\ldots,\xc_{N}(t))$
belongs to the domain of \emph{normal attraction} of $\St{\alpha}{\zLH(\lambda)}$. 

\smallskip

Define a stochastic process $d_{jk}^{\circ,N}(t)=\xc_{j}(t)-\xc_{k}(t),$
$t\geq0.$ According to Lemma~\ref{l-V-V-i-iii}(iii) $d_{jk}^{\circ,N}\sim\LP\left(2\Re\eLH(\lambda)\right)$,
i.e., all $d_{jk}^{\circ,N}(t)$ are $\Levy$ processes in $\RR^{d}$
with the common characteristic function \begin{equation}
\left|\chf{\,}t{\lambda}\right|^{2}=e^{-t\etho(\lambda)}\,\label{eq:f2-exp-eta}\end{equation}
where \[
\etho(\lambda):=-(\eLH(\lambda)+\eLH(-\lambda))=-2\Re\eLH(\lambda).\]
The function (\ref{eq:f2-exp-eta}) is real and, moreover, $\etho(\lambda)\geq0$
for all $\lambda\in\RR^{d}$ since $e^{-t\etho(\lambda)}$ is a characteristic
function of some probability distribution. Hence distributions of
the increments of $\left(d_{jk}^{\circ,N}(t),\, t\geq0\right)$ are
symmetric. 

It follows from Lemma~\ref{l-V-V-i-iii}(i) that if Assumption~D
holds then the process $d_{jk}^{\circ,N}(t)$ belongs to the \emph{domain
of attraction} of the \emph{symmetric} stable law $\St{\alpha}{-\zto(\lambda)}$
with the characteristic function $e^{-\zto(\lambda)}$ where \begin{equation}
\zto(\lambda):=-(\zLH(\lambda)+\zLH(-\lambda))=-2\Re\zLH(\lambda)\label{eq:zeta-zeta0-2Re}\end{equation}
and $\zLH$ is the same as in Assumption~D. It is evident that $\zto(\lambda)\geq0$
for all $\lambda\in\RR^{d}$. According to (\ref{eq:dom-n-attr-zeta-2})
the conclusion that $d_{jk}^{\circ,N}(1)\in\DOA\left(\St{\alpha}{-\zto(\lambda)}\right)$
implies that \[
\exp\left(-n\etho(\lambda/b_{n})\right)\rightarrow\exp(-\zto(\lambda))\quad\quad(n\rightarrow\infty)\,\]
for the same sequence~$\{b_{n}\}$ as in~(\ref{eq:dom-n-attr-zeta-2}).
Applying the logarithmic function to the above convergence we get
that for any $\lambda\in\RR^{d}$\begin{equation}
n\etho(\lambda/b_{n})\rightarrow\zto(\lambda)\quad\quad(n\rightarrow\infty)\,.\label{eq:n-eta-bn-zeta}\end{equation}

Similarly, using the item~(ii) of Lemma~\ref{l-V-V-i-iii} we get
that under Assumption~DN the process $d_{jk}^{\circ,N}(t)$ belongs
to the\emph{ }domain of\emph{ normal attraction} of the same stable
law $\St{\alpha}{-\zto(\lambda)}$. Of course, under Assumption~DN
the condition~(\ref{eq:n-eta-bn-zeta}) takes the following form\[
n\etho(\lambda/n^{1/\alpha})\rightarrow\zto(\lambda)\quad\quad(n\rightarrow\infty)\,.\]

\subsubsection*{The space scaling}

Consider the model $\mathcal{GGS}_{N}(\eLH;F)$, i.e., the $N$-component
synchronization system $x(t)=\left(x_{1}(t),\ldots,x_{N}(t)\right)$
which satisfies the symmetry assumptions of Subsection~\ref{sec:sym-synch-mod-N}.

\begin{mytheo}  \label{t:asympt-char-f}

Let Assumption~D hold with some $\zLH(\lambda)$. Let $\left\{ b_{n}\right\} $
be the normalizing sequence in~(\ref{eq:dom-n-attr-zeta-2}). Rescale
the system $x(t)=\left(x_{1}(t),\ldots,x_{N}(t)\right)$ as follows
\[
y^{(N)}(t)=\frac{x(t)}{b_{N}},\qquad y^{(N)}(t)=\left(y_{1}^{(N)}(t),\ldots,y_{N}^{(N)}(t)\right).\]
Let $Q_{N,t}$ be the probability law of the rescaled differences
$y_{j}^{(N)}(t)-y_{k}^{(N)}(t)$. Then for any fixed $N\geq2$ the
weak limit of $Q_{N,t}$ exists, \[
Q_{N,t}\stackrel{w}{\rightarrow}Q_{N,\infty}\quad as\quad t\rightarrow+\infty,\]
and the characteristic function of the limiting distribution has asymptotically
explicit form as $N\rightarrow\infty$\begin{equation}
\int_{\RR^{d}}\exp\left(i\skpr{\lambda}y\right)\, Q_{N,\infty}(dy)\,\rightarrow\,\frac{1}{1+\frac{1}{2}m\zto(\lambda)}\,.\label{eq:as-har-f}\end{equation}
Here the real function~$\zto=\zto(\lambda)$, $\lambda\in\RR^{d}$,
is the same as in~(\ref{eq:zeta-zeta0-2Re}) and $m$ is defined
by~(\ref{eq:moments-interv}).

\end{mytheo} 

We have an immediate corollary of this theorem under the stronger
condition that the synchronized system $x(t)$ satisfied to Assumption
DN with respect to some stable law $\St{\alpha}{\zLH(\lambda)}$.
In this case $b_{n}=n^{1/\alpha}$ and the statement of Theorem~\ref{t:asympt-char-f}
is true for the rescaled synchronization system \[
y^{(N)}(t)=\frac{x(t)}{N^{1/\alpha}},\qquad y^{(N)}(t)=\left(y_{1}^{(N)}(t),\ldots,y_{N}^{(N)}(t)\right).\]
This result can be interpreted as follows: distances between components
in the synchronized system are of order $N^{1/\alpha}$ provided the
free dynamics belongs to the domain of normal attraction of an $\alpha$-stable
law in the sense of~\cite{Gned-Kolm-1968-eng}.

Note also that the $\Levy$ continuity theorem and (\ref{eq:as-har-f})
imply the weak convergence of $Q_{N,\infty}$ to some probability
law $Q_{\infty,\infty}$ in $\RR^{d}$ having the characteristic function
$\left(1+\frac{1}{2}m\zto(\lambda)\right)^{-1}$.

\subsection{Free dynamics attracting to stable laws. Linnik distributions}

\label{sub:Linnik-distrib}

\subsubsection*{Symmetric stable laws}

It is very useful to illustrate the result of Theorem~\ref{t:asympt-char-f}
by different concrete examples of free dynamics. Before doing this
we need to recall some classical results about representation of stable
laws. It is known~\cite{Applebaum-Levy-book,Sato-Levy} that the
characteristic function of a $d$-dimensional \emph{symmetric} $\alpha$-stable
law has the following form
\begin{itemize}
\item for $0<\alpha<2$~ (the \emph{heavy tail} case): \begin{equation}
e^{-t\zto(\lambda)}=\exp\left(-t\int_{S^{d-1}}\left|\skpr{\lambda}{\xi}\right|^{\alpha}\nu(d\xi)\,\right)\label{eq:stable-heavy-tails}\end{equation}
where $S^{d-1}$ is the unit sphere in~$\RR^{d}$ and $\nu$ is some
finite measure on $S^{d-1}$,
\item for $\alpha=2$~ (the \emph{Gaussian} case): \begin{equation}
e^{-t\zto(\lambda)}=\exp\left(-t\skpr{A\lambda}{\lambda}/2\,\right)\label{eq:al-2-stable-Gauss}\end{equation}
where $A$ is a positive definite symmetric $d\times d$ matrix.
\end{itemize}
Corresponding formula for rotationally invariant $\alpha$-stable
laws, $0<\alpha\leq2$, is simpler:\[
e^{-t\zto(\lambda)}=\exp\left(-tc^{\alpha}\left|\lambda\right|^{\alpha}\,\right),\qquad c>0,\quad\lambda\in\RR^{d},\quad\left|\lambda\right|=\sqrt{\skpr{\lambda}{\lambda}}\,.\]

It is clear from Theorem~\ref{t:asympt-char-f} that any of functions
\begin{equation}
\int_{S^{d-1}}\left|\skpr{\lambda}{\xi}\right|^{\alpha}\nu(d\xi),\quad\skpr{A\lambda}{\lambda}/2,\quad c^{\alpha}\left|\lambda\right|^{\alpha}\label{eq:list-sym-st}\end{equation}
can participate as $\zto(\lambda)$ in the limit (\ref{eq:as-har-f}).
Indeed, to see this one should consider the free dynamics $\xc(t)$
driven by symmetric stable $\Levy$ processes with the $\Levy$ exponent
$\etho^{\circ}(\lambda)=-\frac{1}{2}\zto(\lambda)$ where $\zto(\lambda)$
is taken from the list~(\ref{eq:list-sym-st}). 

\begin{myremark} \label{rem-Geom-stable-d}

Note that the limiting characteristic function in~(\ref{eq:as-har-f})
has the form \[
\frac{1}{1-\log\phi(\lambda)}\,,\qquad\lambda\in\RR^{d},\]
where $\phi(\lambda)$ is a characteristic function of some symmetric
$\alpha$-stable distribution. As it follows from~\cite[Prop.~1]{Mittnik-Rach}
the class of limiting laws in~(\ref{eq:as-har-f}) are exactly the
symmetric \textbf{geometric stable distributions} (GSDs). The GSDs
are obtained as limiting laws of appropriately normalized \emph{random}
sums of i.i.d. random vectors in~$\RR^{d}$ where the \emph{number}
of summands is \emph{geometrically distributed} and independent of
the summands. There is a large bibliography devoted to this topic,
see, for example, \cite{Klebanov,hazod-Khokhlov,Kalash-Geom,Gnedenko-Korolev,Mittnik-Rach,KK-Lapl-distr-book,KMPS-Op-Geo-stable-laws}.

\end{myremark}

\subsubsection*{Free dynamics of the Gaussian type. The Laplace distribution}

Let $\xc_{j}(t)$, $\jolN$, be the same as in example~(\ref{eq:Win-xo}),
\[
d\,\xc_{j}(t)=\sigma dB_{j}(t)+b\, dt,\]
where $\sigma$ is a real $d\times d$ matrix, $b\in\RR^{d}$ and
$B_{j}(t)=(B_{j}^{1}(t),\ldots,B_{j}^{d}(t))\in\RR^{d}$ are independent
standard $d$-dimensional \emph{Brownian motions} . We know from~(\ref{eq:Win-xo})
that any $\xc_{j}(t)$ is a $\Levy$ process $\LP(\eLH)$ determined
by the $\Levy$ exponent $\eLH(\lambda)=i\skpr b{\lambda}-\frac{1}{2}\skpr{\sigma\sigma^{T}\lambda}{\lambda}$.
Using~(\ref{eq:psi-stable}) it is easy to check that $\LP(\eLH)$
is stable with $\alpha=2$. Hence \[
\zto(\lambda)=\etho(\lambda)=-2\Re\eLH(\lambda)=\skpr{\sigma\sigma^{T}\lambda}{\lambda}.\]
By Theorem~\ref{t:asympt-char-f} the proper scaling for differences
$d_{jk}^{(N)}(t)=x_{j}(t)-x_{k}(t),$ is $N^{-1/2}$. Namely, $d_{jk}^{(N)}(t)/\sqrt{N}$
weakly converges to some law $Q_{N,\infty}$ as $t\rightarrow\infty$.
Letting $N\rightarrow\infty$ we get from~(\ref{eq:as-har-f}) that
$Q_{N,\infty}$ weakly converges to the distribution with characteristic
function \[
\,\frac{1}{1+\frac{1}{2}m\skpr{\sigma\sigma^{T}\lambda}{\lambda}}\,,\qquad\lambda\in\RR^{d}\,.\]
In the case $d=1$ this characteristic function takes the form $(1+\frac{1}{2}m\sigma^{2}\lambda^{2})^{-1}$
and corresponds to the Laplace distribution with density \begin{equation}
p_{L}(y)=\frac{1}{2c_{0}}\, e^{-\left|y\right|/c_{0}},\quad\quad-\infty<y<+\infty,\quad c_{0}=\sigma\sqrt{\frac{m}{2}}\,.\label{eq:pL-Lapl-distr-sgm}\end{equation}
This result generalizes the result obtained in~\cite{asmda-2013-man}
and cited in Introduction of the current paper. Indeed, putting $d=1$
and $b=0$ we have equivalence of the following models

\[
\mathcal{GMS}_{N}(-\frac{1}{2}\sigma^{2}\lambda^{2};m)=\mathcal{BM}_{N}(\sigma,N/m).\]
Using self-similarity of the Wiener process one can derive that the
intrinsic scale of $\mathcal{BM}_{N}(\sigma,N/m)$ is $N^{1/2}$ times
smaller than the intrinsic scale of the model $\mathcal{BM}_{N}(\sigma,m^{-1})$
studied in~\cite{asmda-2013-man}.

\subsubsection*{One-dimensional random walks. Linnik distribution}

Let $d=1$ and the free dynamics of each component $\xc_{j}(t)$ be
a continuous time symmetric\emph{ random walk} with the Markov generator
\begin{equation}
\left(Lf\right)(y)=\beta\int_{\RR}\left(f(y+q)-f(y)\right)\,\mu(dq),\qquad f\in C_{b}(\RR,\RR).\label{eq:Lf-beta-int-mu}\end{equation}
Here $\beta>0$ is the intensity of jumps and $\mu(dq)=\frac{1}{2}\bb\left|q\right|^{-1-\bb}1_{\left\{ \left|q\right|\geq1\right\} }dq$
is the distribution of an individual jump $x\mapsto x+q$. This is
a one-dimensional subcase of the example (\ref{eq:mu-j-q-d-al}).
Please, note that the sequence $\{\tT_{n}\}$ is considered here under
general assumptions of Subsections~\ref{sub:syncr-ep-T-M} and~\ref{sub:inter-event-symmetr}.
Hence the synchronization system $x(t)$ is not Markovian while the
free dynamics $\xc(t)$ is a Markov process. 

The jump distribution $\mu(dq)$ has the {}``Pareto tails'' and,
as it will be seen below, the conditions of Theorem~\ref{t:asympt-char-f}
can be easily checked. Let $\mxi$ be a random variable with distribution
$\mu(dq)$. If $\bb>2$ then $\mxi$ has a finite variation $\mVar_{0}=\Var\left(\mxi\right)=a/(a-2)$.
It follows from~\cite[\S~35, Th.~4]{Gned-Kolm-1968-eng} that $\mxi\in\DONA\left(\mathcal{N}(0,\mVar_{0})\right)$
where $\mathcal{N}(0,\mVar_{0})$ is the Gaussian law with zero mean
and variance~$\mVar_{0}$. It follows from~(\ref{eq:Lf-beta-int-mu})
that $\xc_{j}(t)$ is a compound Poisson process, i.e., \[
\xc_{j}(t)\sim\sum_{r=1}^{N_{\beta}(t)}\mxi_{r}\,\]
where $\left(N_{\beta}(t),\, t\geq0\right)$ is the Poisson process
with intencity~$\beta$ and $\mxi_{1},\ldots,\mxi_{r},\ldots$ are
inde\-pen\-dent copies of~$\mxi$. It is well known that $N_{\beta}(t)\sim\beta t$
as $t\rightarrow\infty$. Arguments similar to \cite[\S~4.4]{Meer-Sikor}
show that $\xc_{j}(1)\in\DONA\left(\mathcal{N}(0,\mVar)\right)$ where
$\mVar=\beta a/(a-2)>0$. Hence Assumption~DN holds with $\zLH(\lambda)=-\frac{1}{2}\mVar\lambda^{2}$,
$\alpha=2$, and we can apply Theorem~\ref{t:asympt-char-f}. It
is readily seen that $\zto(\lambda)=\mVar\lambda^{2}$. The distribution
of rescaled differences $d_{jk}^{(N)}(t)/\sqrt{N}$ converges to some
law $Q_{N,\infty}$ as $t\rightarrow\infty$. The sequence of laws
$Q_{N,\infty}$ converges as $N\rightarrow\infty$ to the Laplace
distribution~(\ref{eq:pL-Lapl-distr-sgm}) where $\sigma$ is replaced
by $\sqrt{D}$.

If $0<\bb<2$ (the \emph{heavy tail} case) then by \cite[\S~35, Th.~5]{Gned-Kolm-1968-eng}
the random variable~$\mxi$ with the distribution $\mu(dq)$ belongs
to the domain of \emph{normal} attraction of a symmetric $\bb$-stable
law. As in the above paragraph we conclude that $\xc_{j}(1)$ also
belongs to the domain of normal attraction of some symmetric $\bb$-stable
law. So Assumption~DN holds. Again Theorem~\ref{t:asympt-char-f}
implies that rescaled differences $d_{jk}^{(N)}(t)/N^{1/\bb}$ converge
in law as $t\rightarrow\infty$ to some distribution $Q_{N,\infty}$.
The sequence $Q_{N,\infty}$ converges as $N\rightarrow\infty$ to
a symmetric law which characteristic function is \begin{equation}
\frac{1}{1+c^{\bb}\left|\lambda\right|^{\bb}}\,,\qquad\lambda\in\RR,\label{eq:Linnik-ch-f}\end{equation}
for some $c=c(\bb,\beta,m)>0$. This is characteristic function of
the famous symmetric Linnik distribution~\cite{Linnik53} usually
denoted as $\mathcal{L}_{\bb,c}$. It is known~\cite{Lukacs,KK-Lapl-distr-book}
that this distribution is unimodal, absolutely continuous, geometric
stable (see Remark~\ref{rem-Geom-stable-d}) and infinitely divisible.
If $0<\bb<2$ then the Linnik distribution has heavy tails~\cite{Erdog-Ostrovski,KK-Lapl-distr-book}:~
$q^{\bb}\sP\left(\mathcal{L}_{\bb,c}>q\right)\sim\mbox{const}$~
as $q\rightarrow+\infty$. For $\bb=2$ the law (\ref{eq:Linnik-ch-f})
is the Laplace distribution.

It is not hard to modify this example to obtain domains of non-normal
attraction. Let $0<\bb<2$. Using notation of~\cite[Th.~4]{Khokhlov-FPM}
we introduce a probability measure $\mu(dq)$ on $\RR^{1}$ such that
\begin{eqnarray*}
\mu\left(\,(-\infty,-q]\,\right) & = & q^{-\bb}\left(c_{1}+h_{1}(q)\right)\, L(q),\\
\mu\left(\,[q,+\infty)\,\right) & = & q^{-\bb}\left(c_{2}+h_{2}(q)\right)\, L(q),\end{eqnarray*}
where $c_{1}\geq0$, $c_{2}\geq0$, $c_{1}+c_{2}>0$, $L(q)$ is a
slowly varying function and $h_{i}(q)\rightarrow0,$ $i=1,2$, as
$q\rightarrow\infty.$ Let each component $\xc_{j}(t)$ be a continuous
time\emph{ random walk} with the Markov generator~(\ref{eq:Lf-beta-int-mu}).
Again denote by $\mxi$ a random variable with distribution $\mu(dq)$.
The classical results~\cite{Gned-Kolm-1968-eng} states that there
exists a stable law $U_{\bb,c_{1},c_{2}}$ in $\RR^{1}$ with index
of stability $\bb$ such that $\mxi\in\DOA(U_{\bb,c_{1},c_{2}})$.
This statement implies that $\xc_{j}(1)$ belongs to the $\DOA(U_{\bb,\beta c_{1},\beta c_{2}})$.
We don't need an explicit definition of $U_{\bb,c_{1},c_{2}}$ here.
The law $U_{\bb,c_{1},c_{2}}$ is symmetric iff $c_{1}=c_{2}$. It
is important to note~\cite{Khokhlov-FPM,Sam-Taqqu} that the choice
of normalizing sequence $\left\{ b_{N}\right\} $ arising in Theorem~\ref{t:asympt-char-f}
depends on the function $L(q)$. If there exists a limit \[
L(q)\rightarrow L(\infty),\qquad(q\rightarrow\infty),\qquad0<L(\infty)<\infty,\]
then Assumption~DN holds and $b_{N}=N^{1/\alpha}$. In the other
case $\xc_{j}(1)\in\DOA\backslash\DONA$. In general, $b_{N}=N^{1/\alpha}\ell(N)$
where $\ell=\ell(N)$ is a slowly varying function at infinity. See~\cite{Khokhlov-FPM,Sam-Taqqu}
for details.

In any case the limiting characteristic function in~(\ref{eq:as-har-f})
is from the class of Linnik distributions~(\ref{eq:Linnik-ch-f}).

\subsubsection*{Multi-dimensional random walks with heavy-tailed jumps}

We consider a special subclass of \emph{random walks} $\xc_{j}(t)$
in $\RR^{d}$, $d\geq2$, introduced in~(\ref{eq:Lj-f-y}). According
to the symmetry assumption we put \[
\beta_{j}=\beta,\qquad\mu_{j}(dq)=\mu(dq)\qquad\forall\jolN.\]
As in previous examples we restrict ourself to consideration of power
law jumps. Let $\mxi\in\RR^{d}$ denotes a random vector with distribution
$\mu(dq)$, $q\in\RR^{d}$: $\sP\left(\mxi\in G\right)=\mu(G)$ for
any Borel set $G\in\mathcal{B}\left(\RR^{d}\right)$. Following \cite[\S~6.4]{Meer-Sikor}
we represent it as $\mxi=W\Theta$ where $W$ is a scalar random variable
and $\Theta\in S^{d-1}$ is a random vector taking values on the unit
sphere in~$\RR^{d}$. Assume that $W$ and $\Theta$ are independent,
and \[
\sP\left(W>R\right)=CR^{-\alpha},\quad\quad R\geq R_{0}>0,\qquad\sP\left(W\geq0\right)=1,\]
\[
\sP\left(\Theta\in B\right)=M(B),\qquad B\in\mathcal{B}(S^{d-1}),\]
where $C>0$ is some constant, $C\leq R_{0}^{\alpha}$, and $M(\cdot)$
is a probability measure on~$S^{d-1}$. Consider only the\emph{ heavy
tail} case $\alpha\in(0,2)\backslash\left\{ 1\right\} $ excluding
$\alpha=1$ for breavity of formulae. By Theorem~6.17 in~\cite{Meer-Sikor}
$\mxi\in\DONA(\St{\alpha}{\rho_{\alpha}^{\circ}(\lambda)})$ where
\begin{equation}
\rho_{\alpha}^{\circ}(\lambda)=-C\, K_{\alpha}\int_{S^{d-1}}\left|\skpr{\lambda}{\theta}\right|^{\alpha}\,\left(1-i\,\mbox{sgn}\,\skpr{\lambda}{\theta}\tan\frac{\pi\alpha}{2}\right)\, M(d\theta),\label{eq:rho-0-alpha}\end{equation}
for some $K_{\alpha}>0.$ The same arguments as for the one-dimensional
random walks imply that $\xc_{j}(1)\in\DONA\left(\St{\alpha}{\beta\rho_{\alpha}^{\circ}(\lambda)}\right)$.
Hence Assumption~DN holds with $\zLH(\lambda)=\beta\rho_{\alpha}^{\circ}(\lambda)$.
After calculating \[
\zto(\lambda)=-2\Re\zLH(\lambda)=2C\, K_{\alpha}\int_{S^{d-1}}\left|\skpr{\lambda}{\theta}\right|^{\alpha}\, M(d\theta)\]
(compare with~(\ref{eq:stable-heavy-tails})) we are ready to apply
Theorem~\ref{t:asympt-char-f}. We conclude that the rescaled differences
$d_{jk}^{(N)}(t)/N^{1/\alpha}$ converge in law as $t\rightarrow\infty$
to some distribution $Q_{N,\infty}$. If $N\rightarrow\infty$ then
$Q_{N,\infty}$ is approximated by the distribution with characteristic
function\[
\frac{1}{1+\frac{1}{2}m\int_{S^{d-1}}\left|\skpr{\lambda}{\theta}\right|^{\alpha}\,\nu_{\alpha}(d\theta)}\,\]
where $\nu_{\alpha}(d\theta)=2C\, K_{\alpha}M(d\theta)$. In the case
$\alpha=1$ we have essentially the same final conclusion but the
intermediate formula~(\ref{eq:rho-0-alpha}) is different. 

If $\alpha\geq2$ them the corresponding analysis is based on the
multi-dimensional Central Limit Theorem. Here $b_{N}=N^{1/2}$ and
Assumption~DN holds for $d$-dimensional Gaussian law. We omit details.

\subsection{Intrinsic scales based on matrix transformations. Jurek coordinates}

\label{sub:matrix-Jurek}

Theorem~\ref{t:asympt-char-f} justifies the existence of a natural
space scale for a large $N$-components synchronization system. This
scale is uniform in any of $d$ coordinate axes in $\RR^{d}$ because
the scaling transformation $y=b_{N}^{-1}x$ is the multiplication
by a scalar value $b_{N}^{-1}$. 

It is also interesting to find conditions when large synchronized
systems {}``are con\-cent\-rated'' in space domains which change
non-uniformly in different coordinate direc\-tions as $N\rightarrow\infty$.
Recalling that $\mySect$~\ref{sub:Intr-scales} is related with
attraction to stable laws in~$\RR^{d}$ it is clear that one can
look for generalizations of Theorem~\ref{t:asympt-char-f} by considering
the domains of attraction of \emph{operator stable laws} (OSLs) .

\begin{mydef} \label{def-GDOA-op-stable} 

A random vector $U\in\RR^{d}$ is operator stable if it has a \textbf{generalized}
domain of attraction, i.e., if there is a random vector $V$ and sequences
$\left\{ B_{n}\right\} $ of linear operators $B_{n}:\,\RR^{d}\rightarrow\RR^{d}$
and nonrandom vectors $\left\{ C_{n}\right\} $, $C_{n}\in\RR^{d}$,
such that \begin{equation}
B_{n}\left(V_{1}+\cdots+V_{n}\right)+C_{n}\stackrel{d}{\longrightarrow}U\label{eq:B-VV-CU}\end{equation}
where $V_{1},\ldots V_{n},\ldots$ are independent copies of~$V$.
The random vector~$V$ is said to be in the generalized domain of
attraction of the stable vector~$U$, the short notation for it is
$V\in\GDOA(U)$.

\end{mydef} We see that operator stable laws arise as limiting distributions
of matrix-normalized sums of i.i.d. random vectors. The study of OSLs
was originated by G.N.~Sakovich, the disciple of B.V.~Gnedenko,
and M.~Sharpe. Here we cannot go too deeply in details of this vast
theory and refer to~\cite{Sakov-60,Sakov-diss,Jurek-Mason,Meer-Sikor}.
Below we list a limited number of facts on OSLs which are necessary
to state our result. We will consider only \emph{full} OSLs . The
probability distribution of a random vector $U$ on $\RR^{d}$ is
full if $\skpr{\lambda}U$ is nondegenerate for every $\lambda\in\RR^{d}\backslash\left\{ 0\right\} $.

The simplest examples of OSLs are laws $U$ in $\RR^{d}$ with marginal
stable distributions which possess a stability property similar to~(\ref{eq:UUnUD}),
\begin{equation}
U_{1}+\cdots+U_{n}\eqd n^{E}U+D_{n},\label{eq:UU-nE-U}\end{equation}
where $E$ is a diagonal matrix $E=\mbox{diag}\left(\alpha_{1}^{-1},\ldots,\alpha_{d}^{-1}\right)$,
$\alpha_{i}\in(0,2]$. In this case, evidently, $B_{n}=n^{-E}$. To
have OSLs with dependent coordinated one should replace in~(\ref{eq:UU-nE-U})
the diagonal matrix $n^{E}$ by the multiplier $n^{B}$ where $B$
is a real $d\times d$-matrix with whose eigenvalues all have real
part in $[\frac{1}{2},+\infty)$, see~\cite{Meer-Sikor}. The matrix
$n^{B}$ is defined by using the matrix exponent as $n^{B}=\exp\left(B\log n\right)$.
Any full operator stable~$U$ is infinitely divisible~\cite[\S~4.2]{Jurek-Mason},
hence its characteristic function has the form \[
\myps_{U}(\lambda)=\sE\exp\left(i\skpr{\lambda}U\right)=\exp\zto_{U}(\lambda).\]
We will write $U\sim\mathcal{OS}(\zto_{U}(\lambda))$ to have a short
notation for this situation.

We are ready to state a\textbf{ }\textbf{\emph{result}}\textbf{ generalizing
Theorem~\ref{t:asympt-char-f}}. Consider the symmetric synchronization
model $x(t)\in\RR^{dN}$ with the free dynamics $\xc_{j}\sim\LP(\eLH(\lambda))$,
$\jolN$. \emph{Assume} that there exist an operator stable law $\mathcal{OS}(\zLH(\lambda))$
in~$\RR^{d}$ such that any component $\xc_{j}(1)\in\GDOA\left(\,\mathcal{OS}(\zLH(\lambda))\,\right)$.
According to Definition~\ref{def-GDOA-op-stable} this assumption
means that there exist sequences $\left\{ B_{n}\right\} $ and $\left\{ C_{n}\right\} $
such that for all $\lambda\in\RR^{d}$ \begin{equation}
\exp\left(n\eLH(B_{n}^{T}\lambda)+i\skpr{C_{n}}{\lambda}\right)\rightarrow\exp\zLH(\lambda)\quad\mbox{as}\quad n\rightarrow\infty\,.\label{eq:dom-n-attr-zeta-2-2}\end{equation}
Define a transformed system \begin{equation}
y^{(N)}(t)=\left(B_{N}x_{1}(t)\,,\ldots,\, B_{N}x_{N}(t)\right).\label{eq:y-BN-x}\end{equation}
Then the differences $d_{jk}^{(N)}(t)=y_{j}^{(N)}(t)-y_{k}^{(N)}(t)=B_{N}\left(x_{j}(t)-x_{k}(t)\right)$
converge in law as $t\rightarrow\infty$ to some distribution $Q_{N,\infty}$.
If $N\rightarrow\infty$ then $Q_{N,\infty}$ is approximated by some
distribution with characteristic function given in the explicit form:\begin{equation}
\int_{\RR^{d}}\exp\left(i\skpr{\lambda}y\right)\, Q_{N,\infty}(dy)\,\rightarrow\,\frac{1}{1+\frac{1}{2}m\zto(\lambda)}\,,\qquad\lambda\in\RR^{d},\label{eq:as-har-f-op-stable}\end{equation}
where~$\zto(\lambda)=-2\Re\zLH(\lambda)$.

The proof of this generalization is very similar to the proof of Theorem~\ref{t:asympt-char-f}
and is based on the representation for $\chi_{N}(\infty;\lambda)$
of Theorem~\ref{t:chi-N-inf-lambda}. So we omit it. 

We end this subsection by two remarks. The characterization of GDOA
in the operator stable case and the description of all possible functions~$\zto(\lambda)$
in (\ref{eq:as-har-f-op-stable}) are not easy. They demand many additional
constructions and are out of scope of this paper. We refer interested
readers to~\cite{Jurek-Mason,Meer-Sikor}. 

In the case when $B_{n}=n^{-B}$ the transformation~(\ref{eq:y-BN-x})
is deeply connected with so called \emph{Jurek coordinates}. The Jurek
coordinates in $\RR^{d}$ is a pair $(r,\Theta)$ such that $y=r^{B}\Theta$,
where $y\in\RR^{d}$, $r\geq0$ and $\Theta\in S^{d-1}$. Details
can be found in~\cite{Jurek-Mason,Meer-Sikor}.

\subsection{The Markovian case}

\label{sub:Markov-case}

Here we consider $\mathcal{GMS}_{N}(\eLH;m)$, the symmetric $N$-component
synchronization model in any dimension $d$ with the special choice
of inter-event distribution: \begin{eqnarray}
F(s) & = & \sP\left\{ \Delta_{n}^{(k)}\leq s\right\} =1-\exp\left(-s/m\right),\quad\quad s\geq0,\quad m>0.\label{eq:inter-event-Markov-case}\end{eqnarray}
This is the exponential distribution with the mean $m$. In this case
the sequence $\left\{ \tT_{n}\right\} $ is the Poissonian flow of
intensity~$N/m$ and $x(t)$ is a Markov process.

In the Markovian case it is possible to precise main results of Subsections~\ref{sub:Lim-distrib}
and \ref{sub:Intr-scales}. Theorem~\ref{t:chi-N-inf-lambda} is
replaced by the following one.

\begin{mytheo}  \label{t-Markov-chi-expl}

For the Markovian symmetric synchronization model $\mathcal{GMS}_{N}(\eLH;m)$
\[
\chi_{N}(+\infty;\lambda)=\frac{1}{1+\frac{1}{2}(N-1)m\etho(\lambda)}\]
where the function $\etho(\lambda)$ is the same as in Theorem~\ref{t:chi-N-inf-lambda}.
\end{mytheo}  The proof of this theorem is given at the end of Subsection~\ref{sub:repr-ch-f}.

The next theorem holds for finite $N$. It immediately follows from
Theorem~\ref{t-Markov-chi-expl}.

\begin{mytheo} \label{t-Markov-stable}

Let a Markovian $N$-component symmetric synchronization model $\mathcal{GMS}_{N}(\eLH;m)$
be such that its free dynamics $\xc(t)$ is an $\alpha$-\emph{stable}
$\Levy$ process, $0<\alpha\leq2$. Then for any fixed $N$ the distribution
of rescaled differences $d_{jk}^{(N)}(\infty)/\left(N-1\right)^{1/\alpha}$
does not depend on~$N$.

\end{mytheo} 

The Markov assumption is essential for Theorem~\ref{t-Markov-stable}.
For the non-Markov case the statement~(\ref{eq:as-har-f}) of Theorem~\ref{t:asympt-char-f}
is asymptotic and does not hold for finite $N$.

Theorem~\ref{t-Markov-stable} generalizes results of the paper~\cite{asmda-2013-man}
where the role of the $\alpha$-stable free dynamics was played by
Brownian motions ($\alpha=2$).

For the Markovian symmetric model the function $\chi_{N}(t;\lambda)$
satisfies to the following differential equation\begin{equation}
\frac{d}{dt}\,\chi_{N}(t;\lambda)=-\qn(\lambda)\,\chi_{N}(t;\lambda)+\wn,\label{eq:chi-dif-eq}\end{equation}
where \[
\wn=\frac{2}{(N-1)m},\qquad\qn=\etho(\lambda)+\wn.\]
This equation directly follows from the representation for $\chi_{N}(t;\lambda)$
which will be obtained in Subsection~\ref{sub:repr-ch-f}. In particular,
the statement of Theorem~\ref{t-Markov-chi-expl} easily follows
from this equation.

It is important to note that for non-Markovian models the function
$\chi_{N}(t;\lambda)$ don't satisfy to any differential equation
of such type.

\subsection{Some generalizations}

\label{sub:Some-generalizations}

According to Subsection~\ref{sec:sym-synch-mod-N} and Assumption~P1
the general (non-Markovian) symmetric synchronization model $x(t)$
is determining by the quadruple $\left(N,\,\eLH(\lambda),\, p(s),\Law{x(0)}\right)$.
Here we briefly discuss a possibility to extend our asymptotic results
to the case \[
\left(N,\,\eLH_{N}(\lambda),\, p_{N}(s),\Law{x(0)}\right)\]
 when $\eLH(\lambda)$ and $p(s)$, the functions defining the dynamics,
depend on~$N$. The main task is \emph{to generalize} Theorem~\ref{t:chi-N-inf-lambda}.
Note that this problem is interesting only for \emph{non-Markovian
models}. Indeed, in the Markovian case Theorem~\ref{t-Markov-chi-expl}
already gives the exact and explicit answer to the question.

We will restrict ourself to the special situation when \begin{equation}
p_{N}(s)=\skb p(\skb s)\label{eq:pN-p-scale}\end{equation}
for some sequence $\left\{ \skb\right\} $, $\skb>0$. This situation
corresponds to the rescaling of the time~$t$ and is quite simple.
Obviously, \begin{equation}
m_{N,1}=m/\skb\,\label{eq:mN1-m1}\end{equation}
where \[
m_{N,1}=\int_{0}^{\infty}sp_{N}(s)\, ds,\qquad m_{1}=\int_{0}^{\infty}sp(s)\, ds\,.\]
The main idea is to compare models with different quadruples. Indeed,
in distributional sense \[
x(t)\,\left|_{\left(N,\,\eLH_{N}(\lambda),\, p_{N}(s),\Law{x(0)}\right)}\right.=x(\skb t)\,\left|_{\left(N,\,\eLH_{N}(\lambda)/\skb,\, p(s),\Law{x(0)}\right)}\right.\,.\]
Hence \[
\chi_{N}(t;\lambda)\,\left|_{\left(N,\,\eLH_{N}(\lambda),\, p_{N}(s),\Law{x(0)}\right)}\right.=\chi_{N}(\skb t;\lambda)\,\left|_{\left(N,\,\eLH_{N}(\lambda)/\skb,\, p(s),\Law{x(0)}\right)}\right..\]
Let the probability density function $p(s)$ satisfies to Assumptions~P1--P3.
Then by Lemma~\ref{l-ch-f-limit}\[
\chi_{N}(+\infty;\lambda)\,\left|_{\left(N,\,\eLH_{N}(\lambda),\, p_{N}(s),\Law{x(0)}\right)}\right.=\chi_{N}(+\infty;\lambda)\,\left|_{\left(N,\,\eLH_{N}(\lambda)/\skb,\, p(s),\Law{x(0)}\right)}\right..\]
Note that these limiting characteristic functions do not depend on
the initial distribution $\Law{x(0)}$ so we can omit it in the notation.
From Theorem~\ref{t:chi-N-inf-lambda} and remark~(\ref{eq:theta-uniforn-in-eta})
we get the following representation\[
\chi_{N}(+\infty;\lambda)\,\left|_{\left(N,\,\eLH_{N}(\lambda),\, p_{N}(s)\,\right)}\right.=\frac{1}{1+\ggn\etho_{N}(\lambda)/\skb}+\rho_{2,N}(\lambda).\]
Here $\etho_{N}(\lambda)=-2\Re\eLH_{N}(\lambda)$, the real sequence
$\left\{ \ggn\right\} $ is such that $\ggn\sim\frac{1}{2}mN$ as
$N\rightarrow\infty$ and the sequence of functions $\left\{ \rho_{2,N}(\lambda)\right\} $
vanishes uniformly in $\lambda$. Taking into account~(\ref{eq:mN1-m1})
we can rewrite this representation as follows\begin{equation}
\chi_{N}(+\infty;\lambda)\,\left|_{\left(N,\,\eLH_{N}(\lambda),\, p_{N}(s)\,\right)}\right.=\frac{1}{1+\rho_{1,N}\etho_{N}(\lambda)}+\rho_{2,N}(\lambda)\label{eq:decomp-quadr}\end{equation}
where the real sequence $\left\{ \rho_{1,N}\right\} $ is such that
$\rho_{1,N}\sim\frac{1}{2}m_{N,1}N$ as $N\rightarrow\infty$. Using
this result one can study intrinsic scales of the corresponding synchronization
models with large number of components similarly to Theorem~\ref{t:asympt-char-f}.

It would be interesting to know if the decomposition~(\ref{eq:decomp-quadr})
holds for other sequences $\left\{ p_{N}(s)\right\} $ different from~(\ref{eq:pN-p-scale}).

\section{Proofs}

\label{sec:Proofs-step-by-step}

\subsection{Lemmas of dynamics}

\label{sub:Lemmas-dynam}

As in paper~\cite{asmda-2013-man} we start from introducing useful
functions. Fix some even function $\fg=g(a)$ on $\RR^{d}$: \[
\fg:\,\,\RR^{d}\rightarrow\CC,\qquad\fg(a)=\fg(-a).\]
Consider also $\fg_{0}(a)=g(a)-g(0)$. Now define the following functions
on the configuration space~$\RR^{Nd}$\[
V(x):=\frac{2}{(N-1)N}\,\sum_{j_{1}<j_{2}}\fg\left(x_{j_{1}}-x_{j_{2}}\right),\qquad V_{0}(x):=\frac{2}{(N-1)N}\,\sum_{j_{1}<j_{2}}\fg_{0}\left(x_{j_{1}}-x_{j_{2}}\right)\]
where $x=\left(x_{1},\ldots,x_{N}\right)$, $x_{j}\in\RR^{d}$. Evidently,
$V(x)=V_{0}(x)+g(0)$. Note that \[
x_{1}=\cdots=x_{N}\,\quad\Rightarrow\quad V_{0}(x)=0.\]

Keeping in mind notation of Subsections~\ref{sub:perturb-indep}
and~\ref{sub:syncr-ep-T-M} we introduce a map $S_{(i,j)}:$ $\RR^{N}\rightarrow\RR^{N}$,
as follows $S_{(i,j)}x\,:=\, x\circ S_{i,j}$. In other words, \vspace{-3ex}
 \begin{equation}
S_{(i,j)}:\quad(x_{1},\ldots,x_{i},\ldots,x_{j},\ldots,x_{N})\mapsto(x_{1},\ldots,x_{i},\ldots,\begin{array}{c}
\\x_{i}\\
j\end{array},\ldots,x_{N})\,.\label{eq:Sij-map-nm}\end{equation}
Define a map-valued random variable $\mathcal{S}$ such that \begin{equation}
\sP\left\{ \mathcal{S}=S_{(i,j)}\right\} =\frac{1}{(N-1)N},\qquad i\not=j.\label{eq:m-v-rw-S}\end{equation}

\begin{mylem}

\label{lem:V0-Sx-nm}There exists $\kap>0$ such that for any $x\in\RR^{N}$
\[
\sE V_{0}(\mathcal{S}x)=\mkn V_{0}(x),\]
where $\mkn=1-\kap/\left((N-1)N\right).$

\end{mylem} Lemma~\ref{lem:V0-Sx-nm} was proved in~\cite{manita-TVP-large-ident-eng}
for much more general synchronization jumps. For the pair-wise synchronization
interaction considered in the current paper in the framework of the
symmetric model the value of $\varkappa$ is known: $\varkappa=2$.

From this point we take the following concrete even function $\fg(y):=\cos\skpr y{\lambda}$.
Its dependence on the variable $\lambda\in\RR^{d}$ will usually be
omitted. Consider the function \begin{equation}
V(x):=\frac{2}{(N-1)N}\,\sum_{j_{1}<j_{2}}\cos\skpr{\lambda}{x_{j_{1}}-x_{j_{2}}}\label{eq:Vx-def-cos}\end{equation}
corresponding to this choice of $\fg$. It follows form Lemma~\ref{lem:V0-Sx-nm}
that \begin{equation}
\sE V(\mathcal{S}x)=\mkn V(x)+l_{N},\label{eq:E-V-kV-l}\end{equation}
where \begin{equation}
\mkn=1-\frac{\kap}{(N-1)N}\,,\qquad l_{N}:=1-\mkn=\frac{\kap}{(N-1)N}\,.\label{eq:notation-ln-kn-kappa}\end{equation}

\begin{mylem}

\label{lem:V-free-mult-nm}For $s>0$, $x\in\RR^{Nd}$

\begin{equation}
\sE V(x+\xc(s))=V(x)\, e^{-s\etho(\lambda)}\label{eq:E-V-xc-s}\end{equation}
where $\etho(\lambda)=-2\Re\eLH(\lambda)$ and $V$ is defined in~(\ref{eq:Vx-def-cos}). 

\end{mylem}

\noindent \emph{Proof of Lemma~}\ref{lem:V-free-mult-nm}. \[
\cos\skpr{\lambda}y=\frac{\exp\left(i\skpr{\lambda}y\right)+\exp\left(-i\skpr{\lambda}y\right)}{2},\]

\begin{eqnarray*}
\sE\exp\left(i\skpr{\lambda}{x_{j_{1}}+\xc_{j_{1}}(s)-x_{j_{2}}-\xc_{j_{2}}(s)}\right) & = & \exp\left(i\skpr{\lambda}{x_{j_{1}}-x_{j_{2}}}\right)\,\sE\exp\left(i\skpr{\lambda}{\xc_{j_{1}}(s)-\xc_{j_{2}}(s)}\right)\\
 & = & \exp\left(i\skpr{\lambda}{x_{j_{1}}-x_{j_{2}}}\right)\,\chf{j_{1}}s{\lambda}\chf{j_{2}}s{-\lambda}\\
 & = & \exp\left(i\skpr{\lambda}{x_{j_{1}}-x_{j_{2}}}\right)\,\left|\chf{\,}s{\lambda}\right|^{2}\\
\sE\exp\left(-i\skpr{\lambda}{x_{j_{1}}+\xc_{j_{1}}(s)-x_{j_{2}}-\xc_{j_{2}}(s)}\right) & = & \exp\left(-i\skpr{\lambda}{x_{j_{1}}-x_{j_{2}}}\right)\,\left|\chf{\,}s{-\lambda}\right|^{2}\end{eqnarray*}
 Note that $\left|\chf{\,}s{\lambda}\right|^{2}$ is the real symmetric
characteristic function and \[
\left|\chf{\,}s{\lambda}\right|^{2}=\left|\chf{\,}s{-\lambda}\right|^{2}=\left|\exp(s\eLH(\lambda))\right|^{2}=\exp(2\Re\eLH(\lambda)s).\]
So 

\[
\sE\cos\skpr{\lambda}{x_{j_{1}}+\xc_{j_{1}}(s)-x_{j_{2}}-\xc_{j_{2}}(s)}=\cos\skpr{\lambda}{x_{j_{1}}-x_{j_{2}}}e^{-s\etho(\lambda)}.\]
Summing over $j_{1}<j_{2}$ as in~(\ref{eq:Vx-def-cos}) we get (\ref{eq:E-V-xc-s}).
\hfill  $\square$

\smallskip 

The function $V$ defined by (\ref{eq:Vx-def-cos}) is very important
because \begin{equation}
\sE V(x(t))=\chi_{N}(t;\lambda)\label{eq:E-V-x-chi}\end{equation}
where $\chi_{N}(t;\lambda)$ is the characteristic function of $d_{j_{1}j_{2}}^{(N)}(t)=x_{j_{1}}(t)-x_{j_{2}}(t)$
for the symmetric synchronization model $x(t)$ of Subsection~\ref{sec:sym-synch-mod-N}.
Indeed, in symmetric model random variables $d_{jk}^{(N)}(t)$ are
symmetrically distributed hence $\chi_{N}(t;\lambda)$ is real and
\[
\chi_{N}(t;\lambda)=\sE\exp\left(i\skpr{\lambda}{x_{j_{1}}(t)-x_{j_{2}}(t)}\right)=\sE\cos\skpr{\lambda}{x_{j_{1}}(t)-x_{j_{2}}(t)}.\]
Now (\ref{eq:E-V-x-chi}) easily follows from~(\ref{eq:Vx-def-cos}).

\subsection{Recurrent equations}

\label{sub:rec-equats}

Recall that the symmetric $N$-component synchronization model $x(t)$,
$t\geq0$, is the stochastic process with values in~$\mathbb{R}^{Nd}$.
Let $f=f(x)$ be some function on the configuration space~$\mathbb{R}^{Nd}$.
Put \begin{equation}
\emf n=\sE\left(f(x(\tT_{n}+0))\,|\,\left\{ \tT_{q}\right\} _{q=1}^{\infty}\right),\qquad n=1,2,\ldots\,.\label{eq:fn-f-tau-j}\end{equation}
Hence $\emf n$ is a random variable functionally depending on the
sequence $\ptT:=\left\{ \tT_{q}\right\} _{q=1}^{\infty}$. In particular,
we may consider $\left\{ \emV n\right\} $ where $V$ is defined in~(\ref{eq:Vx-def-cos}).
Main result of this subsection will be given in Lemma~\ref{l:rec-eq-Vn-k}
below.

\begin{myremark} 

Note that conditional expectations \[
\sE\left(\cdot\,|\,\left\{ \tau_{l}^{(j)}\right\} _{l=1}^{\infty},\,\jolN\right)\quad\mbox{and}\quad\sE\left(\cdot\,|\,\left\{ \tT_{q}\right\} _{q=1}^{\infty}\right)\]
are different. The first one carries the total information about senders
at epochs $\tT_{q}$ but in the second conditional expectation such
information is unavailable.

\end{myremark} 

Below we will use the telescopic property of the conditional expectation
\[
\sE\left(\sE\left(\cdot\,|\,\xi,\ptT\right)\,|\,\ptT\right)=\sE\left(\cdot\,|\,\ptT\right)\,\]
where $\xi$ is some random variable. Let $V$ be as in~(\ref{eq:Vx-def-cos}).
Then \begin{equation}
\emV n=\sE\left(V(x(\tT_{n}+0))\,|\,\ptT\right)=\sE\left(\sE\left(V(x(\tT_{n}+0))\,|\, x(\tT_{n}),\ptT\right)\,|\,\ptT\right).\,\label{eq:Vn-E-cnd-1}\end{equation}
Consider now $\sE\left(V(x(\tT_{n}+0))\,|\, x(\tT_{n}),\ptT\right)$.
What is the difference between configurations $x(\tT_{n})$ and $x(\tT_{n}+0)$?
This difference is produced by a single message~$(j_{1},j_{2})$
sent from some component~$j_{1}$ to another component~$j_{2}$.
Obviously, the index~$j_{1}$ of the sender is random. What is the
distribution of $j_{1}$? For the symmetric model the answer is simple:
since the dynamics of the stochastic process $x(t)$ is invariant
with respect to permutations of indices the distribution of $j_{1}$
is uniform: \[
\sP\left\{ j_{1}=k\right\} =\frac{1}{N}\,,\qquad k=\olN.\]
In symmetric model the recipient of the message is chosen with probability
$\frac{1}{N-1}$ among the components different from the sender. So
in the symmetric model all messages $(j_{1},j_{2})$ have the same
probability $\frac{1}{(N-1)N}$ to be sent at epoch~$\tT_{n}$. This
means that \[
\sE\left(V(x(\tT_{n}+0))\,|\, x(\tT_{n}),\ptT\right)=\mathsf{E}_{S}V(\mathcal{S}x(\tT_{n}))\]
where averaging $\mathsf{E}_{S}$ is taken over distribution of the
map-valued random variable $\mathcal{S}$ introduced in~(\ref{eq:m-v-rw-S}).
Hence by~(\ref{eq:E-V-kV-l}) we get \begin{equation}
\sE\left(V(x(\tT_{n}+0))\,|\, x(\tT_{n}),\ptT\right)=\mkn V(x(\tT_{n}))+l_{N}.\label{eq:E-V-T-kV-l}\end{equation}

\textbf{}

Consider now \[
\sE\left(V(x(\tT_{n}))\,|\,\ptT\right)=\sE\left(\sE\left(V(x(\tT_{n}))\,|\, x(\tT_{n-1}+0),\ptT\right)\,|\,\ptT\right).\]
There are no synchronization jumps inside the time interval $(\tT_{n-1},\tT_{n})$
hence by Lemma~\ref{lem:V-free-mult-nm} \[
\sE\left(V(x(\tT_{n}))\,|\, x(\tT_{n-1}+0),\ptT\right)=V(x(\tT_{n-1}+0))\,\exp\left(-(\tT_{n}-\tT_{n-1})\etho(\lambda)\right).\]
Applying the conditional averaging $\sE\left(\cdot\,|\,\ptT\right)$
we get \begin{equation}
\sE\left(V(x(\tT_{n}))\,|\,\ptT\right)=\emV{n-1}\exp\left(-(\tT_{n}-\tT_{n-1})\etho(\lambda)\right).\label{eq:EV-n-n-1-free}\end{equation}
Collecting~(\ref{eq:Vn-E-cnd-1}), (\ref{eq:E-V-T-kV-l}) and (\ref{eq:EV-n-n-1-free})
together we obtain

\begin{mylem}  \label{l:rec-eq-Vn-k}\begin{equation}
\emV n=\mkn\emV{n-1}\exp\left(-(\tT_{n}-\tT_{n-1})\etho(\lambda)\right)+l_{N}\,.\label{eq:rec-eq-Vn-k}\end{equation}

\end{mylem}  

On the time interval $(\tT_{\mPiT t},t)$ there are no synchronization
jumps, so similar arguments give \begin{equation}
\sE\left(V(t)\,|\,\ptT\right)=\emV{\mPiT t}\exp\left(-(t-\tT_{\mPiT t})\etho(\lambda)\right).\label{eq:t-Pi-S-t}\end{equation}

\subsection{Representations for the characteristic function}

\label{sub:repr-ch-f}

Recall notation: the point process $\left\{ \tT_{q}\right\} $ is
a superposition of the renewal processes $\left\{ \tau_{l}^{(j)}\right\} $,
$j=1,\ldots,N$, \[
\mPiT t=\sum_{j=1}^{N}\mPi tj=\max\left\{ q\geq0:\,\,\tT_{q}\leq t\,\right\} .\]
 Let $\mkn$ and $\ln$ be the same as in~(\ref{eq:notation-ln-kn-kappa}). 

\begin{mylem}\label{l-repr-chi} For any $t>0$ and $\lambda\in\RR^{d}$\begin{equation}
\chi_{N}(t;\lambda)=\chi_{N}(0;\lambda)\exp(-t\etho(\lambda))\sE\mkn^{\mPiT t}+l_{N}\sE\sum_{q=1}^{\mPiT t}\exp\left(-(t-\tT_{q})\etho(\lambda)\right)\mkn^{\mPiT t-q}\,.\label{eq:chi-N}\end{equation}

\end{mylem}

Similar decompositions were used in \cite{asmda-2013-man} and \cite{manita-questa}. 

\noindent \emph{Proof} \emph{of Lemma}~\ref{l-repr-chi}. Denote
$\sDlt_{q}=\tT_{q}-\tT_{q-1}$. Iterating~(\ref{eq:rec-eq-Vn-k})
we get\[
\emV n=\mkn^{2}\emV{n-2}\exp\left(-(\sDlt_{n-1}+\sDlt_{n})\etho(\lambda)\right)+\mkn\exp\left(-\sDlt_{n}\etho(\lambda)\right)\ln+\ln\,,\]
\begin{eqnarray*}
\emV n & = & \mkn^{n}\emV 0\exp\left(-(\sDlt_{1}+\cdots+\sDlt_{n})\etho(\lambda)\right)+\mkn^{n-1}\exp\left(-(\sDlt_{2}+\cdots+\sDlt_{n})\etho(\lambda)\right)\ln+\cdots+\,\\
 &  & \,+\mkn\exp\left(-\sDlt_{n}\etho(\lambda)\right)\ln+\ln\,.\end{eqnarray*}
Taking into account (\ref{eq:t-Pi-S-t}) and using identity $\sum\limits _{i=q}^{n}\sDlt_{i}=\tT_{n}-\tT_{q-1}$
we come to the following representation\begin{eqnarray*}
\sE\left(V(x(t))\,|\,\,\ptT\right) & = & \emV{\mPiT t}\exp\left(-(t-\tT_{\mPiT t})\etho(\lambda)\right)\,=\,\\
 & = & \mkn^{\mPiT t}\emV 0\exp(-t\etho(\lambda))+\mkn^{\mPiT t-1}\exp\left(-(t-\tT_{1})\etho(\lambda)\right)\ln+\cdots\,\\
 &  & \,+\mkn\exp\left(-(t-\tT_{\mPiT t-1})\etho(\lambda)\right)\ln+\exp\left(-(t-\tT_{\mPiT t})\etho(\lambda)\right)\ln\,.\end{eqnarray*}
The statement of Lemma\emph{~}\ref{l-repr-chi}\emph{ }will now follow
from~(\ref{eq:E-V-x-chi}) if we apply the unconditional expectation
$\sE$ to the both sides of this representation. \hfill $\square$

We introduce some notation. Since for the symmetric model all renewal
processes $\left\{ \tau_{m}^{(k)}\right\} $, $k=1,\ldots,N$, are
equally distributed they have the common renewal function \begin{equation}
H(t)=\sE\mPi tk\,.\label{eq:H-t-E-Pi}\end{equation}
 Similarly, $\left(\mPi tk,\, t\geq0\right)$ have the same moment
generating function \begin{equation}
\prf(u,v)=\sE\left(v^{\mPi uk}\right),\quad u\geq0,\quad v\in\RR.\label{eq:prf-fi-uv}\end{equation}
Denote also $\hvF(s):=1-F(s)$ where $F(s)$ is the common inter-event
probability distribution function of the renewal processes $\left\{ \tau_{m}^{(k)}\right\} $
(see Subsection~\ref{sub:inter-event-symmetr}). If $f_{1}=f_{1}(t)$
and $f_{2}=f_{2}(t)$ are two functions vanishing for $t<0$ then
their convolution $\left(f_{1}\mst f_{2}\right)(t)$ is the function
defined as \[
\left(f_{1}\mst f_{2}\right)(t)=\int_{0}^{t}f_{1}(s)f_{2}(t-s)\, ds\:\]
for $t\geq0$ and $\left(f_{1}\mst f_{2}\right)(t)=0$ for $t<0$.

Define the following functions \begin{eqnarray}
\prfs(u,v) & = & \int_{0}^{s}dH(y)\,\hvF(s+u-y)+\hvF(s+u)+\,\nonumber \\
 &  & \,+v\,\left(g_{s}\mst\prf(\cdot,v)\right)(u)\label{eq:fi1s}\end{eqnarray}
\[
g_{s}(w)=\int_{0}^{s}dH(y)\, p(s+w-y)+p(s+w)\,.\]
Here $s,u,w\geq0$, $v\in\RR$. By definition $\prfs(u,v)=0$ for
$u<0$ and $g_{s}(w)=0$ for $w<0$.

The function $\prfs(u,v)$ has a very clear meaning: it is the moment
generating function for the number renewals in $\left\{ \tau_{m}^{(k)}\right\} $
happened on the interval $[s,s+u]$. Note that the first two summands
in~(\ref{eq:fi1s}) is the probability that the flow $\left\{ \tau_{m}^{(k)}\right\} $
has no renewal on $[s,s+u]$. The probability that the first renewal
in $\left\{ \tau_{m}^{(k)}\right\} $ fits to a small interval $[s+w,\, s+w+dw]$
is equal to $g_{s}(w)\, dw+o(dw)$.

\begin{mylem} \label{l-IN-t-lambda}

The expectation in the second term of (\ref{eq:chi-N}) is \[
I_{N}(t,\lambda):=N\int_{0}^{t}dH(s)\, e^{-(t-s)\etho(\lambda)}\left(\prfs(t-s,\mkn)\right)^{N-1}\prf(t-s,\mkn)\,.\]

\end{mylem}

\emph{Proof} of Lemma~\ref{l-IN-t-lambda}. For any $A\subset\RR_{+}$
denote by $\nEp{(j)}A$ a random variable {}``the number of epochs
of the point process $\left\{ \tau_{l}^{(j)}\right\} _{l=0}^{\infty}$
belonging to~$A$'': \[
\nEp{(j)}A:=\sum_{l=1}^{\infty}\boldsymbol{1}_{\left\{ \tau_{l}^{(j)}\in A\right\} }.\]
Denote also $\nEp SA:=\sum\limits _{j=1}^{N}\nEp{(j)}A$. In particular,
$\mPi tj=\nEp{(j)}[0,t]$, $\mPiT t=\nEp S[0,t]$. Define a function
\begin{equation}
\fa(s)=\left\{ \begin{array}{cc}
\exp(-s\etho(\lambda)), & \quad s\geq0\\
0, & \quad s<0.\end{array}\right.\label{eq:a-s-exp-def}\end{equation}
Then \[
I_{N}(t,\lambda):=\sE\sum_{q=1}^{\mPiT t}\exp\left(-(t-\tT_{q})\etho(\lambda)\right)\mkn^{\mPiT t-q}=\sE\sum_{q=1}^{\infty}\fa(t-\tT_{q})\mkn^{\mPiT t-q}\]
Consider a single summand in these sums \begin{eqnarray*}
\exp\left(-(t-\tT_{q})\etho(\lambda)\right)\mkn^{\mPiT t-q} & = & \fa(t-\tT_{q})\mkn^{\nEp S(\tT_{q},t)}\\
 & = & \fa(t-\tT_{q})\prod_{j=1}^{N}\mkn^{\nEp{(j)}(\tT_{q},t)}\,\\
 & = & \sum_{r=1}^{N}\boldsymbol{1}_{\left\{ \tT_{q}\in\pTau^{(r)}\right\} }\fa(t-\tT_{q})\,\mkn^{\nEp{(r)}(\tT_{q},t)}\prod_{j\not=r}^{N}\mkn^{\nEp{(j)}(\tT_{q},t)}\end{eqnarray*}
So \[
I_{N}(t,\lambda)=\sum_{r=1}^{N}\sE\sum_{q=1}^{\infty}\boldsymbol{1}_{\left\{ \tT_{q}\in\pTau^{(r)}\right\} }\fa(t-\tT_{q})\,\mkn^{\nEp{(r)}(\tT_{q},t)}\prod_{j\not=r}^{N}\mkn^{\nEp{(j)}(\tT_{q},t)}.\]
The point process $\ptT$ is the superposition of the point processes
$\pTau^{(j)},$ $j=\olN$. Hence the sum $\sum\limits _{q=1}^{\infty}\boldsymbol{1}_{\left\{ \tT_{q}\in\pTau^{(r)}\right\} }$
is the summation over all point of $\pTau^{(r)}=\left\{ \tau_{n}^{(r)}\right\} _{n=1}^{\infty}$.
Therefore \begin{eqnarray*}
I_{N}(t,\lambda) & = & \sum_{r=1}^{N}\sE\sum_{n=1}^{\infty}\fa(t-\tau_{n}^{(r)})\,\mkn^{\nEp{(r)}(\tau_{n}^{(r)},t)}\prod_{j\not=r}^{N}\mkn^{\nEp{(j)}(\tau_{n}^{(r)},t)}\\
 & = & N\sE\sum_{n=1}^{\infty}\fa(t-\tau_{n}^{(1)})\,\mkn^{\nEp{(1)}(\tau_{n}^{(1)},t)}\prod_{j=2}^{N}\mkn^{\nEp{(j)}(\tau_{n}^{(1)},t)}\end{eqnarray*}
since in the symmetric model all renewal processes $\pTau^{(j)},$
$j=\olN$, are independent and identically distributed. Note also
that the random variables $\nEp{(1)}(\tau_{n}^{(1)},t)$ and $\nEp{(j)}(\tau_{n}^{(1)},t)$,
$j=\overline{2,N}$, are conditionally independent when the value
of $\tau_{n}^{(1)}$ is known. So we can proceed with our calculation
as follows \begin{eqnarray*}
I_{N}(t,\lambda) & = & N\sE\sum_{n=1}^{\infty}\sE\left(\fa(t-\tau_{n}^{(1)})\,\mkn^{\nEp{(1)}(\tau_{n}^{(1)},t)}\prod_{j=2}^{N}\mkn^{\nEp{(j)}(\tau_{n}^{(1)},t)}\,|\,\tau_{n}^{(1)}\right)\\
 & = & N\sE\sum_{n=1}^{\infty}\fa(t-\tau_{n}^{(1)})\sE\left(\mkn^{\nEp{(1)}(\tau_{n}^{(1)},t)}\,|\,\tau_{n}^{(1)}\right)\prod_{j=2}^{N}\sE\left(\mkn^{\nEp{(j)}(\tau_{n}^{(1)},t)}\,|\,\tau_{n}^{(1)}\right).\end{eqnarray*}
Note that \[
\sE\left(\mkn^{\nEp{(1)}(\tau_{n}^{(1)},t)}\,|\,\tau_{n}^{(1)}\right)=\prf(t-\tau_{n}^{(1)},\mkn)\]
where $\prf$ is the generating function~(\ref{eq:prf-fi-uv}). Here
we have used the fact that there is a renewal at point $\tau_{n}^{(1)}$.
If $j\geq2$ then $\sE\left(\mkn^{\nEp{(j)}(\tau_{n}^{(1)},t)}\,|\,\tau_{n}^{(1)}\right)$
differs from the generating function~$\prf$ because the point process
$\pTau^{(j)}$ has a memory and $\tau_{n}^{(1)}$ is not a renewal
point for $\pTau^{(j)}$. If we denote \begin{equation}
\prfs(u,v):=\sE\left(v^{\nEp{(j)}(s,s+u)}\right)\,.\label{eq:fi1s-def-as-E}\end{equation}
then $\sE\left(\mkn^{\nEp{(j)}(\tau_{n}^{(1)},t)}\,|\,\tau_{n}^{(1)}\right)=\prf_{1,\tau_{n}^{(1)}}(t-\tau_{n}^{(1)},\mkn)$.
Hence the following representation \[
I_{N}(t,\lambda)=N\sE\sum_{n=1}^{\infty}f(\tau_{n}^{(1)})\]
holds with the function \[
f(s):=\,\fa(t-s)\,\prf(t-s,\mkn)\left(\prfs(t-s,\mkn)\right)^{N-1}.\]
It follows from the renewal theory that \[
\sE\sum_{n=1}^{\infty}f(\tau_{n}^{(1)})=\int_{0}^{\infty}f(s)\, dH(s)\]
where $H(s)$ is the renewal function of the point process~$\pTau^{(1)}$
(see~(\ref{eq:H-t-E-Pi})). Recalling notation (\ref{eq:a-s-exp-def})
we conclude that the proof of Lemma~\ref{l-IN-t-lambda} is almost
done. The only thing remains to be proved is the formula~(\ref{eq:fi1s})
for the function defined as~(\ref{eq:fi1s-def-as-E}). This is a
standard exercise form the renewal theory so we leave it to readers.
\hfill  $\square$

The representation~(\ref{eq:chi-N}) will be very useful for Subsection~\ref{sub:around-KRT}.
At this point we discuss the next two immediate corollaries of Lemma~\ref{l-IN-t-lambda}.

Recall that $\etho(\lambda)\geq0$. The first summand in~(\ref{eq:chi-N})
tends to 0 as $t\rightarrow+\infty$ uniformly in~$\lambda\in\RR^{d}$.
Indeed, for any fixed $v\in(0,1)$ the generating function \[
\phi_{S}(t,v)=\sE\left(v^{\mPiT t}\right)\]
tend to 0 as $t\rightarrow+\infty$ since $\mPiT t\rightarrow\infty$
(a.s.)~\cite{Cox}. Hence we come to the following result.

\medskip 

\noindent \textbf{Corollary. }\emph{For any fixed~$N$} \begin{equation}
\lim_{t\rightarrow\infty}\chi_{N}(t;\lambda)=l_{N}\lim_{t\rightarrow\infty}I_{N}(t,\lambda).\label{eq:lim-chi-I-N-lim}\end{equation}

\smallskip 

\noindent The existence of these limits will be proved in Subsection~\ref{sub:around-KRT}.

\smallskip 

The \emph{second corollary} of Lemma~\ref{l-IN-t-lambda} is a short
proof of Theorems~\ref{t-Markov-chi-expl} and \ref{t-Markov-stable}
for the Markovian model. In the Markovian situation the representation~(\ref{eq:chi-N})
turns in a simple explicit formula. Details are given in a separate
subsection.

\subsection{The Markovian case: proofs of Theorems~\ref{t-Markov-chi-expl}
and \ref{t-Markov-stable}}

\label{sub:the-Mark-case}

Assume that (\ref{eq:inter-event-Markov-case}) holds. The inter-event
distribution is exponential and has the ``lack of memory'' property.
Hence the generating functions $\prf(u,v)$ and $\prfi(u,v)$ are
equal and, moreover, \[
\prf(u,v)=\prfi(u,v)=e^{t(v-1)/m},\]
 the probability generating function of the Poisson law with the mean
$t/m$. Since the renewal processes $\mPi tj$ are Poissonian we have
$H(t)=t/m$ (see~\cite{Cox}). Then\begin{eqnarray}
I_{N}(t,\lambda) & = & \frac{N}{m}\int_{0}^{t}ds\, e^{-(t-s)\etho(\lambda)}\left(\prf(t-s,\mkn)\right)^{N}=\nonumber \\
 & = & \frac{N}{m}\int_{0}^{t}ds\, e^{-(t-s)\etho(\lambda)}\exp\left((t-s)(\mkn-1)N/m\right)=\nonumber \\
 & = & \frac{N}{m}\int_{0}^{t}ds\,\exp\left(-A\cdot(t-s)\right),\label{eq:IN-exp-A}\end{eqnarray}
where $A=\etho(\lambda)+l_{N}N/m$. By (\ref{eq:lim-chi-I-N-lim})
we have \begin{eqnarray*}
\chi_{N}(\infty;\lambda) & = & l_{N}I_{N}(\infty,\lambda)=\frac{l_{N}N}{m}\int_{0}^{\infty}du\,\exp\left(-Au\right)\,=\,\\
 & = & \frac{l_{N}N}{m}\,\frac{1}{\etho(\lambda)+l_{N}N/m}=\frac{1}{\etho(\lambda)\cdot m/(l_{N}N)+1}=\\
 & = & \frac{1}{\etho(\lambda)\cdot(N-1)m/\varkappa+1}\,.\end{eqnarray*}
In these calcucations the notation (\ref{eq:notation-ln-kn-kappa})
were used. Theorem~\ref{t-Markov-chi-expl} is proved.

Theorem~\ref{t-Markov-stable} immediately follows from Theorem~\ref{t-Markov-chi-expl}
and the definition of a stable law.

Using~(\ref{eq:IN-exp-A}) we derive an explicit formula for the
characteristic function $\chi_{N}(t;\lambda)$: \[
\chi_{N}(t;\lambda)=\exp\left(-At\right)\chi_{N}(0;\lambda)+\frac{l_{N}N}{m}\int_{0}^{t}ds\,\exp\left(-A\cdot(t-s)\right).\]
It is straightforward to check that $\chi_{N}(t;\lambda)$ satisfies
to the differential equation~(\ref{eq:chi-dif-eq}).

\subsection{The general case: around the Key Renewal Theorem}

\label{sub:around-KRT}

We go back to the general non-Markovian synchronization model. In
this subsection $N$ is fixed. Define functions

\begin{eqnarray}
\prfi(u,v) & = & \frac{1}{m}\int_{0}^{\infty}\,\hvF(w+u)\, dw\,+v\cdot\left(\frac{1}{m}\hvF\mst\prf(\cdot,v)\right)(u),\quad u\geq0,\quad v\in\RR,\label{eq:fi2uv}\end{eqnarray}
and

\begin{equation}
J_{N}(t,\lambda)=N\int_{0}^{t}dH(s)\, e^{-(t-s)\etho(\lambda)}\left(\prfi(t-s,\mkn)\right)^{N-1}\prf(t-s,\mkn)\,.\label{eq:JN-t-l-def}\end{equation}

Roughly speaking, the function~$J_{N}(t,\lambda)$ differs from the
function $I_{N}(t,\lambda)$ by the formal replacing $\prfs(u,v)$
by $\prfi(u,v)$. 

We are going to give a probabilistic interpretation to (\ref{eq:fi2uv}).
Following~\cite[p.~62]{Cox} we denote by $V_{t}$ the forward recurrence-time
in $\mPi tj$, defined as the time measured from $t$ to the next
renewal. It is well known~\cite{Cox} that the law of $V_{t}$ converges
to some absolutely continuous distribution as $t\rightarrow\infty.$
Moreover, the probability density function of the limiting law is
\begin{equation}
p_{V_{\infty}}(w)=\frac{1}{m}\hvF(w),\qquad w\geq0.\label{eq:p-V-inf}\end{equation}
We see that $\prfi(u,v)$ is the generating function for the number
of renewals on $[0,u]$ in the\emph{ modified} renewal process for
which the distribution of the first interval $\Delta_{1}$ is~(\ref{eq:p-V-inf}),
\begin{equation}
\sP\left(\Delta_{1}\leq s\right)=\int_{0}^{s}p_{V_{\infty}}(w)\, dw\,,\label{eq:Delta1-V-inf}\end{equation}
but the intervals $\Delta_{2},\Delta_{3},\ldots\,$ have the same
distribution as before $\sP\left(\Delta_{n}\leq s\right)=F(s)$, $n\geq2$.
Note that the modified renewal process with the first interval $\Delta_{1}$
distributed as (\ref{eq:p-V-inf})--(\ref{eq:Delta1-V-inf}) is a\emph{
stationary} renewal process. We recall also that $\prf(u,v)=\sE v^{\Pi_{u}}$
is the corresponding generating function for the \emph{ordinary} renewal
process.

As it was explained in~(\ref{eq:lim-chi-I-N-lim}) the main task
is to study the limit of the function $I_{N}(t,\lambda)$. Our idea
is to reduce this problem to the analysis of the function~$J_{N}(t,\lambda)$.

\begin{mylem}  \label{l-IN-infty-1}

Let $N$ be fixed and $t\rightarrow\infty$. Then \[
\sup_{\lambda\in\RR^{d}}\left|I_{N}(t,\lambda)-J_{N}(t,\lambda)\right|\rightarrow0\qquad(t\rightarrow\infty).\]

\end{mylem}  

\begin{mylem}  \label{l-IN-infty-2}

Let $N$ be fixed and $t\rightarrow\infty$. The family of functions
$\left\{ J_{N}(t,\lambda),\, t\geq0\right\} $ converges to some limit
$J_{N}(\infty,\lambda)$ as $t\rightarrow\infty$. This limit is uniform
in $\lambda\in\RR^{d}$, \[
\sup_{\lambda\in\RR^{d}}\left|J_{N}(t,\lambda)-J_{N}(\infty,\lambda)\right|\rightarrow0\,,\]
and the limiting function is \begin{equation}
J_{N}(\infty,\lambda)=\frac{N}{m}\int_{0}^{\infty}du\, e^{-u\etho(\lambda)}\left(\prfi(u,\mkn)\right)^{N-1}\prf(u,\mkn)\,.\label{eq:J-infty-f2-u}\end{equation}

\end{mylem}  

Hence using~(\ref{eq:lim-chi-I-N-lim}), Lemma~\ref{l-IN-infty-1}
and Lemma~\ref{l-IN-infty-2} we get the following representation
for the limiting characteristic function: \[
\lim_{t\rightarrow\infty}\chi_{N}(t;\lambda)=l_{N}J_{N}(\infty,\lambda).\]

Lemmas~\ref{l-IN-infty-1} and \ref{l-IN-infty-2} will be proven
in Subsection~\ref{sub:class-Kt}. Before proceeding with the proofs
we want to discuss some connection of the above results with the classical
renewal processes theory.

The Smith theorem states that if a function $f:\,\RR_{+}\rightarrow\RR$
satisfies to certain \emph{sufficient conditions} then

\begin{equation}
\int_{0}^{t}dH(s)\, f(t-s)\rightarrow\frac{1}{m}\int_{0}^{\infty}f(t)\, dt\,\quad(t\rightarrow+\infty)\,\label{eq:uzl-th-ma-vosst}\end{equation}
Here $H(s)$ is the renewal function of the ordinary renewal process
and $m$ is the expectation of the inter-event interval. This statement
is known also as the Key Renewal Theorem (KRT).

Lemma~\ref{l-IN-infty-2} looks like a formal application of the
KRT to the function \begin{equation}
f(t)=e^{-t\etho(\lambda)}\left(\prfi(t,\mkn)\right)^{N-1}\prf(t,\mkn).\label{eq:f-t-exp-fi2-fi}\end{equation}

It is well known that the KRT holds if, for example, any of the next
\emph{sufficient conditions} SC1 or SC2 is satisfied.
\begin{description}
\item [{SC1:}] the function $f(t)$ in (\ref{eq:uzl-th-ma-vosst}) is nonnegative,
nonincreasing and integrable (the Smith's conditions, \cite{Gned-Kovalenko});
\item [{SC2:}] the function $f(t)$ is directly integrable (the Feller's
condition, \cite{Feller-I-II}). 
\end{description}
They should be verified for any fixed~$\lambda$ and $N$. The function~defined
by (\ref{eq:f-t-exp-fi2-fi}) is not convenient to check SC2. 

Consider now SC1. The first condition is evidently satisfied. Since
$\mkn\in(0,1)$ the both generating functions $\prf(t,\mkn)$ and
$\prfi(t,\mkn)$ are nonincreasing in $t$ and the second condition
of SC1 is true. To have integrability of $f(t)$ it is sufficient
to assume that \begin{equation}
\int_{0}^{\infty}\prf(t,\mkn)\, dt\,<\,+\infty.\label{eq:fi-int-finite}\end{equation}
So we are interested in conditions on the inter-event distribution
(on the function $F(s)$ and $p(s)$) that ensure (\ref{eq:fi-int-finite}).
To find such conditions one need to study behavior of $\prf(t,\mkn)$
when $t\rightarrow+\infty$. Similar problems arise in the renewal
theory~\cite{Cox}. It is natural to attack them by using the classical
analytic methods involving the Laplace transform or Tauber theorems.
Based on the experience existing in this field we can imagine that
it would be rather hard to get exhaustive general description of such
distributions in simple and concise terms. In the present paper we
would like to avoid too heavy analytical considerations. From the
other side there is a hope that for many concrete inter-event distributions
the condition (\ref{eq:fi-int-finite}) could be checked by direct
methods. So we chose a {}``happy medium'' and adopt the strategy
followed in the classical book~\cite{Cox}. We consider distributions
with\emph{ rational Laplace transforms} (ME distributions) which are
sufficient for most applications and very convenient in the context
of the current study. In Subsection~\ref{sub:class-Kt} we construct
special classes of functions (we call them the $\mathcal{K}$-classes)
and propose a method based on a set of rules for manipulation of these
functions. Functions of the form~(\ref{eq:f-t-exp-fi2-fi}) belong
to these classes. Moreover, this technique is also very efficient
for proving Lemma~\ref{l-IN-infty-1}. Proofs can be obtained in
a transparent {}``algebraic'' way.

Briefly speeking, we will prove here~(\ref{eq:uzl-th-ma-vosst})
under assumptions different from SC1 and SC2.

\begin{mylem} 

Assume that the function $f(t)$, $t\geq0$, is such that its Laplace
transform \begin{equation}
f^{*}(z)=\int_{0}^{\infty}f(t)\exp\left(-zt\right)\, dt\,\label{eq:Lapl-transf-def}\end{equation}
 is a RPFN-function. Then (\ref{eq:uzl-th-ma-vosst}) holds. 

\end{mylem}  

To conclude this discussion one should mention the paper \cite{Lam-Lehoczky}.
It contains an interesting approach based on the derivation of an
analog of the KRT for the superposition of renewal processes. Nevertheless,
we cannot use results of \cite{Lam-Lehoczky} because they exploite
sufficient conditions similar to the direct intergrability (SC2) which
is very hard to verify.

\subsection{Algebra of functions $\Cl t$ }

\label{sub:class-Kt}

We introduce some notation. Let $\Cl t$ be a linear space of functions
$f=f(t)$, $f:\,\RR_{+}\rightarrow\CC$, having the following form
\[
f(t)=\sum_{j}P_{n_{j}}(t)e^{\lambda_{j}t}\,,\]
where the sum is taken over a finite set of indices, $\lambda_{j}\in\CC$
are such that $\Re\lambda_{j}<0$, $P_{n_{j}}(t)$ are polynomials
with complex coefficients, $n_{j}$ is a degree of the polynomial~$P_{n_{j}}(t)$.
It is important to note that if $f\in\Cl t$ then $f(t)\rightarrow0$
as $t\rightarrow+\infty$.

It is easy to see that the Laplace transform~(\ref{eq:Lapl-transf-def})
maps the set $\Cl t$ to a set $\zCl z$ of complex-valued functions
$f^{*}=f^{*}(z)$, $z\in\CC$, which is exactly the set of RPFN-function.
In other words, the Laplace transform provides a one-to-one correspondence
between the sets $\Cl t$ and $\zCl z$. Other properties of these
sets are listed below.
\begin{description}
\item [{$\prop$1.}] The set $\Cl t$ is an algebra over the field~$\CC$
with the usual operations {}``$+$'' and {}``$\cdot$'', the summation
and the pointwise multiplication of functions. In particular, if functions
$f_{1}=f_{1}(t)$ and $f_{2}=f_{2}(t)$ belong to $\Cl t$ then the
functions $c_{1}f_{1}(t)+c_{2}f_{2}(t)$ and $f_{1}(t)f_{2}(t)$ also
belong to $\Cl t$ for all $c_{1},c_{2}\in\CC$.
\item [{$\prop$2.}] Similarly, the set $\zCl z$ is also an algebra over
$\CC$ with the operations {}``$+$'' and {}``$\cdot$''.
\end{description}
\begin{myremark}The Laplace transform is a one-to-one correspondence
between the vector spaces $\Cl t$ and $\zCl z$ but it is not an
homomorphism of the algebras $\Cl t$ and $\zCl z$. \end{myremark}
\begin{description}
\item [{$\prop$3.}] The set $\Cl t$ is closed with respect to the convolution,
i.e., \[
f_{1},f_{2}\in\Cl t\quad\Rightarrow\quad f_{1}\mst f_{2}\in\Cl t\,,\]
where \[
\left(f_{1}\mst f_{2}\right)(t)=\int_{0}^{t}f_{1}(s)f_{2}(t-s)\, ds\,.\]

\item [{$\prop$4.}] If $f\in\Cl t$ then \[
\int_{0}^{t}dH(s)\, f(t-s)\,=\,\frac{1}{m}\int_{0}^{\infty}f(s)\, ds+\kg_{1}(t),\]
where $\kg_{1}$ is some function from $\Cl t$. In particular, $\kg_{1}(t)\rightarrow0$
as $t\rightarrow+\infty$.
\end{description}
\emph{Proof of} $\prop$4. Denote $q(t)=\int_{0}^{t}dH(s)\, f(t-s)$.
Under Assumptions~P1 there exists a renewal density function $h(s)$
corresponding to $H(s)$: $dH(s)=h(s)\, ds$. It follows from the
classic results \cite[p.~54]{Cox} that \begin{equation}
q^{*}(z)=\left(\int_{0}^{t}h(s)f(t-s)\, ds\right)^{*}(z)=\frac{p^{*}(z)}{1-p^{*}(z)}\, f^{*}(z)\,.\label{eq:lapl-H-f}\end{equation}
By Assumption~2 $p^{*}(z)$ is a RPFN-function. By Lemma~\ref{l--roots-1-p-0}
the equation $p^{*}(z)-1=0$ has a simple root at $z=0$. Hence the
r.h.s of (\ref{eq:lapl-H-f}) has $z_{0}=0$ as a simple pole. Again
by Lemma~\ref{l--roots-1-p-0} there is no other singularities in
the half-plane $\Re z\geq0$ and there is a finite number of other
poles in $\Re z<0$. A pole $z_{j}$ of the order $n_{j}$ of the
Laplace transform $q^{*}(z)$ corresponds to a summand $P_{j}(t)\exp\left(z_{j}t\right)$
in the original function~$q(t)$ where $P_{j}(t)$ is some polynomial
of degree $n_{j}$. Hence \[
q(t)=q_{0}+\sum_{j\not=0}P_{j}(t)\exp\left(z_{j}t\right)\]
where $q_{0}$ is the residual $\ds\res_{z=0}q^{*}(z)$. Recall that
$p^{*}(0)=1$ and $\left(p^{*}\right)'(0)=-m<0$ . Taking the limit
of $q(t)$ as $t\rightarrow+\infty$ we get \begin{eqnarray*}
\lim_{t\rightarrow+\infty}\,\int_{0}^{t}dH(s)\, f(t-s) & = & q(+\infty)=q_{0}\\
 & = & \res_{z=0}\frac{p^{*}(z)}{1-p^{*}(z)}\, f^{*}(z)=\frac{1}{m}f^{*}(0)=\frac{1}{m}\int_{0}^{\infty}f(t)\, dt\,.\qquad\qquad\square\end{eqnarray*}

\medskip 
\begin{description}
\item [{$\prop$5.}] If $f\in\Cl t$, $n\in\mathbb{Z}_{+},$ $\beta\in\CC$
and $\Re\beta>0$ then \[
\kg_{2}(t)=\int_{0}^{t}dH(s)\, s^{n}e^{-\beta s}f(t-s)\,\in\Cl t\,.\]
In particular, $\kg_{2}(t)\rightarrow0$ as $t\rightarrow+\infty$.
\end{description}
\emph{Proof of} $\prop$5. Consider the Laplace transform $\kg_{2}^{*}(z)$:
\[
\left(\int_{0}^{t}dH(s)\, s^{n}e^{-\beta s}f(t-s)\right)^{*}(z)=f^{*}(z)\,(-1)^{m}\frac{d^{m}}{d\beta^{m}}\left(\frac{p^{*}(z+\beta)}{1-p^{*}(z+\beta)}\right)\,.\]
All poles of the function in the r.h.s.\  belong to the left half-plane
$\Re z<0$. \hfill $\square$

\medskip 

Denote by $\tCl sw$ a class of functions $a=a(s,w)$, $s,w\in\RR_{+}$,
of the following form \begin{equation}
a(s,w)=\sum_{k,l}a_{k,l}e_{k}(s)f_{l}(w),\qquad e_{k},f_{l}\in\Cl{\,}\label{eq:aef}\end{equation}
where the sum is taken over a finite set of indices, $a_{k,l}\in\CC$.
In other words, the set $\tCl sw$ is a tensor product of $\Cl s$
and $\Cl w$: \[
\tCl sw=\Cl s\otimes\Cl w\,.\]

\begin{myremark} \label{rem-f-Kt} Note that for any $f\in\Cl t$
the function $f(s+w)$ belongs to~$\tCl sw$. \end{myremark}
\begin{description}
\item [{$\prop$6.}] The set of functions $\tCl sw$ is an algebra, in
particular, $a(s,w)b(s,w)\in\tCl sw$ if $a\in\tCl sw$ and $b\in\tCl sw$.
\item [{$\prop$7.}] If $f=f(t)\in\Cl t$ and $a=a(s,w)\in\tCl sw$ then 

\begin{description}
\item [{a)}] $\ds\quad a(s,w)f(w)\in\tCl sw\,,\,$
\item [{b)}] $\ds\quad\left(a(s,\cdot)\mst f(\cdot)\right)(w)=\int_{0}^{w}a(s,w-y)f(y)\, dy\,\in\tCl sw\,,$
\item [{c)}] $\ds\quad\int_{0}^{+\infty}a(s,w)\, ds\,\in\Cl w\,,\quad\int_{0}^{+\infty}f(s+w)\, ds\,\in\Cl w\,.$
\end{description}
\end{description}
The item $\prop$7b follows from $\prop$3 and (\ref{eq:aef}). The
next two properties are corollaries of $\prop$4 and $\prop$5. Assume
that $a(s,w)\in\tCl sw$ and $H(y)$ is the renewal function~(\ref{eq:H-t-E-Pi}). 
\begin{description}
\item [{$\prop$8.}] $\ds\int_{0}^{s}dH(y)\, a(s-y,w)\,=\,\frac{1}{m}\int_{0}^{\infty}a(s,w)\, ds+\kg_{3}(s,w)$
\\
where $\kg_{3}(s,w)\in\tCl sw$. 
\item [{$\prop$9.}] $\ds\int_{0}^{t}dH(s)\, a(s,t-s)\,\in\Cl t$~. 
\item [{$\prop$10.}] ~~ 

\begin{description}
\item [{a)}] $\forall\, f(t)\in\Cl t$ $\exists\, g(t)\in\Cl t\,:$ $\,\left|f(t)\right|\leq g(t)$
\item [{b)}] $\forall\, a(s,w)\in\tCl sw\,$ $\exists\, b(s,w)\in\tCl sw:$
$\,\left|a(s,w)\right|\leq b(s,w)$~.
\end{description}
\end{description}
The item $\prop$10a can be proved using the following simple bounds:
\begin{eqnarray*}
\left|t^{2k-1}\right| & \leq & 1+t^{2k},\quad k\in\mathbb{N},\\
\left|e^{-\mu t}\right| & \leq & e^{-t\Re\mu}.\end{eqnarray*}
A proof of the item $\prop$10b is similar. \hfill  $\square$

\medskip 

\noindent \emph{Proof of Lemma~\ref{l-IN-infty-1}. }The idea is
to prove that $I_{N}(t,\lambda)=J_{N}(t,\lambda)+\psi_{N}(t,\lambda)\,$
where $\left|\psi_{N}(t,\lambda)\right|\leq\psi_{1,N}(t)$ for some
function $\psi_{1,N}(t)$ from $\Cl t$. This is easy to do by using
the above properties $\prop$1--$\prop$10. Below we give the chain
of conclusions with minor comments.

\begin{doublespace}
To analyze $I_{N}(t,\lambda)$ consider first the function $\prfs(u,v)$
defined by the formula~(\ref{eq:fi1s}).

$\overline{F}(t)\in\Cl t$

$\overline{F}(s+u)\in\tCl su$ (Remark \ref{rem-f-Kt})

Below we use notation $\kg_{4}(s,u)$, $\ldots$, $\kg_{9}(s,u)$
for functions belonging to $\tCl su$. 

$\ds\int_{0}^{s}dH(y)\,\hvF(s-y+u)=\frac{1}{m}\int_{0}^{\infty}\hvF(y+u)\, dy+\kg_{4}(s,u)$,
$\quad\kg_{4}(s,u)\in\tCl su$ (see $\prop$8)

$p(s+w)\in\tCl su$(Remark $\prop$9.)

$\ds g_{s}(w)=\frac{1}{m}\int_{0}^{\infty}p(y+w)\, dy+\kg_{5}(s,u)$,
$\quad\kg_{5}(s,u)\in\tCl su$ (see $\prop$8)

For any fixed $v\in(0,1)$ we have $\prf^{*}(z,v)\in\zCl z$ by the
formula~(\ref{eq:fi-zv-classic}). Hence $\prf(u,v)\in\Cl u$ for
any fixed $v\in(0,1)$. It follows from $\prop$7 and $\prop$8

\[
\left(g_{s}\mst\prf(\cdot,v)\right)(u)=\ds\frac{1}{m}\int_{0}^{\infty}\left(p(y+\cdot)\mst\prf(\cdot,v)\right)(u)\, dy+\kg_{6}(s,u)\,.\]
Note that the function $\kg_{6}(s,u)$ depends on the variable~$v$
but in the current lemma its value is fixed ($v=\mkn$) so we skip
this dependence in the notation~$\kg_{6}(s,u)$. So we get \begin{eqnarray*}
\prfs(u,\hn) & = & \frac{1}{m}\int_{0}^{\infty}\hvF(y+u)\, dy+\hn\,\frac{1}{m}\int_{0}^{\infty}\left(p(y+\cdot)\mst\prf(\cdot,\mkn)\right)(u)\, dy+\kg_{7}(s,u)=\\
 & = & \prfi(u,\hn)+\kg_{7}(s,u)\end{eqnarray*}
where $\prfs$ and $\prfi$ are defined in~(\ref{eq:fi1s}) and~(\ref{eq:fi2uv}).

$\ds\int_{0}^{\infty}\hvF(y+u)\, dy\,$, $\ds\int_{0}^{\infty}p(y+u)\, dy$
$\in\Cl u$ (by $\prop$7c). Hence $\prfi(u,\hn)\in\Cl u$ ($\prop$3).

$\left(\prfs(u,\hn)\right)^{N-1}=\left(\prfi(u,\hn)\right)^{N-1}+\kg_{8}(s,u)$
(by $\prop$7a)\\
Since $\prf(u,\mkn)\in\Cl u$ by $\prop$7a we get

$\left(\prfs(u,\hn)\right)^{N-1}\prf(u,\mkn)=\left(\prfi(u,\hn)\right)^{N-1}\prf(u,\mkn)+\kg_{9}(s,u)$.
\end{doublespace}

Hence \begin{equation}
I_{N}(t)=J_{N}(t)+N\int_{0}^{t}dH(s)\, e^{-(t-s)\etho(\lambda)}\kg_{9}(s,t-s)\,.\label{eq:In-Jn-int-gam}\end{equation}
Denoting $\psi_{N}(t,\lambda)=N\int_{0}^{t}dH(s)\, e^{-(t-s)\etho(\lambda)}\kg_{9}(s,t-s)$
and recalling that $\etho(\lambda)\geq0$ we have\begin{eqnarray*}
\left|\psi_{N}(t,\lambda)\right| & \leq & N\int_{0}^{t}dH(s)\,\left|\kg_{9}(s,t-s)\right|\,\,\\
 & \stackrel{\prop10b}{\leq} & N\int_{0}^{t}dH(s)\,\kg_{10}(s,t-s)\,\in\Cl t\quad(\mbox{by \ensuremath{K}9) .}\end{eqnarray*}
Hence for any fixed~$N$ \[
\sup_{\lambda\in\RR^{d}}\left|\psi_{N}(t,\lambda)\right|\rightarrow0\qquad(t\rightarrow\infty).\qquad\qquad\square\]

\noindent \emph{Proof of Lemma~\ref{l-IN-infty-2}. }In proving Lemma~\ref{l-IN-infty-1}
we obtained inclusions: $\prf(u,v)\in\Cl u$ and $\prfi(u,\hn)\in\Cl u$.
Since $\etho(\lambda)\geq\etho(0)=0$ we have that for any fixed $\lambda\in\RR^{d}$
and $N$ the function \[
f_{\lambda,N}(u):=e^{-u\etho(\lambda)}\left(\prfi(u,\hn)\right)^{N-1}\prf(u,\mkn),\quad\quad u\geq0,\]
belongs to the class $\Cl u$. Recalling the definition~(\ref{eq:JN-t-l-def})
and using~$\prop$4 we conclude that for fixed $\lambda\in\RR^{d}$
and $N$ \begin{equation}
J_{N}(t,\lambda)\rightarrow J_{N}(\infty,\lambda)\quad\quad(t\rightarrow\infty)\label{eq:J-J-t-inf}\end{equation}
and the limit $J_{N}(\infty,\lambda)$ is given by the formula~(\ref{eq:J-infty-f2-u}).
Let us show that this convergence in uniform in~$\lambda$.

\noindent 

We use the following properties of the renewal density function $h(s)$:
\[
h(t)\rightarrow m^{-1},\quad(t\rightarrow\infty)\qquad M:=\sup_{t\geq0}h(t)\,<\,+\infty.\]
Their proof is similar to the proof of $K$4 (see also \cite[\S~4.4]{Cox}).
Note that the functions $\prf$ and $\prfi$ are non-negative hence
$f_{\lambda,N}(u)\geq0$.  Fix some $A>0$ and consider $t>A$. Then
\begin{eqnarray*}
J_{N}(t,\lambda) & = & \int_{0}^{t}h(t-s)\, f_{\lambda,N}(s)\, ds=\\
 & = & \int_{0}^{A}h(t-s)\, f_{\lambda,N}(s)\, ds+\int_{A}^{t}h(t-s)\, f_{\lambda,N}(s)\, ds.\end{eqnarray*}
The second summand can be bounded uniformly in $\lambda$ as\begin{equation}
M\int_{A}^{\infty}\left(\prfi(u,\hn)\right)^{N-1}\prf(u,\mkn)\, du.\label{eq:int-A-inf-for-J}\end{equation}
The intergrand is $f_{0,N}(u)\in\Cl u$ hence (\ref{eq:int-A-inf-for-J})
goes to $0$ as $A\rightarrow+\infty$. Consider \begin{eqnarray*}
\left|\int_{0}^{A}h(t-s)\, f_{\lambda,N}(s)\, ds-\frac{1}{m}\int_{0}^{A}f_{\lambda,N}(s)\, ds\right| & \leq & \int_{0}^{A}\left|h(t-s)-m^{-1}\right|\, f_{\lambda,N}(s)\, ds\\
 & \leq & \int_{0}^{A}\left|h(t-s)-m^{-1}\right|\, f_{0,N}(s)\, ds.\end{eqnarray*}
By the Lebesque domination theorem the last integral vanishes as $t\rightarrow+\infty$.
Now it is readily seen that the convergence~(\ref{eq:J-J-t-inf})
is uniform in~$\lambda\in\RR^{d}$. Lemma~\ref{l-IN-infty-2} is
proved. \hfill $\square$

\begin{doublespace}
%%%
\end{doublespace}

\subsection{Proofs of Theorems~\ref{t:chi-N-inf-lambda} and~\ref{t:asympt-char-f}}

\label{sub:chi-f-decomp}

In this subsection we study asymptotic behavior of the characteristic
function \[
\chi_{N}(\infty;\lambda)=l_{N}J_{N}(\infty,\lambda),\qquad\lambda\in\RR^{d},\]
when $N$ tends to infinity. We will use the representation~(\ref{eq:J-infty-f2-u}).
We start with detailed considerations of the functions $\prf^{*}$
and $\prfi^{*}$.

\subsubsection{Laplace transforms: decompositions and bounds}

Here we obtain decompositions of the functions $\prf^{*}$ and $\prfi^{*}$.
Let $\prf^{*}(z,\mkn)$ and $\prfi^{*}(z,\mkn)$ be their Laplace
transforms: \[
\prf^{*}(z,\mkn)=\int_{0}^{+\infty}e^{-zu}\prf^{*}(u,\mkn)\, du\,,\qquad\prfi^{*}(z,\mkn)=\int_{0}^{+\infty}e^{-zu}\prfi(u,\mkn)\, du\,.\]
Recall \cite[\S~3.2]{Cox} that \begin{equation}
\prf^{*}(z,\mkn)=\frac{1-p^{*}(z)}{z\left(1-\mkn p^{*}(z)\right)}\,\label{eq:fi-zv-classic}\end{equation}
as the Laplace transform of the generating function for the ordinary
renewal process. Similarly, $\prfi^{*}(z,\mkn)$ is the Laplace transform
of the generating function for the modified renewal process (\ref{eq:Delta1-V-inf}):
\[
\prfi^{*}(z,\mkn)=\frac{1-p_{V_{\infty}}^{*}(z)}{z\left(1-\mkn p^{*}(z)\right)}\,.\]

It follows from basic properties of the Laplace transform \cite[\S~1.3]{Cox}
that \[
p_{V_{\infty}}^{*}(z)=\frac{1-p^{*}(z)}{mz}\,.\]
Hence the Laplace transform of the function~(\ref{eq:fi2uv}) is\begin{equation}
\prfi^{*}(z,\mkn)=\frac{p^{*}(z)-1+mz}{mz^{2}}+\hn\frac{1-p^{*}(z)}{mz}\,\prf^{*}(z,\mkn).\label{eq:fi2-zv-1}\end{equation}
It can be rewritten as \[
\prfi^{*}(z,\mkn)=\prf^{*}(z,\mkn)+\left(1-\mkn\right)\,\frac{p^{*}(z)-1+mz\cdot p^{*}(z)}{mz^{2}\cdot\left(1-\mkn p^{*}(z)\right)}\,.\]
Finally we get \begin{equation}
\prfi^{*}(z,\mkn)=\prf^{*}(z,\mkn)\left(1+\left(1-\mkn\right)\,\mtet(z)\right)\label{eq:fi2-fi-theta}\end{equation}
where \[
\mtet(z)\,=\,\frac{p^{*}(z)-1+mz\cdot p^{*}(z)}{mz\cdot\left(1-p^{*}(z)\right)}\,.\]

\begin{myremark}  \label{rem:Mark-Lapl}

While the Markovian case~(\ref{eq:inter-event-Markov-case}) was
completely discussed in Subsection~\ref{sub:the-Mark-case} is is
interesting to see its exceptionality in the formulae derived for
the general situation. If $p(s)=m^{-1}\exp(-s/m),$ $s\geq0,$ the
density of exponential distribution with the mean $m^{-1}$, then
one can easily check that \[
p^{*}(z)=\frac{1}{1+mz}\,,\qquad\mtet(z)=0\,,\qquad\prf^{*}(z,\mkn)=\prfi^{*}(z,\mkn)=\left(z+\frac{1-\mkn}{m}\right)^{-1}\,.\]
If inter-event intervals have a non-exponential distribution then
$\mtet(z)\not=0$.

\end{myremark} 

In the general case we see from formulae~(\ref{eq:fi-zv-classic})--(\ref{eq:fi2-fi-theta})
that $\prf^{*}(z,\mkn)$ and $\prfi^{*}(z,\mkn)$ are RPF-functions.
Our goal is to obtain representation~(\ref{eq:Laurent-RPF}) for
these functions. First of all we will find their poles.

We use the following notation. $\poles g$ denotes the set of poles
of a rational function~$g=g(z)$ and $\roots g$ denotes the set
of its roots: $\roots g=\left\{ z:\, g(z)=0\right\} $. From~(\ref{eq:fi2-zv-1})
we see that $\prfi^{*}$ has the same singularities as $\prf^{*}$.
Hence \begin{equation}
\poles{\prfi^{*}}\,=\,\poles{\prf^{*}}.\label{eq:poles-fi2-union}\end{equation}

Let all assumptions of Subsection~\ref{sub:inter-event-symmetr}
hold. Recall that $r_{0}=0$ is a simple root of the equation $1-p^{*}(z)=0$.
If \[
\roots{1-p^{*}(z)}=\left\{ 0,\, r_{1},\ldots,r_{q}\right\} \]
denotes the set of different roots of the equation $1-p^{*}(z)=0$
then by Lemma~\ref{l--roots-1-p-0} all numbers $r_{1},\ldots,r_{q}$
belong to the subplane $\Re z<0$. By Assumption~P3 the roots $r_{1},\ldots,r_{q}$
are simple that is $\left(p^{*}\right)'(r_{j})\not=0$.

It is well known that roots of a polynomial depend continuously on
its coefficients~(see, for example, \cite{Kato-book} or \cite[Th.~2.7.1]{Watkins2007}).
The coefficients of the equation $1-\mkn p^{*}(z)=0$ are analytic
in~$\mkn$ in the vicinity of~$1$. Hence for sufficiently large
$N$ the {}``perturbed'' equation $1-\mkn p^{*}(z)=0$ has $q+1$
different roots \begin{equation}
\roots{1-\hn p^{*}(z)}=\left\{ \kapn,\, r_{1}^{(N)},\ldots,r_{q}^{(N)}\right\} .\label{eq:rts-1-kn-pzv}\end{equation}
It follows from the general theory~\cite[Ch.~9, \S~2]{Lax-Linalg}
that the roots~(\ref{eq:rts-1-kn-pzv}) are also simple. Any root
$r_{j}^{(N)}$ is close to the root $r_{j}$ in the following sense
\begin{equation}
r_{j}^{(N)}\rightarrow r_{j}\quad as\quad N\rightarrow\infty.\label{eq:rjkN-rj}\end{equation}
 It is straightforward to check that $\kapn$ is real and, moreover,
\begin{equation}
\kapn=-\frac{\gamn}{m_{1}}+\frac{m_{2}}{2m_{1}^{3}}\gamn^{2}+o(\gamn^{2})\quad\quad(N\rightarrow\infty)\label{eq:kapn-gamn}\end{equation}
where $\gamn=\hn^{-1}-1$, $m_{n}=\sE\prir^{n}=\int x^{n}p(x)\, dx$.
In particular, $\kapn<0$ and $\kapn\rightarrow r_{0}=0$. Hence for
sufficiently large $N$ all roots listed in~(\ref{eq:rts-1-kn-pzv})
belong to the subplane $\Re z<0$. Moreover, the real parts of $r_{1}^{(N)},\ldots,r_{q}^{(N)}$
are separated from~$0$. Namely, for sufficiently large $N$ \begin{equation}
r_{1}^{(N)},\ldots,r_{q}^{(N)}\in\left\{ z:\,\Re z<\frac{1}{2}\max_{j=1,\ldots,q}\Re r_{j}<0\right\} .\label{eq:roots-pert-half-plane}\end{equation}

The representation~(\ref{eq:Laurent-RPF}) for the function $\prf^{*}$
takes the following form\begin{equation}
\prf^{*}(z,\hn)=\frac{c_{0}^{(N)}}{z-\kapn}+\sum_{j=1}^{q}\frac{c_{j}^{(N)}}{z-r_{j}^{(N)}}\,.\label{eq:fz-ryad-Lor}\end{equation}
 Since the function $\prfi^{*}$ has the same poles as the function
$\prf^{*}$ we obtain also \begin{equation}
\prfi^{*}(z,\hn)=\frac{d_{0}^{(N)}}{z-\kapn}+\sum_{j=1}^{q}\frac{d_{j}^{(N)}}{z-r_{j}^{(N)}}\,.\label{eq:f2z-ryad-Lor}\end{equation}
 We need some bounds for the coefficients of these decompositions.

\begin{mylem} \label{l-c0-d0-bounds}

There exist $C_{1}>0$ and $C_{2}>0$ such that for sufficiently large
$N$ \[
\left|c_{0}^{(N)}-1\right|<C_{1}\gamn,\quad\left|c_{j}^{(N)}\right|<C_{1}\gamn,\quad\]
\[
\left|d_{0}^{(N)}-1\right|<C_{2}\gamn^{2},\quad\left|d_{j}^{(N)}\right|<C_{2}\gamn^{2},\quad j=1,\ldots,q\,.\]

\end{mylem}

\noindent \emph{Proof of Lemma~\ref{l-c0-d0-bounds}}. All we need
to prove the lemma is a careful calculation of residuals. Consider~(\ref{eq:fi-zv-classic}).
We have

\[
c_{0}^{(N)}=\res_{z=\kapn}\prf^{*}(z,\mkn)=\frac{1-p^{*}(\kapn)}{-\mkn\left(p^{*}\right)'(\kapn)\,\kapn}\,=\,\frac{1-\hn^{-1}}{-\mkn\left(p^{*}\right)'(\kapn)\,\kapn}\,.\]
Using the Taylor\textquoteright{}s theorem with the Lagrange form
of the remainder we have \begin{eqnarray}
p^{*}(0)-p^{*}(\kapn) & = & \left(p^{*}\right)'(\kapn)\,(-\kapn)+\frac{1}{2!}\left(p^{*}\right)''(\kapn)\,(-\kapn)^{2}+\frac{1}{3!}\left(p^{*}\right)^{(3)}(\kapn)\,(-\kapn)^{3}\nonumber \\
 &  & \,+\frac{1}{4!}\left(p^{*}\right)^{(4)}(\xi_{N})\,(-\kapn)^{4}\label{eq:p0-pkapN}\end{eqnarray}
for some $\xi_{N}\in[\kapn,0]$. Expanding \begin{eqnarray*}
\left(p^{*}\right)''(\kapn) & = & \left(p^{*}\right)''(0)+\left(p^{*}\right)^{(3)}(0)\,\kapn+O\left(\kapn^{2}\right),\\
\left(p^{*}\right)^{(3)}(\kapn) & = & \left(p^{*}\right)^{(3)}(0)+O\left(\kapn\right),\qquad N\rightarrow\infty,\end{eqnarray*}
we get from~(\ref{eq:p0-pkapN}) \[
\left(p^{*}\right)'(\kapn)\,\kapn=\left(\hn^{-1}-1\right)+\frac{1}{2}\left(p^{*}\right)''(0)\kapn^{2}+\left(\frac{1}{2}-\frac{1}{6}\right)\left(p^{*}\right)^{(3)}(0)\kapn^{3}+O(\kapn^{4}).\]
Taking into account that $\left(p^{*}\right)^{(n)}(0)=(-1)^{n}m_{n}$
and $\hn=\left(1+\gamn\right)^{-1}$ we obtain \[
c_{0}^{(N)}=\frac{\left(1+\gamn\right)\gamn}{\gamn+\frac{1}{2}m_{2}\kapn^{2}-\frac{1}{3}m_{3}\kapn^{3}+O(\kapn^{4})}\,.\]
Using (\ref{eq:kapn-gamn}) we come to the expansion \begin{equation}
c_{0}^{(N)}=1+\left(1-\frac{m_{2}}{2m_{1}^{2}}\right)\gamn+\frac{9m_{2}^{2}-6m_{1}^{2}m_{2}-4m_{1}m_{3}}{12m_{1}^{4}}\,\gamn^{2}+O(\gamn^{3}).\label{eq:c0-decomp-gamN}\end{equation}

It is seen from~(\ref{eq:fi2-fi-theta}) that \begin{equation}
d_{0}^{(N)}=c_{0}^{(N)}\left(1+(1-\hn)\theta(\kapn)\right)\label{eq:d0-c0-skob}\end{equation}
It is easy to check that \[
\theta(\kapn)=-1+\frac{{m_{2}}}{2\,{m_{1}}^{2}}+\frac{\left(3\,{m_{2}}^{2}-2\,{m_{1}}\,{m_{3}}\right)\,\kapn}{12\,{m_{1}}^{3}}+\frac{\left(3\,{m_{2}}^{3}-4\,{m_{1}}\,{m_{2}}\,{m_{3}}\right)\,\kapn^{2}}{24\,{m_{1}}^{4}}+O(\kapn^{3}).\]

Taking into account that $(1-\hn)=\gamn-\gamn^{2}+O(\gamn^{3})$ and
using (\ref{eq:kapn-gamn}) we get \[
1+(1-\hn)\theta(\kapn)=1+\left(-1+\frac{{m_{2}}}{2\,{m_{1}}^{2}}\right)\,\gamn+\left(1-\frac{{m_{2}}}{2\,{m_{1}}^{2}}-\frac{{m_{2}}^{2}}{4\,{m_{1}}^{4}}+\frac{{m_{3}}}{6\,{m_{1}}^{3}}\right)\,\gamn^{2}+O(\gamn^{3}).\]
Combining the latter decomposition with (\ref{eq:d0-c0-skob}) we
see that $d_{0}^{(N)}$ does not contain a term proportional to the
first power of $\gamn$: \begin{equation}
d_{0}^{(N)}=1+\frac{\left(3\,{m_{2}}^{2}-2\,{m_{1}}\,{m_{3}}\right)\,\gamn^{2}}{12\,{m_{1}}^{4}}+O(\gamn^{3}).\label{eq:d0-decomp-gamN}\end{equation}

\begin{myremark} 

If $p(x)$ is an exponential p.d.f.\  then coefficients in front
of $\gamn$ and $\gamn^{2}$ in (\ref{eq:c0-decomp-gamN}) and~(\ref{eq:d0-decomp-gamN})
vanish.

\end{myremark}

Let us estimate $c_{j}^{(N)}$ and $d_{j}^{(N)}$, $j=1,\ldots,q$.
They are residuals of the first order poles. Hence \[
c_{j}^{(N)}=\res_{z=r_{j}^{(N)}}\prf^{*}(z,\mkn)\,=\,\frac{1-\hn^{-1}}{-\mkn\left(p^{*}\right)'(r_{j}^{(N)})\, r_{j}^{(N)}}\,\,=\,\frac{\gamn}{\mkn\left(p^{*}\right)'(r_{j}^{(N)})\, r_{j}^{(N)}}\,,\]
\[
d_{j}^{(N)}=\res_{z=r_{j}^{(N)}}\prfi^{*}(z,\mkn)\,=\hn\frac{1-p^{*}(r_{j}^{(N)})}{mr_{j}^{(N)}}\, c_{j}^{(N)}=-\frac{\hn\gamn}{mr_{j}^{(N)}}\, c_{j}^{(N)}.\]

By (\ref{eq:rts-1-kn-pzv})--(\ref{eq:rjkN-rj}) there exists $N_{0}>0$
such that the numbers $r_{j}^{(N)}$ and $\left(p^{*}\right)'(r_{j}^{(N)})$,
$j=1,\ldots,q$, are separated from $0$ uniformly in $N\geq N_{0}$.
So we come to the conclusion that for some $C_{1},C_{2}>0$ \[
\left|c_{j}^{(N)}\right|<C_{1}\gamn,\quad\left|d_{j}^{(N)}\right|<C_{2}\gamn^{2},\quad j=1,\ldots,q\,.\]
Lemma~\ref{l-c0-d0-bounds} is proved. \hfill  $\square$

\smallskip 

Using simple properties of the Laplace transform and the decompositions~(\ref{eq:fz-ryad-Lor})--(\ref{eq:f2z-ryad-Lor})
we come to the following representations of the functions $\prf(t,\mkn)$
and $\prfi(t,\mkn)$: 

\[
\prf(t,\mkn)=c_{0}^{(N)}\exp\left(\kapn t\right)+\sum_{j=1}^{q}c_{j}^{(N)}\exp\left(r_{j}^{(N)}t\right)\,,\]
\[
\prfi(t,\mkn)=d_{0}^{(N)}\exp\left(\kapn t\right)+\sum_{j=1}^{q}d_{j}^{(N)}\exp\left(r_{j}^{(N)}t\right)\,.\]
Combining the above formulae with (\ref{eq:kapn-gamn})--(\ref{eq:roots-pert-half-plane})
and Lemma~\ref{l-c0-d0-bounds} we get the next lemma. 

\begin{mylem}\textbf{\label{l-fi-fi2-princ-exp} }There exist $N_{0}\in\mathbb{N},$
$\delta>0$ and $C>0$ such that for all $t\geq0$

\[
\sup_{N\geq N_{0}}\left|\prf(t,\mkn)-c_{0}^{(N)}\exp\left(-\left|\kapn\right|t\right)\right|<C\gamn\exp(-\delta t)\]
\[
\sup_{N\geq N_{0}}\left|\prfi(t,\mkn)-d_{0}^{(N)}\exp\left(-\left|\kapn\right|t\right)\right|<C\gamn^{2}\exp(-\delta t)\]

\end{mylem}

This lemma is very essential for the further proof. Moreover, in order
to prove our main results under assumptions weaker than P2 and P3
one should first derive Lemma~\ref{l-fi-fi2-princ-exp} under that
new assumptions.

\subsubsection{Asymptotics for large $N$}

\label{sub:chi-f-large-N}

From Lemma~\ref{l-IN-infty-2} we know that \begin{equation}
\chi_{N}(\infty;\lambda)=l_{N}J_{N}(\infty,\lambda)\label{eq:chi-N-lN-JN-lambda}\end{equation}
where \begin{equation}
J_{N}(\infty,\lambda)=\frac{N}{m}\int_{0}^{\infty}du\, e^{-u\etho(\lambda)}\left(\prfi(u,\mkn)\right)^{N-1}\prf(u,\mkn).\label{eq:JN-inf-lambda-integr}\end{equation}
Now we will study the limit of $\chi_{N}(\infty;\lambda)$ as $N\rightarrow\infty$.
Recall that $l_{N}\sim c/N^{2}$. Main idea is to show that \begin{equation}
\sup_{\lambda\in\RR^{d}}\left|l_{N}J_{N}(\infty,\lambda)-l_{N}J_{N}^{\circ}(\lambda)\right|\rightarrow0\quad as\quad N\rightarrow\infty\label{eq:lJ-LJ}\end{equation}
where \[
J_{N}^{\circ}(\lambda)=\frac{N}{m}\int_{0}^{\infty}du\, e^{-u\etho(\lambda)}\left(d_{0}^{(N)}\right)^{N-1}c_{0}^{(N)}\exp\left(-N\left|\kapn\right|u\right).\]
In other words $J_{N}^{\circ}(\lambda)$ is obtained from $J_{N}(\infty,\lambda)$
by the formal replacement of the functions $\prfi(u,\mkn)$ and $\prf(u,\mkn)$
by their principal asymptotics (see Lemma~\ref{l-fi-fi2-princ-exp}).
To prove (\ref{eq:lJ-LJ}) we will use the following bounds\[
\left|\prf(t,\mkn)\right|\leq\left(1+C_{3}\gamn\right)\exp\left(-\left|\kapn\right|t\right)\]

\begin{equation}
\left|\prfi(t,\mkn)\right|\leq\left(1+C_{3}\gamn^{2}\right)\exp\left(-\left|\kapn\right|t\right)\label{eq:fi2-C3-exp}\end{equation}
where $C_{3}=C+\max\left(C_{1},C_{2}\right)$ and $N\geq N_{0}$.
Then \[
A_{N}(\lambda):=\left|J_{N}(\infty,\lambda)-\frac{N}{m}\int_{0}^{\infty}du\, e^{-u\etho(\lambda)}\left(\prfi(u,\mkn)\right)^{N-1}c_{0}^{(N)}\exp\left(-\left|\kapn\right|u\right)\right|\leq\]
\begin{eqnarray*}
\, & \leq & \,\frac{N}{m}\int_{0}^{\infty}du\,\left(\left(1+C_{3}\gamn^{2}\right)\exp\left(-\left|\kapn\right|u\right)\right)^{N-1}C\gamn\exp(-\delta u)\,\\
 & \leq & \,\frac{N}{m}\left(1+C_{3}\gamn^{2}\right)^{N-1}\left(C\gamn\right)\frac{1}{\delta+(N-1)\left|\kapn\right|}\,.\end{eqnarray*}
Recall that $\gamn=\hn^{-1}-1\sim\varkappa/N^{2}$, therefore $\left(1+C_{3}\gamn^{2}\right)^{N-1}\rightarrow1$
as $N\rightarrow\infty$. Hence \[
A_{N}(\lambda)\,\leq\,\frac{N}{m}\,\delta^{-1}\,,\quad\forall N\geq N_{1},\]
 for some specially chosen $N_{1}\geq N_{0}$. Consider now\[
B_{N}(\lambda):=\left|\frac{N}{m}\int_{0}^{\infty}du\, e^{-u\etho(\lambda)}\left(\prfi(u,\mkn)\right)^{N-1}c_{0}^{(N)}\exp\left(-\left|\kapn\right|u\right)-J_{N}^{\circ}(\lambda)\right|\,.\]
Denote $a_{N}(u)=\prfi(u,\mkn)$ and $b_{N}(u)=d_{0}^{(N)}\exp\left(-\left|\kapn\right|u\right)$.
We have \[
\left(a_{N}(u)\right)^{N-1}-\left(b_{N}(u)\right)^{N-1}=\left(a_{N}(u)-b_{N}(u)\right)\sum_{i=0}^{N-2}\left(a_{N}(u)\right)^{i}\left(b_{N}(u)\right)^{N-1-i}\,.\]
By (\ref{eq:fi2-C3-exp}) and Lemmas~\ref{l-c0-d0-bounds} and~\ref{l-fi-fi2-princ-exp}
the following bounds hold\[
\max\left(\left|a_{N}(u)\right|,\left|b_{N}(u)\right|\right)\leq\left(1+C_{3}\gamn^{2}\right)\exp\left(-\left|\kapn\right|u\right)\]
\[
\sup_{N\geq N_{0}}\left|a_{N}(u)-b_{N}(u)\right|<C\gamn^{2}\exp(-\delta u).\]
Therefore \begin{eqnarray*}
\left|\left(a_{N}(u)\right)^{N-1}-\left(b_{N}(u)\right)^{N-1}\right| & \leq & C\gamn^{2}e^{-\delta u}\,(N-1)\,\left(1+C_{3}\gamn^{2}\right)^{N-1}\exp\left(-(N-1)\left|\kapn\right|u\right)\end{eqnarray*}
So \begin{eqnarray*}
B_{N}(\lambda) & \leq & \frac{N}{m}\,\underbrace{C\gamn^{2}\,(N-1)\,\left(1+C_{3}\gamn^{2}\right)^{N-1}\left(1+C_{3}\gamn\right)}\int_{0}^{\infty}du\,\exp\left(-\delta u-N\left|\kapn\right|u\right).\end{eqnarray*}
Since the underbraced expression vanishes as $N\rightarrow\infty$
we get $B_{N}\leq\frac{N}{m}\,\delta^{-1}$ for sufficiently large~$N\geq N_{2}$.

We see that the following estimate \begin{equation}
\left|l_{N}J_{N}(\infty,\lambda)-l_{N}J_{N}^{\circ}(\lambda)\right|\leq\left(A_{N}(\lambda)+B_{N}(\lambda)\right)l_{N}\leq\,\frac{2}{m\delta}\, Nl_{N}\,\label{eq:Jinf-Jo-N1}\end{equation}
holds for $N\geq\max(N_{1},N_{2})$. Now the statement~(\ref{eq:lJ-LJ})
easily follows because $Nl_{N}\rightarrow0$ as $N\rightarrow\infty$. 

\begin{myremark} 

It is easy to see that we are able to get a bound even better than
(\ref{eq:Jinf-Jo-N1}), namely, $cl_{N}/(m\delta N)$.

\end{myremark}

We just proved that $\chi_{N}(\infty;\lambda)=l_{N}J_{N}^{\circ}(\lambda)\,+\,\tetdn^{\circ}(\lambda)$
for some function $\tetdn^{\circ}(\lambda)$ such that the bound \[
\left|\tetdn^{\circ}(\lambda)\right|\leq\,\frac{2Nl_{N}}{m\delta}\,,\qquad N\geq\max(N_{1},N_{2}),\]
holds for any function $\etho=\etho(\lambda)\geq0$. 

Let us calculate $l_{N}J_{N}^{\circ}(\lambda)$. We have\begin{eqnarray*}
l_{N}J_{N}^{\circ}(\lambda) & = & \left(d_{0}^{(N)}\right)^{N-1}c_{0}^{(N)}\frac{l_{N}N}{m}\int_{0}^{\infty}du\, e^{-u\etho(\lambda)}\exp\left(-N\left|\kapn\right|u\right)\\
 & = & \left(d_{0}^{(N)}\right)^{N-1}c_{0}^{(N)}m^{-1}\,\frac{l_{N}N}{N\left|\kapn\right|+\etho(\lambda)}\,\\
 & = & \frac{1+\teton}{1+\ggn\etho(\lambda)}\end{eqnarray*}
where\[
\ggn:=\left(N\left|\kapn\right|\right)^{-1},\qquad\teton:=\left(d_{0}^{(N)}\right)^{N-1}c_{0}^{(N)}m^{-1}l_{N}/\left|\kapn\right|-1.\]
 As it is seen from the above estimates (Lemma~\ref{l-c0-d0-bounds})\[
\left(d_{0}^{(N)}\right)^{N-1}c_{0}^{(N)}\rightarrow1\qquad(N\rightarrow\infty).\]
Recall that \[
l_{N}=\frac{\varkappa}{(N-1)N}\,,\quad\hn=1-l_{N},\quad\kapn\sim-\frac{\hn^{-1}-1}{m}=-\frac{l_{N}}{\hn m}\,.\]
We see that if $N\rightarrow\infty$ then $\teton\rightarrow0$ and
$\ggn\sim mN/\varkappa.$ Since $\etho(\lambda)\geq0$ we can write
\[
\chi_{N}(\infty;\lambda)=\frac{1}{1+\ggn\etho(\lambda)}\,+\,\tetdn(\lambda)\]
where the function $\tetdn(\lambda)$ is bounded by\[
\left|\tetdn(\lambda)\right|\leq\left|\teton\right|+\left|\tetdn^{\circ}(\lambda)\right|.\]
It is readily seen that $\tetdn(\lambda)$ satisfies to the conditions~(\ref{eq:theta-2-N})
and~(\ref{eq:theta-uniforn-in-eta}). Theorem~\ref{t:chi-N-inf-lambda}
is proved. \hfill  $\square$

Theorem~\ref{t:asympt-char-f} easily follows from Theorem~\ref{t:chi-N-inf-lambda}
and definitions of domains of attraction to a stable law ($\mySect$~\ref{sub:Intr-scales}).

\section{Conclusions}

We presented a wide class of stochastic synchronization systems whose
dynamics was constructed by means of $\Levy$ processes and superposition
of renewal processes. Such systems can be used after minor modification
to build non-Markovian mathematically tractable models for various
applications in parallel computing, wireless networks etc. For the
symmetric $N$-component models we showed the long time synchronization
in the stochastic sense and proved some limit theorems for the synchronized
systems as $N\rightarrow\infty$. It is interesting to note that the
limiting distributions depend on very few parameters (the $\Levy$
exponent~$\etho(\lambda)$ and the mean~$m$ of an inter-event interval
for a single component). This suggests that Theorems~\ref{t:chi-N-inf-lambda}--\ref{t:asympt-char-f}
hold true under more general assumptions.

Future research could be directed at realistic non-Markovian synchronization
models generalizing already existing studies of WSNs~\cite{manita-questa,lmanita-controlling}.
Methods of the present paper can also be adapted for studying correlations
between components of synchronization systems.

% end of small 

\end{document}